\documentclass[preprint,aos]{imsart}
\usepackage{amsthm,amsmath, amssymb}
\usepackage{amsfonts,dsfont} 
\usepackage{epsfig}
\psfigdriver{dvips}
\arxiv{0802.0566}
\startlocaldefs

\newcommand{\meil}[1]{\mathbf{#1}} 
\newcommand{\latin}[1]{\emph{#1}}

\theoremstyle{plain}
\newtheorem{theorem}{Theorem}
\newtheorem{lemma}{Lemma}

\newtheorem{proposition}[lemma]{Proposition}
\theoremstyle{definition}

\newtheorem{proc}{Algorithm}
\theoremstyle{remark}
\newtheorem{remark}{Remark}

\newcommand{\egaldef}{:=} 
\newcommand{\flens}{\mapsto} 
\newcommand{\flapp}{\mapsto} 
\newcommand{\telque}{\, \mbox{ s.t. } \,} 

\newcommand{\1}{\mathds{1}} 

\newcommand{\R}{\mathbb{R}} 
\newcommand{\N}{\mathbb{N}} 
\newcommand{\X}{\mathcal{X}}
\newcommand{\Y}{\mathcal{Y}}

\DeclareMathOperator{\card}{Card} 

\newcommand{\mini}[2]{#1 \wedge #2}

\newcommand{\maxi}[2]{#1 \vee #2}

\newcommand{\nupl}{_{1 \ldots n}} 

\newcommand{\paren}[1]{\left( \left. #1 \right. \right)} 
\newcommand{\croch}[1]{\left[ \left. #1 \right. \right]} 
\newcommand{\set}[1]{\left\{ \left. #1 \right. \right\}}
\newcommand{\absj}[1]{\left\lvert #1 \right\rvert} 
\providecommand{\norm}[1]{\left \lVert #1 \right\rVert}
\newcommand{\carre}[1]{\left(#1\right)^2}

\renewcommand{\P}{\mathbb{P}}
\newcommand{\Prob}{\mathbb{P}} 
\DeclareMathOperator{\supp}{supp} 
\newcommand{\E}{\mathbb{E}} 
\newcommand{\sachant}{\, \right| \left. \,} 
\newcommand{\sachantd}{\, \right. \left| \,} 
\newcommand{\loi}{\mathcal{L}} 
\DeclareMathOperator{\Leb}{Leb} 

\newcommand{\bayes}{s}
\newcommand{\perte}[1]{l(\bayes , #1 )}
\newcommand{\ERM}{\widehat{s}}

\newcommand{\M}{\mathcal{M}}
\newcommand{\Mh}{\widehat{\mathcal{M}}} 
\newcommand{\mM}{m \in \M}
\newcommand{\mMh}{m \in \Mh}
\newcommand{\mo}{m^{\star}} 
\newcommand{\mh}{\widehat{m}}
\newcommand{\mhVF}{\mh_{\mathrm{VFCV}}}

\DeclareMathOperator{\pen}{pen}
\DeclareMathOperator{\crit}{crit}
\newcommand{\critVF}{\crit_{\mathrm{VFCV}}} 
\newcommand{\penVFCV}{\pen_{\mathrm{VFCV}}} 
\newcommand{\penVF}{\pen_{\mathrm{VF}}} 
\newcommand{\penid}{\pen_{\mathrm{id}}} 

\newcommand{\ERMb}[1][]{\widehat{s}^{W #1}} 
\newcommand{\Pnb}[1][]{{P_n^{W #1}}}
\newcommand{\Es}{{\E_W}}
\newcommand{\CWinf}{\ensuremath{C_{W,\infty}}}

\newcommand{\hyp}[1]{\ensuremath{\mathbf{(#1)}}} 
\newcommand{\hypAb}{\hyp{Ab}} 
\newcommand{\hypAn}{\hyp{An}} 
\newcommand{\hypAp}{\hyp{Ap}} 
\newcommand{\hypArXl}{\hyp{Ar^{X}_{\ell}}} 
\newcommand{\hypPpoly}{\hyp{P1}} 
\newcommand{\hypPrich}{\hyp{P2}} 
\newcommand{\hypPconst}{\hyp{P3}} 
\newcommand{\hypPtoutVF}{\hyp{P1-3}} 
\newcommand{\hypVF}{\hyp{pVF}} 

\newcommand{\crXl}{\ensuremath{c_{\mathrm{r},\ell}^X}}
\newcommand{\cbiasmaj}{C_{\mathrm{b}}^{+}} 
\newcommand{\cbiasmin}{C_{\mathrm{b}}^{-}} 
\newcommand{\aM}{\alpha_{\M}}
\newcommand{\cM}{c_{\M}}
\newcommand{\sigmin}{\sigma_{\min}} 

\newcommand{\einv}[1]{e^+_{#1}}
\newcommand{\einvz}[1]{e^0_{#1}}

\newcommand{\delc}{\overline{\delta}} 
\newcommand{\punmin}{\widetilde{p_1}} 
\newcommand{\punzero}{\punmin^{(0)}} 

\newcommand{\El}{\E^{\Lambda_m}}
\newcommand{\Il}{I_{\lambda}}
\newcommand{\lamm}{\lambda \in \Lambda_m} 

\newcommand{\pl}{p_{\lambda}} 
\newcommand{\betl}{\beta_{\lambda}} 
\newcommand{\sigl}{\sigma_{\lambda}}

\newcommand{\phl}{\widehat{p}_{\lambda}} 
\newcommand{\bethl}{\widehat{\beta}_{\lambda}} 

\newcommand{\phlW}[1][]{\ensuremath{\widehat{p}^{W #1}_{\lambda}}} 
\newcommand{\bethlW}[1][]{\ensuremath{\widehat{\beta}^{W #1}_{\lambda}}} 
\newcommand{\Wl}{W_{\lambda}}

\endlocaldefs

\begin{document}
\begin{frontmatter}

\title{$V$-fold cross-validation improved: $V$-fold penalization}
\runtitle{$V$-fold penalization}

\begin{aug}
\author{\fnms{Sylvain} \snm{Arlot}\ead[label=e1]{sylvain.arlot@math.u-psud.fr}}
\runauthor{Arlot, S.}

\affiliation{Universit\'e Paris-Sud\vspace{1cm}\\}

\address{Sylvain Arlot\\
Univ Paris-Sud, UMR 8628, \\Laboratoire de Math\'ematiques,\\ Orsay, F-91405 ;
CNRS, Orsay, F-91405 ;\\
INRIA-Futurs, Projet Select \\
\printead{e1}\\
\phantom{E-mail: sylvain.arlot@math.u-psud.fr\ }}
\end{aug}

\runauthor{Arlot, S.}

\begin{abstract}
We study the efficiency of $V$-fold cross-validation (VFCV) for model selection from the non-asymptotic viewpoint, and suggest an improvement on it, which we call ``$V$-fold penalization''. 

Considering  a particular (though simple) regression problem, we prove that VFCV with a bounded $V$ is suboptimal for model selection, because it ``overpenalizes'' all the more that $V$ is large. Hence, asymptotic optimality requires $V$ to go to infinity.
However, when the signal-to-noise ratio is low, it appears that overpenalizing is necessary, so that the optimal $V$ is not always the larger one, despite of the variability issue. This is confirmed by some simulated data.

In order to improve on the prediction performance of VFCV, we define a new model selection procedure, called ``$V$-fold penalization'' (penVF). It is a $V$-fold subsampling version of Efron's bootstrap penalties, so that it has the same computational cost as VFCV, while being more flexible. In a heteroscedastic regression framework, assuming the models to have a particular structure, we prove that penVF satisfies a non-asymptotic oracle inequality with a leading constant that tends to 1 when the sample size goes to infinity. In particular, this implies adaptivity to the smoothness of the regression function, even with a highly heteroscedastic noise.
Moreover, it is easy to overpenalize with penVF, independently from the $V$ parameter. A simulation study shows that this results in a significant improvement on VFCV in non-asymptotic situations.
\end{abstract}



\end{frontmatter}

\label{VFCV.sec}

\section{Introduction}
%
%
There are typically two kinds of model selection criteria. On the one-hand, penalized criteria are the sum of an empirical loss and some penalty term, often measuring the complexity of the models. This is the case of AIC (Akaike \cite{Aka:1973}), Mallows' $C_p$ or $C_L$ (Mallows \cite{Mal:1973}) and BIC (Schwarz \cite{Sch:1978}), to name but a few.
On the other hand, cross-validation (Allen \cite{All:1974}, Stone \cite{Sto:1974}, Geisser \cite{Gei:1975}) and related criteria are based on the idea of data splitting. Part of the data (the training set) is used for fitting each model, and the rest of the data (the validation set) is used to measure the performance of the models. There are several versions of cross-validation (CV), \latin{e.g.} leave-one-out (LOO, also called ordinary CV), leave-$p$-out (LPO, also called delete-$p$ CV) and generalized CV (Craven and Wahba \cite{Cra_Wah:1979}).
In practical applications, cross-validation is often computationally very expensive. This is why less greedy CV algorithms have been proposed, among which $V$-fold cross-validation (VFCV, Geisser \cite{Gei:1975}) and repeated learning testing methods (Breiman \latin{et al.} \cite{Bre_etal:1984}). 
In this article, we mainly consider VFCV --- which seems to be the most widely used nowadays --- when the goal of model selection is to be {\em efficient}, \latin{i.e.} to minimize the prediction risk among a family of estimators. Let us emphasize that this is quite different from picking up the ``true model'', which is often recalled as the {\em identification} or {\em consistency} issue.

%
%
The properties of CV (in particular leave-$p$-out) for prediction and model identification have been widely studied from the asymptotical viewpoint. It typically depends on the splitting ratio, \latin{i.e.} the ratio between the sizes of the validation and training sets ($p/(n-p)$ in the leave-$p$-out case; $1/(V-1)$ for $V$-fold cross-validation). This has been shown for instance by Shao \cite{Sha:1997} (for regression on linear models) and by van der Laan, Dudoit and Keles \cite{vdL_Dud_Kel:2004} (for density estimation).
Asymptotic optimality occurs when this ratio goes to zero at infinity, as shown by Li \cite{KCLi:1987} for the leave-one-out, and generalized by Shao \cite{Sha:1997} for the leave-$p$-out with $p \ll n$, both in the regression setting, when all the models are linear. Other asymptotic results about CV in regression can be found in the book by Gy\"orfi \latin{et al.} \cite{Gyo_etal:2002}, and in the paper of van der Laan, Dudoit and Keles \cite{vdL_Dud_Kel:2004} for density estimation. 
Notice that the behaviour of these procedures changes completely when the goal is consistency; we refer to Yang \cite{Yan:2007b} and Sect.~\ref{VFCV.sec.disc.consist} below for references on this problem. 

%
%
When it comes to practical application, a major question is how to choose the tuning parameters of CV procedures, since their performance strongly depend on them. In the case of VFCV, this means choosing $V$. Basically, there are three competing factors. 
First, the VFCV estimator of the prediction error, $\critVF$, is biased, and its bias decreases with $V$. As shown by Burman \cite{Bur:1989,Bur:1990}, it is possible to correct this bias; otherwise, $V$ should not be taken too small. 
Second, the variance of $\critVF$ depends on $V$: it is always decreasing for small values of $V$, but then it can either stay decreasing (as in the linear regression case \cite{Bur:1989}) or start to increase before $V=n$ (as in some classification problems \cite{Bre:1996,Has_Tib_Fri:2001,Mol_Sim_Pfe:2005} or in density estimation \cite{Cel_Rob:2006}; see Sect.~\ref{VFCV.sec.subopt.VF.choix}). 
Third, the computational cost of VFCV is proportional to $V$, so that the theoretic optimum (taking only bias and variability into account) can not always be computed. More precisely, it is necessary to understand well how the performance of VFCV depends on $V$ before taking into account the computational cost. This is one of the purposes of this article.

\medskip

We here aim at providing a better understanding of some CV procedures (including VFCV) from the {\em non-asymptotic viewpoint}. This may have two major implications.
First, non-asymptotic results are made to handle collections of models which may depend on the sample size $n$: their sizes may typically be a power of $n$, and they may contain models whose complexities grow with $n$. Such collections of models are particularly significant for designing adaptive estimators of a function which is only assumed to belong to some h\"olderian ball, which may require an arbitrarily large number of parameters.
Second, in several practical applications, we are in a ``non-asymptotic situation'' in the sense that the signal-to-noise ratio is low. We shall see in the following that it should really be taken into account for an optimal tuning of $V$.
It is worth noticing that such a non-asymptotic approach is not common in the literature, since most of the results already mentioned are asymptotic, and none is considering our second point above. 

Another important point in our approach is that our framework includes several kinds of heteroscedastic data. We only assume that the observations $(X_i,Y_i)_{1 \leq i \leq n}$ are i.i.d. with 
\[ Y_i = \bayes(X_i) + \sigma(X_i) \epsilon_i \enspace , \]
where $\bayes: \X \flens \R$ is the (unknown) regression function, $\sigma: \X \flens \R$ is the (unknown) noise-level, and $\epsilon_i$ has a zero mean and a unit variance conditionally to $X_i$. In particular, the noise-level $\sigma(X)$ can be strongly dependent from $X$, and the distribution of $\epsilon$ can itself depend from $X$. 
Such data are generally considered as very difficult to handle, because we have no information on $\sigma$, making irregularities of the signal harder to distinguish from noise. Then, simple model selection procedures such as Mallows' $C_p$ may not work (see Chap.~4 of \cite{Arl:2007:phd} for a theoretical argument), and it is natural to hope that VFCV or other resampling methods may be robust to heteroscedasticity. In this article, both theoretical and simulation results confirm this fact. 

\medskip

%
%
In Sect.~\ref{VFCV.sec.subopt.VF}, we provide a non-asymptotic analysis of the performance of VFCV. The aforementioned bias turns out into a non-asymptotic negative result (Thm.~\ref{VFCV.the.subopt}), showing a rather simple problem for which VFCV can not satisfy an oracle inequality with leading constant smaller than $\kappa(V)-\epsilon_n$, with $\kappa(V)>1$ for any $V \geq 2$ and $\epsilon_n \rightarrow 0$. In particular, VFCV with a bounded $V$ can not be asymptotically optimal. 
But our analysis also has a major positive consequence in some ``non-asymptotic'' situations. Indeed, by considering VFCV as a penalization procedure, our previous result can be interpretated as an overpenalization property of VFCV. This should be related to the fact that the efficiency of penalization methods (like Mallows' $C_p$) is often improved by overpenalization, when the signal-to-noise ratio is small. Then, one can expect the optimal $V$ for VFCV to be smaller than $n$, even for least-squares regression, which is confirmed by the simulation study of Sect.~\ref{VFCV.sec.simus}. 
So, it appears that choosing the optimal $V$ for VFCV may be quite hard. In addition, the optimal choice may not be satisfactory when it corresponds to a highly variable criterion such as the $2$-fold CV one. It is likely that there is some room left here to improve on VFCV.

\medskip

%
%
This is why we propose in Sect.~\ref{VFCV.sec.penVFCV} another $V$-fold algorithm, that we call ``$V$-fold penalization'' (penVF). It is based upon Efron's resampling heuristics \cite{Efr:1979}, in the same way as Efron's bootstrap penalty \cite{Efr:1983}, but with a $V$-fold subsampling scheme instead of the bootstrap. It thus has exactly the same computational cost as the classical VFCV, and our results show that is has a similar robustness property, in some heteroscedastic regression framework.
In addition, it turns out to be a generalization of Burman's corrected VFCV \cite{Bur:1989,Bur:1990} (at least when the splitting into $V$ blocks is regular). The main advance of penVF being that it is straightforward to overpenalize within any factor when this is required, for instance when the signal-to-noise ratio seems low.

In the least-square regression framework, when we have to select among histogram models (see Sect.~\ref{VFCV.sec.VFCV.calculs} for an accurate definition), we prove that penVF satisfies a non-asymptotic oracle inequality with a leading constant almost one (Thm.~\ref{VFCV.the.histos.penVFCV}). To our knowledge, such a non-asymptotic result is new for any $V$-fold model selection procedure. One of its strengths is that it requires very few assumptions on the noise, allowing in particular heteroscedasticity. It is a strong result for penVF --- which was not built for this particular setting at all --- to improve on VFCV for such difficult problems, where VFCV is among the best procedures overall.
As a consequence of Thm.~\ref{VFCV.the.histos.penVFCV}, one can use penVF with the family of regular histograms in order to obtain an estimator adaptive to the smoothness of the regression function, when the noise is heteroscedastic (while having no information at all on the distribution of the noise).
Notice that we only consider this result as a first step towards a more general theorem, without the restriction to histograms, as discussed in Sect.~\ref{VFCV.sec.disc.pred}. The main interest of this toy framework is that we can study it deeply, and then derive general heuristics for practical use.

\medskip

As an illustration to our theoretical study, we provide the results of a simulation study in Sect.~\ref{VFCV.sec.simus}. It confirms the good performances of penVF against both VFCV and the simpler Mallows' $C_p$ criterion, in particular for difficult heteroscedastic problems. We also show how useful may be the flexibility of penVF when the signal-to-noise ratio is low. By decoupling $V$ from the overpenalization factor, we allowed a significant improvement of the performance of both VFCV and its bias-corrected version.

Finally, our results are discussed in Sect.~\ref{VFCV.sec.disc}. The remaining of the paper is devoted to some probabilistic tools (App.~\ref{VFCV.sec.proba}) and proofs (App.~\ref{VFCV.sec.proofs}).

\section{Performance of $V$-fold cross-validation}
\label{VFCV.sec.subopt.VF}

In this section, we provide a non-asymptotic study of $V$-fold cross-validation (VFCV) in the least-squares regression framework. In order to make explicit computations possible, we focus on the case where each model is an ``histogram model'', \latin{i.e.} the vector space of piecewise constant functions on some fixed partition of the feature space. This is only a first theoretical step. We use it to derive heuristics, that should help the practical user of VFCV in any framework.  Notice also that we do not assume that the regression function itself is piecewise constant.

\subsection{General framework} \label{VFCV.sec.subopt.VF.gal}

First consider the general prediction setting: $\X \times \Y$ is a measurable space, $P$ an unknown probability measure on it and we observe some data $(X_1,Y_1), \ldots, (X_n,Y_n) \in \X \times \Y$ of common law $P$. Let $\mathcal{S}$ be the set of predictors (measurable functions $\X \mapsto \Y$) and $\gamma: \mathcal{S} \times (\X \times \Y) \mapsto \R$ a contrast function. Given a family $(\ERM_m)_{\mM_n}$ of data-dependent predictors, our goal is to find the one minimizing the prediction loss $P\gamma(t) \egaldef \E_{(X,Y)\sim P}[\gamma(t,(X,Y))]$. Notice that the expectation here is only taken w.r.t. $(X,Y)$, so that $P \gamma(t)$ is random when $t$ is random (\latin{e.g.} data-driven). Assuming that there exists a minimizer $\bayes \in \mathcal{S}$ of the loss (the Bayes predictor), we will often consider the excess loss $\perte{t} = P\gamma(t) - P\gamma(\bayes) \geq 0$ instead of the loss.

\medskip

We assume that each predictor $\ERM_m$ can be written as a function $\ERM_m(P_n)$ of the empirical distribution of the data $P_n = n^{-1} \sum_{i=1}^n \delta_{(X_i,Y_i)}$. The case-example of such a predictor is the empirical risk minimizer $\ERM_m \in \arg\min_{t \in S_m} \set{ P_n \gamma(t)}$, where $S_m$ is any set of predictors (called a {\em model}).
In the classical version of VFCV, we first choose some partition $(B_j)_{1\leq j \leq V}$ of the indexes $\set{1, \ldots, n}$. Then, we define 
\begin{gather*}
 P_n^{(j)} = \frac{1}{\card(B_j)} \sum_{i \in B_j} \delta_{(X_i,Y_i)} \qquad \ERM_m^{(j)} = \ERM_m \paren{ P_n^{(j)} } \\
 P_n^{(-j)} = \frac{1}{n-\card(B_j)} \sum_{i \notin B_j} \delta_{(X_i,Y_i)} \qquad \ERM_m^{(-j)} = \ERM_m \paren{ P_n^{(-j)} } \enspace .
\end{gather*}
The final VFCV estimator is $\ERM_{\mhVF} (P_n) $ with 
\begin{equation} \label{VFCV.def.mh.VFCvclass}
\mhVF \in \arg\min_{\mM_n} \left\{  \critVF(m) \right\} \quad \mbox{and} \quad \critVF(m) \egaldef  \frac{1}{V} \sum_{j=1}^V P_n^{(j)} \gamma \paren{\ERM_m^{(-j)} }  \enspace .
\end{equation}

It is classical to assume that the partition $(B_j)_{1\leq j \leq V}$ is regular, \latin{i.e.} that $\forall j$, $\absj{ \card(B_j) -  n/V} < 1$. 
In order to understand deeply the properties of VFCV, we have to compare precisely $\critVF(m)$ to the excess loss $\perte{\ERM_m}$. A crucial point is to compare their expectations, which is quite hard in general. This is why we restrict ourselves to a particular framework, namely the histogram regression one. We describe it in the next subsection.

\subsection{The histogram regression case}
\label{VFCV.sec.VFCV.calculs}
In the regression framework, the data $(X_i,Y_i) \in \X \times \R$ are i.i.d. of common law $P$. Denoting by $\bayes$ the regression function, we have 
\begin{equation} \label{VFCV.eq.donnees.reg} Y_i = \bayes(X_i) + \sigma(X_i) \epsilon_i \end{equation} 
where $\sigma: \X \mapsto \R$ is the heteroscedastic noise-level and $\epsilon_i$ are i.i.d. centered noise terms, possibly dependent from $X_i$, but with mean 0 and variance 1 conditionally to $X_i$.
In order to simplify the theory, we will make two main assumptions on the data throughout this paper:
\[ \sigma(X) \geq \sigmin > 0 \quad \mbox{a.s.} \qquad \mbox{and} \qquad \norm{Y}_{\infty} \leq A < + \infty \enspace . \]
Notice that we do not assume $\sigmin$ and $A$ to be known from the statistician. Moreover, those two assumptions can be relaxed, as shown by Chap.~6 and Sect.~8.3 of \cite{Arl:2007:phd}. %
The feature space $\X$ is typically a compact subset of $\R^d$. We use the least-squares contrast $\gamma: (t,(x,y)) \mapsto (t(x)-y)^2$ to measure the quality of a predictor $t: \X \mapsto \Y$. As a consequence, the Bayes predictor is the regression function $\bayes$, and the excess loss is $\perte{t} = \E_{(X,Y)\sim P} \carre{t(X)-\bayes(X)}$. 
To each model $S_m$, we associate the {\em empirical risk minimizer}
\[ \ERM_m \egaldef \ERM_m(P_n) = \arg\min_{t \in S_m} \set{ P_n \gamma(t) } \] (when it exists and is unique). Define also $\bayes_m \egaldef \arg\min_{t \in S_m} P \gamma\paren{t}$.

We now focus on histograms. Each model in $(S_m)_{\mM_n}$ is the set of piecewise constant functions (histograms) on some partition $(\Il)_{\lamm}$ of $\X$. It is thus a vector space of dimension $D_m = \card(\Lambda_m)$, spanned by the family $(\1_{\Il})_{\lamm}$. As this basis is orthogonal in $L^2(\mu)$ for any probability measure $\mu$ on $\X$, we can make explicit computations. The following notations will be useful throughout this article.
\begin{gather*}
\pl \egaldef P(X\in \Il) \qquad \phl \egaldef P_n(X \in \Il)  \qquad
\sigl^2 \egaldef \E \croch{ \paren{Y - \bayes(X)}^2 \sachant X \in \Il } \enspace .
\end{gather*}
Remark that $\ERM_m$ is uniquely defined if and only if each $\Il$ contains at least one of the $X_i$, \latin{i.e.} $\min_{\lamm} \set{\phl}>0$. 
Prop.~\ref{VFCV.pro.EcritVFCV-EcritID} below compares the $V$-fold criterion and the ideal criterion $P \gamma(\ERM_m)$ in expectation.
\begin{proposition} \label{VFCV.pro.EcritVFCV-EcritID}
Let $S_m$ be the model of histograms associated with the partition $(\Il)_{\lamm}$ and $\paren{B_j}{1 \leq j \leq V}$ some ``almost regular'' partition of $\set{1, \ldots, n}$, \latin{i.e.} such that 
\[ \max_j \set{ \frac{ \card(B_j) } {n} }  \leq c_B < 1 \qquad \mbox{and} \qquad \sup_j \set{ \absj{ \frac{\card(B_j)}{n} - \frac{1}{V}} } \leq \epsilon^{reg}_n \xrightarrow[n \rightarrow \infty]{} 0 
\enspace . \]
Then, the expectation of the ideal and $V$-fold criteria are respectively equal to 
\begin{align} \label{VFCV.eq.EcritID}
\E \croch{ P \gamma(\ERM_m) } &= P \gamma(\bayes_m) + \frac{1}{n} \sum_{\lamm} \paren{1 + \delta_{n, \pl}} \sigl^2 
\\ \label{VFCV.eq.EcritVFCV}
\E \croch{ \critVF(m) } &= P \gamma(\bayes_m) + \frac{V}{V-1} \times \frac{1}{n} \sum_{\lamm} \paren{1 + \delta^{(VF)}_{n,\pl}} \sigl^2 
\end{align}
where $\delta_{n,p}$ only depends on $(n,p)$, $\delta_{n,p}^{(VF)}$ depends on $(n,p)$ and the partition $(B_j)_{1 \leq j \leq V}$, but both are small when the product $np$ is large: 
\[ \absj{\delta_{n, p} } \leq L_1 \qquad \mbox{and} \qquad \absj{\delta_{n,p}^{(VF)}} \leq L_2 \croch{ \epsilon_n^{reg} + \max \paren{ (np)^{-1/4} , e^{-np (1-c_B)} } } \enspace , \]
where $L_1$ is a numerical constant ant $L_2$ only depends on $c_B$.
%
\end{proposition}
\begin{remark} Since we deal with histograms, $\ERM_m$ is not defined when $\min_{\lamm} \phl=0$, which occurs with positive probability. We then have to take a convention for $P\gamma\paren{\ERM_m}$ (on the event $\min_{\lamm} \set{\phl}=0$, which has generally a very small probability) so that it has a finite expectation. The same kind of problem occur with $\critVF$. See the proof of Prop.~\ref{VFCV.pro.EcritVFCV-EcritID}.
\end{remark}

\medskip

Prop.~\ref{VFCV.pro.EcritVFCV-EcritID} is consistent with Burman's asymptotic estimate of the bias of VFCV \cite{Bur:1989}. The major advance here is that it is non-asymptotic, and we have explicit upper bounds on the remainder terms (see the proof of Prop.~\ref{VFCV.pro.EcritVFCV-EcritID} in App.~\ref{VFCV.sec.proofs.E}).
It shows that the classical $V$-fold cross-validation overestimates the variance term $n^{-1} \sum_{\lamm} \sigl^2$, because it estimates the generalization ability of $\ERM^{(-j)}_m$, which is built upon less data than $\ERM_m$. This interpretation is consistent with the results of Shao \cite{Sha:1997} on linear regression, and van der Laan, Dudoit and Keles \cite{vdL_Dud_Kel:2004} in the density estimation framework. 

When $V$ stays bounded as $n$ grows to infinity, it is then natural to think that VFCV is underfitting, and thus be suboptimal for prediction. Since Prop.~\ref{VFCV.pro.EcritVFCV-EcritID} is non-asymptotic and quite accurate, we are now in position to prove such a result.
\begin{theorem} \label{VFCV.the.subopt}
Let $n \in \N$, $(X_i,Y_i)_{1 \leq i \leq n}$ be i.i.d. random variables, with $X \sim \mathcal{U}([0,1])$ and $Y = X + \sigma \epsilon$ with $\sigma>0$, $\E\croch{\epsilon\sachant X}=0$, $\E\croch{\epsilon^2\sachant X}=1$ and $\norm{\epsilon}_{\infty} < + \infty$. 
Let $\M_n = \set{1, \ldots , n}$ and $\forall \mM_n$, $S_m$ be the model of regular histograms with $D_m = m$ pieces on $\X = [0,1]$.
Let $V \in \set{2, \ldots, n}$ and $\paren{B_j}_{1 \leq j \leq V}$ be some partition of $\set{1, \ldots , n}$ such that for every $j$, $\absj{\card(B_j) - n V^{-1} } < 1$.

Then, there is an event of probability at least $1 - K_1 n^{-2}$ on which 
\begin{equation}
\label{VFCV.eq.subopt}
\perte{\ERM_{\mhVF}} \geq (1 + \kappa(V) - \ln(n)^{-1/5}) \inf_{\mM_n} \set{\perte{\ERM_m}} \enspace ,
\end{equation}
for some constant $\kappa(V)>0$ depending only on $V$ (and decreasing as a function of $V$), and a constant $K_1$ which depends on $\sigma$, $A$ and $V$. 
\end{theorem}
We now make a few comments:
\begin{itemize}
\item 
In the same framework, using similar arguments, we can prove an upper bound on $\perte{\ERM_{\mhVF}}$ showing that the constant $1+ \kappa(V)$ is exact (up to the $\ln(n)^{-1/5}$ term). In particular, 
\[ \frac{\perte{\ERM_{\mhVF}} } { \inf_{\mM_n} \set{\perte{\ERM_m}} } \xrightarrow[n \rightarrow + \infty]{a.s.} 1 + \kappa(V) = 1 + \frac{2^{2/3}}{3} \croch{ 1 - \paren{\frac{V-1}{V}}^{1/3} }^2   > 1 \enspace . \]
\item When $\paren{B_j}_{1 \leq j \leq V}$ is not assumed regular, the proof of Prop.~\ref{VFCV.pro.EcritVFCV-EcritID} shows that the factor $V/(V-1)$ becomes $\sum_{j=1}^V n/(n-\card(B_j))$ which is always larger, because $x \flapp (n-x)^{-1}$ is convex.
On the other hand, if one chooses a $(X_i)_{1 \leq i \leq n}$-dependent partition such that for every $\lamm$, $\card\set{X_i \in \Il  \, \mbox{ and } \, i \in B_j}$ is (almost) independent from $j$, then a similar proof shows that $\delta_{n,p}^{(VF)}$ is made much smaller than the previous upper bound. 
In a nutshell, it seems that the best performance of VFCV corresponds in general to the regular partition case, for which \eqref{VFCV.eq.subopt} holds.
\item Although we restrict in Thm.~\ref{VFCV.the.subopt} to a very particular problem, a similar result stays valid much more generally, possibly with a different value for the constant $\kappa(V)$. The only purpose of our assumptions is to compare very precisely $\critVF(m)$ and $P\gamma\paren{\ERM_m}$ as functions of $m$. Since $D_{\mhVF}$ is smaller than the optimum from a multiplicative factor independent from $n$ only, this analysis strongly depends on how $P \gamma\paren{\ERM_m}$ varies with $m$.
\item One can easily extend this result to any cross-validation like method, when two conditions are satisfied. First, the ratio between the size of the training set and $n$ has to be upperbounded by $1 - V^{-1} < 1$ (uniformly in $n$). Second, the number of training sets considered has to be bounded by $B_{\max}$ (from which $K_1$ may depend). This includes for instance the hold-out case, and repeated learning-testing methods. 
Notice that the second assumption is mainly technical; if we were able to prove the corresponding concentration inequalities, the leave-$p$-out with $p \sim n/V$ should have approximately the same properties.
\end{itemize}

\subsection{How to choose $V$} \label{VFCV.sec.subopt.VF.choix}
\subsubsection{Classical analysis}
There are three well-known factors to take into account in order to choose $V$:
\begin{itemize}
\item {\em bias}: when $V$ is too small, $\critVF$ overestimates the variance term in $P \gamma\paren{\ERM_m}$, which leads to underfitting and suboptimal model selection (Thm.~\ref{VFCV.the.subopt}).
\item {\em variability}: the variance of $\critVF(m)$ is a decreasing function of $V$, at least in the linear regression framework (see Burman \cite{Bur:1989} for an asymptotic expansion of this variance).
%
In general, $V=2$ is known to be quite variable because of the single split. When the prediction algorithm $\paren{X_i,Y_i}_{1 \leq i \leq n} \flapp \ERM_m$ is unstable (\latin{e.g.} classification with CART, as noticed by Hastie, Tibshirani and Friedman \cite{Has_Tib_Fri:2001}; see also Breiman \cite{Bre:1996}), the leave-one-out criterion (\latin{i.e.} $V=n$) is also known to be quite variable, but this phenomenon seems to disappear when $\ERM_m$ is more stable (Molinaro, Simon and Pfeiffer \cite{Mol_Sim_Pfe:2005}). In particular, in the least-squares regression framework, the variance of $\critVF(m)$ should decrease with $V$.
\item {\em computational complexity}: $V$-fold cross-validation needs to compute at least $V$ empirical risk minimizers for each model.
\end{itemize}
In the least-squares regression setting, $V$ has to be chosen large in order to improve accuracy (by reducing bias and variability); on the contrary, computational issues arise when $V$ is too big. This is why $V=5$ and $V=10$ are very classical and popular choices.

\subsubsection{The non-asymptotic need for overpenalization}
We now come to some particularity of the non-asymptotic viewpoint. Indeed, our proof of Thm.~\ref{VFCV.the.subopt} shows that the asymptotic behaviour of hold-out and cross-validation criterions only depend on their bias, because all these criterions are sufficiently close to their expectations asymptotically. However, this is not true when the sample size is fixed, and even the less variable criterions are far from being deterministic. As a consequence, using an unbiased estimator is no longer a guarantee of being optimal, since it can still lead to choosing a very poor model with a positive probability.

\begin{figure}
\centerline{\epsfig{file=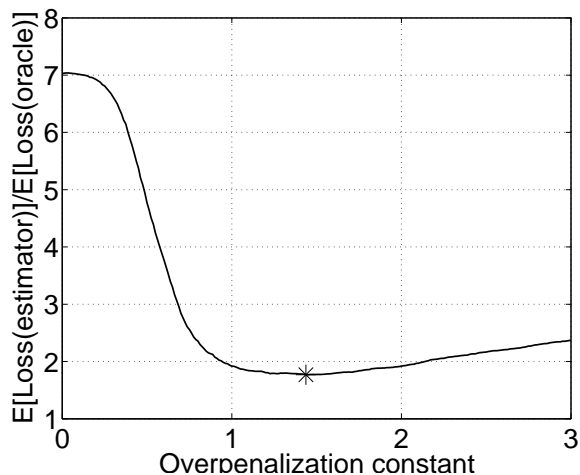,width=0.48\textwidth}}
\caption{The non-asymptotic need for overpenalization: the prediction performance $C_{\mathrm{or}}$ (defined in Sect.~\ref{VFCV.sec.simus.setup}) of the model selection procedure \eqref{def.mh.pen} with $\pen(m) = C_{\mathrm{ov}} \E\croch{\penid(m)}$ is represented as a function of $C_{\mathrm{ov}}$. Data and models are the ones of experiment (S1): $n=200$, $\sigma\equiv 1$, $\bayes(x) = \sin(\pi x)$. See Sect.~\ref{VFCV.sec.simus} for more details. \label{fig.S1.surpen.Epenid}}
\end{figure}

In order to analyze this phenomenon, it is useful to take the penalization viewpoint. The idea of penalization for model selection is to define 
\begin{equation} \label{def.mh.pen} \mh \in \arg\min_{\mM_n} \set{ P_n \gamma\paren{\ERM_m} + \pen(m) } \enspace , \end{equation}
where $\pen(m)$ is chosen so that $P_n \gamma\paren{\ERM_m} + \pen(m)$ is close to the prediction error $P \gamma\paren{\ERM_m}$. In other words, the ``ideal penalty'' is 
\begin{equation} \label{VFCV.def.penid}
 \penid(m) \egaldef (P-P_n) \gamma(\ERM_m) \enspace . \end{equation}
According to Prop.~\ref{VFCV.pro.EcritVFCV-EcritID} and \eqref{VFCV.eq.p2.hist} (which follows its proof), in the histogram regression case, we can compute the expectation of the ideal penalty:
\begin{equation} \label{VFCV.eq.Epenid.histo}
\E\croch{\penid(m)} = \frac{1}{n} \sum_{\lamm} \paren{2 + \delta_{n, \pl}} \sigl^2 \enspace ,
\end{equation}
which is close to Mallows' $C_p$ penalty $2\sigma^2 D_m n^{-1}$ in the homoscedastic case.
The point is that overpenalization (that is, taking $\pen$ larger than $\penid$, even in expectation) can improve the prediction performance of $\ERM_{\mh}$ when the signal-to-noise ratio is small. This can be seen on Fig.~\ref{fig.S1.surpen.Epenid}, according to which the optimal overpenalization constant $C_{\mathrm{ov}}^{\star}$ seems to be between $1.2$ and $1.7$ for this particular model selection problem. See also \cite{Arl:2007:phd} for a longer discussion of this problem.

\subsubsection{Choosing $V$ in the non-asymptotic framework}
Since $V$-fold cross-validation is choosing the model $\mhVF$ which minimizes some criterion $\critVF$, it can be written as a penalization procedure: it satisfies \eqref{def.mh.pen} with 
\[ \penVFCV(m) \egaldef \critVF(m) - P_n \gamma\paren{\ERM_m} \enspace . \]
Using again Prop.~\ref{VFCV.pro.EcritVFCV-EcritID} and \eqref{VFCV.eq.p2.hist}, we can compute its expectation: 
\[ \E\croch{ \penVFCV(m) } = \frac{1}{n} \sum_{\lamm} \croch{1 + \frac{V}{V-1} \paren{1 + \delta^{(VF)}_{n,\pl}} } \sigl^2 .
\]
Compared to \eqref{VFCV.eq.Epenid.histo}, this shows that $V$-fold cross-validation is overpenalizing within a factor $1 + 1/(2(V-1))$.

We can now revisit the question of choosing $V$ for optimal prediction, in such a non-asymptotic situation:
\begin{itemize}
\item the {\em overpenalization} factor is $1 + 1/(2(V-1))$.
\item the {\em variance} of $\critVF$ roughly decreases with $V$.
\item the {\em computational complexity} of computing $\critVF$ is roughly proportional to $V$.
\end{itemize}
First, take only the prediction performance into account. The variability question should be less crucial than overpenalization, because the variance of $\critVF$ depends only on $V$ through second order terms, according to the asymptotic computations of Burman \cite{Bur:1989}. Since the optimal overpenalization constant is $C_{\mathrm{ov}}^{\star}>1$, the performance of $V$-fold cross-validation should be optimal for some $V^{\star} < n$. This analysis is confirmed by the simulation study of Sect.~\ref{VFCV.sec.simus}, where $V=2$ provides better performance than $V=5$ and $V=10$ for several different experiments.

Now, if computational cost comes into the balance, or if we consider less stable prediction algorithms than least-squares regression estimators, the optimal $V$ may be even smaller. Whatever the framework, it seems quite difficult to find the optimal $V$, even if $C_{\mathrm{ov}}^{\star}$ was known (which is far from being the case in general). 
It would be at least necessary to understand well how the variance of $\critVF$ depends on $V$ in the non-asymptotic framework. This is a difficult practical problem, since ``there is no universal
(valid under all distributions) unbiased estimator of the variance of $V$-fold cross-validation'' (Bengio and Grandvalet \cite{Ben_Gra:2004}).
In the density estimation framework, this question has been tackled recently by Celisse and Robin \cite{Cel_Rob:2006}.

\medskip

The conclusion of this section is that choosing $V$ for $V$-fold is a very complex issue in practice, even independently from the cost of computing $\critVF$.
Moreover, it seems unsatisfactory to select a model according to a criterion as variable as the 2-fold cross-validation one when $V^{\star}=2$ because of the need for overpenalization. 
Finally, when the signal-to-noise ratio is large, we would like to obtain a nearly unbiased procedure without having to take $V$ very large, which can be computationally too heavy. 

In other words, we would like to decouple the choice of an overpenalization factor from the variability issue (which is essentially linked with complexity). The drawback of $V$-fold cross-validation is that they both depend on the $V$ parameter. As we shall see in the next section, such a decoupling can be naturally obtained through the use of penalization.

\section{An alternative $V$-fold algorithm: $V$-fold penalties} \label{VFCV.sec.penVFCV}
There are several ways to define $V$-fold cross-validation like penalization procedures with a tunable overpenalization factor, independent from the $V$ parameter. 
A first idea may be to multiply $\penVFCV (m) $ by a constant 
\latin{i.e.} to use \eqref{def.mh.pen} with the penalty 
\[ \pen(m) = C_{\mathrm{ov}} \paren{1 + \frac{1}{2 (V-1)}}^{-1} \paren{ \critVF(m) - P_n\gamma\paren{\ERM_m} } \enspace . \]
From the proof of Thm.~\ref{VFCV.the.subopt} (see also the one of Thm.~\ref{VFCV.the.histos.penVFCV} below), it is clear that when $C_{\mathrm{ov}} \sim 1$, this procedure satisfies with large probability a non-asymptotic oracle inequality with leading constant $1 + \epsilon_n$, and more generally an oracle inequality with leading constant $K(C_{\mathrm{ov}}) \geq 1$. However, this may seem a little artificial, and strongly dependent from the histogram regression framework in which the computations of Prop.~\ref{VFCV.pro.EcritVFCV-EcritID} work.

In this section, we consider another approach, that we call ``$V$-fold penalization'', which seems more natural to us. We shall see below that it is closely related to an idea of Burman \cite{Bur:1989,Bur:1990} for correcting the bias of $V$-fold cross-validation. However, Burman did not consider his method as a penalization one. His goal was only to obtain an unbiased estimate of the prediction error, so that it is not straightforward to choose an overpenalization factor different from~1 with his method. This is a major difference with our approach.

\subsection{Definition of $V$-fold penalties} \label{VFCV.sec.penVFCV.def}
\subsubsection{General framework}
We come back to the general setting of Sect.~\ref{VFCV.sec.subopt.VF.gal}. Recall that each predictor $\ERM_m$ can be written as a function $\ERM_m(P_n)$ of the empirical distribution of the data $P_n = n^{-1} \sum_{i=1}^n \delta_{(X_i,Y_i)}$. We want to build a penalization method, \latin{i.e.} choose $\mh$ according to \eqref{def.mh.pen}, so that the prediction error of $\ERM_{\mh}$ is as small as possible. This could be done exactly if we knew the ideal penalty $\penid(m) = (P-P_n) \gamma(\ERM_m(P_n)) $, but this quantity depends on the unknown distribution $P$. Following a heuristics due to Efron \cite{Efr:1979}, we propose to define $\pen$ as the resampling estimate of $\penid$, according to a $V$-fold subsampling scheme. We first recall the general form of this heuristics.

Basically, the {\em resampling heuristics} tells that one can mimic the relationship between $P$ and $P_n$ by building a $n$-sample of common distribution $P_n$ (the ``resample''). $\Pnb$ denoting the empirical distribution of the resample, the pair $(P,P_n)$ should be close (in distribution) to the pair $(P_n,\Pnb)$ (conditionally to $P_n$ for the latter distribution). Then, the expectation of any quantity of the form $F(P,P_n)$ can be estimated by $\Es \croch{ F(P_n, \Pnb) }$, where $\Es\croch{\cdot}$ denotes expectation w.r.t. the resampling randomness. In the case of $\penid$, this leads to Efron's bootstrap penalty \cite{Efr:1983}. 
Later on, this heuristics has been generalized to other resampling schemes, with the exchangeable weighted bootstrap (Mason and Newton \cite{Mas_New:1992}, Pr\ae stgaard and Wellner \cite{Pra_Wel:1993}). The empirical distribution of the resample then has the general form 
\[ \Pnb \egaldef \frac{1}{n} \sum_{i=1}^n W_i \delta_{(X_i,Y_i)} \qquad \mbox{with} \quad W \in \R^n \quad \mbox{an exchangeable weight vector,} \]
independent from the data ($W$ is said to be {\em exchangeable} when its distribution is invariant by any permutation of its coordinates). Fromont \cite{Fro:2007} used it successfully (with a particular upper bound on $\penid$) to build global penalties in the classification framework. Exchangeable resampling penalties (generalizing Efron's bootstrap penalty) have also been recently proposed, and studied in the regression framework \cite{Arl:2007:phd}. 
The idea of $V$-fold penalties is to use a $V$-fold subsampling scheme 
instead, \latin{i.e.} take $W_i = \frac{V}{V-1} \1_{i \notin B_J}$ with $J \sim \mathcal{U}(\set{1, \ldots, V})$ independent from the data ($\mathcal{U}(E)$ denotes the uniform distribution over the set $E$). Then, $\Pnb = P_n^{(-J)}$ and we obtain the following algorithm.
\begin{proc}[$V$-fold penalization] \label{VFCV.def.proc.gal.VFold}
\hfill
\begin{enumerate}
\item Choose a partition $(B_j)_{1 \leq j \leq V}$ of $\set{1,\ldots,n}$, as regular as possible.
\item Choose a constant $C \geq \CWinf = V-1$.
\item Compute the following resampling penalty for each $\mM_n$:
\begin{equation*} 
\pen(m) = \penVF(m) \egaldef \frac{C}{V} \sum_{j=1}^V \croch{  P_n \gamma \paren{ \ERM_m \paren{ P_n^{(-j)} } } - P_n^{(-j)} \gamma \paren{\ERM_m \paren{ P_n^{(-j)} }} } \enspace .
\end{equation*}
\item Choose $\mh$ according to \eqref{def.mh.pen}.
\end{enumerate}
\end{proc}
%
%
\begin{remark}[About the constant $C$] \label{VFCV.rem.Burman=penVF}
Contrary to Efron's resampling heuristics, we have to put a constant $C \neq 1$ in front of the penalty ($\pen$ being an unbiased estimator of $\penid$ when $C=\CWinf$). This is because each $W_i$ has a variance $(V-1)^{-1} \neq 1$ (we only normalized $W$ so that $\E\croch{W_i}=1$ for every $i$). According to Lemma~8.4 of \cite{Arl:2007:phd}, the right normalizing constant can be derived from the exchangeable case. As a consequence, from Theorem 3.6.13 in \cite{vdV_Wel:1996}, 
\[ \CWinf \sim_{n \rightarrow \infty} \left( n^{-1} \sum_{i=1}^n \E\paren{W_i-1}^2  \right)^{-1} \sim_{n \rightarrow \infty} V-1 \enspace . \]
The asymptotic value of $\CWinf$ can also be derived from the computations of Burman \cite{Bur:1989} in the linear regression framework.
Indeed, with our notations, Burman's criterion (formula (2.3) in \cite{Bur:1989}) is
\begin{align*}
\crit_{\mathrm{corr. VF}}(m) &\egaldef \critVF(m) + P_n \gamma\paren{\ERM_m} - \frac{1}{V} \sum_{j=1}^V P_n \gamma\paren{\ERM_m^{(-j)}} \\
&= P_n \gamma\paren{\ERM_m} + \frac{1}{V} \sum_{j=1}^V \croch{ \paren{ P_n^{(j)}  - P_n} \gamma\paren{\ERM_m^{(-j)}} } \enspace .
\end{align*}
If all the blocks of the partition have the same size $n/V$, then $P_n^{(j)} - P_n = (V-1) (P_n - P_n^{(-j)})$, so that Burman's corrected VFCV coincides exactly with $V$-fold penalization when $C = V-1$. Since $\crit_{\mathrm{corr. VF}}(m)$ is an asymptotically unbiased estimator of $P\gamma\paren{\ERM_m}$ (at least for linear regression), the result follows.
From the non-asymptotic viewpoint, we prove in Sect.~\ref{VFCV.sec.penVFCV.exp} below that $V-1$ also leads to an unbiased estimator of $\penid$ in the histogram regression case.
\end{remark}

Notice also that we do not assume that $C = \CWinf$, but only $C \geq \CWinf$. This is a major quality of $V$-fold penalization (penVF): it is straightforward to choose any overpenalization factor, independently from $V$.
Further comments about the choice of $C$ and $V$ are made in Sect.~\ref{VFCV.sec.disc}.

\subsubsection{The histogram regression case}
We now come back to the framework of Sect.~\ref{VFCV.sec.VFCV.calculs}, in which we can analyze deeper Algorithm~\ref{VFCV.def.proc.gal.VFold}. Remind that histograms are not our final goal, but only a convenient setting from which we can derive heuristics for practical use of penVF in any framework. 
From now on, $\paren{S_m}_{\mM_n}$ is a collection of histogram models  and $\paren{\ERM_m}_{\mM_n}$ the associated collection of least-squares estimators.
We first introduce some more notations:
\begin{gather*}
\bayes_m = \sum_{\lamm} \betl \1_{\Il} \mbox{~~and~~} \ERM_m = \sum_{\lamm} \bethl \1_{\Il} \quad \mbox{with} \quad \betl = \E \croch{ Y \sachant X \in \Il } \mbox{~~and~~} \bethl = \frac{1}{n\phl} \sum_{X_i \in \Il} Y_i   \\
\phlW \egaldef \Pnb(X \in \Il) = \phl \Wl \quad \mbox{with} \quad \Wl \egaldef \frac{1}{n \phl} \sum_{X_i \in \Il} W_i \\
\mbox{and} \qquad \ERMb_m \egaldef \arg\min_{t \in S_m} \Pnb \gamma(t) = \sum_{\lamm} \bethlW \1_{\Il} \quad \mbox{with} \quad \bethlW \egaldef \frac{1}{n\phlW} \sum_{X_i \in \Il} W_i Y_i \enspace . 
\end{gather*}

Assuming that $\min_{\lamm} \phl >0$ (otherwise, the model $m$ should clearly not be chosen), we can compute the ideal penalty (see \eqref{VFCV.eq.p1.hist} and \eqref{VFCV.eq.p2.hist} in Sect.~\ref{VFCV.sec.proofs.E}) and its resampling estimate:
\begin{align} \notag 
\penid(m) = (P-P_n) \gamma(\ERM_m) &= \sum_{\lamm} \paren{ \pl + \phl } \paren{  \bethl - \betl }^2 + (P-P_n) \gamma(\bayes_m) \\
\Es \croch{ (P_n-\Pnb ) \gamma(\ERMb_m) } 
&= \sum_{\lamm} \Es \croch{ \paren{\phl + \phlW} \paren{ \bethlW - \bethl }^2 } \enspace , 
\label{VFCV.eq.pen.pb.his} \end{align}
since $\sum_i \E[W_i]=1$ implies that $\Es \croch{ (P_n-\Pnb) \gamma(\ERM_m) }=0$.
The penalty \eqref{VFCV.eq.pen.pb.his} is well-defined if and only if $\ERMb_m$ is a.s. uniquely defined, \latin{i.e.} $\Wl >0$ for every $\lamm$ a.s. This is why we modified the definition of the weights in algorithm~\ref{VFCV.def.proc.gal.VFold}, so that this problem does not occur.
\begin{proc}[$V$-fold penalization for histograms] \label{VFCV.def.proc.his} 
\hfill
\begin{enumerate}
\item Replace $\M_n$ by 
$\Mh_n = \set{ \mM_n \telque \min_{\lamm} \set{ n\phl } \geq 3 }$.
\item Choose a constant $C \geq \CWinf = V-1$.
\item For every $\mMh_n$, choose a partition $(B_j)_{1 \leq j \leq V}$ of $\set{1,\ldots,n}$ such that 
\[ \forall \lamm, \, \forall 1 \leq j \leq V, \quad \absj{ \card \paren{ B_j \cap \set{i \telque X_i \in \Il} } - \frac{n\phl}{V} } < 1 \enspace . \]
%
\item Compute the following resampling penalty for each $\mM_n$:
\begin{equation} \label{VFCV.def.penVFCV.his}
\pen(m) = \penVF(m) \egaldef \frac{C}{V} \sum_{j=1}^V \croch{  P_n \gamma \paren{ \ERM_m^{(-j)} } - P_n^{(-j)} \gamma \paren{\ERM_m^{(-j)} } } \enspace .
\end{equation}
\item Choose $\mh$ according to \eqref{def.mh.pen}.
\end{enumerate}
\end{proc}
At step 3, we choose a different partition for each model $m$. Our choice is consistent with the proposal of Breiman \latin{et al.} \cite{Bre_etal:1984} (see also Burman \cite{Bur:1990}, Sect.~2) to stratify the data and choose a partition which respects the stratas. In the histogram case, natural stratas are the sets $\set{ i \telque X_i \in \Il}$. In particular, steps~1 and~3 of Algorithm~\ref{VFCV.def.proc.his}  ensure that $\min_{\lamm} \Wl >0$ for every $m \in \Mh_n$, so that \eqref{VFCV.def.penVFCV.his} is well-defined.

Other modifications of algorithm~\ref{VFCV.def.proc.gal.VFold} are possible. For instance, keep the same regular partition $(B_j)_{1 \leq j \leq V}$ for all the models, and take
\begin{equation}  \label{VFCV.def.penVFCV.his.2}
\penVF (m) = C \sum_{\lamm} \paren{ \Es \croch{ \phl \paren{ \bethlW - \bethl }^2 \sachant \Wl >0 } + \Es \croch{ \phlW \paren{ \bethlW - \bethl }^2 }} 
\end{equation}
instead of \eqref{VFCV.eq.pen.pb.his}. This is what we did in the simulations of Sect.~\ref{VFCV.sec.simus}, and a short theoretical study of this method is done in Sect.~8.4.1 of \cite{Arl:2007:phd}.
It confirms that the two algorithms should have very similar performances in practical applications.
 
\subsection{Expectations} \label{VFCV.sec.penVFCV.exp}
We now come to the expectation of $V$-fold penalties, in the histogram regression framework.
\begin{proposition} \label{VFCV.pro.EcritpenVFCV}
Let $S_m$ be the model of histograms associated with some partition $\paren{\Il}_{\lamm}$ and $\pen=\penVF$ be defined as in Algorithm~\ref{VFCV.def.proc.his}. Then, if $\min_{\lamm} \set{n\phl} \geq 3$,
\begin{equation} \label{VFCV.eq.EcritpenVFCV}
\El \croch{ \penVF (m) } = \frac{1}{n} \sum_{\lamm} \paren{ \frac{2C}{V-1} + \frac{C}{V-1} \delta_{n,\phl}^{(\mathrm{penV})} } \sigl^2  
\end{equation}
with $\El \croch{\cdot} = \El \croch{ \cdot \sachantd \paren{\1_{X_i \in \Il}}_{1 \leq i \leq n, \, \lamm} }$ and $\frac{2}{n\phl - 2} \geq \delta_{n,\phl}^{(\mathrm{penV})} \geq 0 $.
\end{proposition}

Comparing \eqref{VFCV.eq.EcritpenVFCV} with \eqref{VFCV.eq.Epenid.histo}, it appears that $\penVF$ is an (almost) unbiased estimator of $\penid$ when $C=V-1$. 
Indeed, when $\min_{\lamm} \set{n \pl}$ goes to infinity faster than some constant times $\ln(n)$, so does $\min_{\lamm} \set{n \phl}$ with a large probability. Moreover, following the proof of Lemma~\ref{VFCV.le.einv.binom}, we can show that 
\[ \E \croch{ \delta_{n,\phl}^{(\mathrm{penV})} \1_{n \phl \geq 3} } \leq \kappa \min\paren{1, (n\pl)^{-1/4} } \xrightarrow[n \pl \rightarrow \infty]{} 0 \] for some absolute constant $\kappa>0$. 
This is consistent with the asymptotic computations of Burman \cite{Bur:1989}. The main novelty of Prop.~\ref{VFCV.pro.EcritpenVFCV} is that we have an explicit non-asymptotic upperbound on the remainder term. This is crucial to derive oracle inequalities for Algorithm~\ref{VFCV.def.proc.his}.

\subsection{Oracle inequalities and asymptotic optimality} \label{VFCV.sec.theorie}
We are now in position to state the main result of this section: $V$-fold penalties (Algorithm~\ref{VFCV.def.proc.his}) satisfy a non-asymptotic oracle inequality with a leading constant close to 1, on a large probability event. This implies the asymptotic optimality of Algorithm~\ref{VFCV.def.proc.his} in terms of excess loss. For this, we assume the existence of some non-negative constants $\aM$, $\cM$, $c_{\mathrm{rich}}$, $\eta$ such that:
\begin{enumerate}
\item[\hypPpoly] Polynomial complexity of $\M_n$: $\card(\M_n) \leq \cM n^{\aM}$.
\item[\hypPrich] Richness of $\M_n$: $\exists m_0 \in \M_n$ s.t. $D_{m_0} \in \croch{ \sqrt{n}; c_{\mathrm{rich}} \sqrt{n} }$.
\item[\hypPconst] The constant $C$ is well chosen: $\eta (V-1) \geq C \geq V-1$.
\end{enumerate}

\begin{theorem} \label{VFCV.the.histos.penVFCV}
Assume that the $(X_i,Y_i)$'s satisfy the following:
\begin{enumerate}
\item[\hypAb] Bounded data: $\norm{Y_i}_{\infty} \leq A < \infty$.
\item[\hypAn] Noise-level bounded from below: $\sigma(X_i) \geq \sigmin>0$ a.s.
\item[\hypAp] Polynomial decreasing of the bias: there exists $\beta_1\geq\beta_2>0$ and $\cbiasmaj,\cbiasmin>0$ such that \[ \cbiasmin D_m^{-\beta_1} \leq \perte{\bayes_m} \leq \cbiasmaj D_m^{-\beta_2} \enspace . \]
\item[\hypArXl] Lower regularity of the partitions for $\loi(X)$: $ D_m \min_{\lamm} \pl \geq \crXl >0$.
\end{enumerate}

Let $\mh$ be the model chosen by algorithm~\ref{VFCV.def.proc.his} (under restrictions \hypPtoutVF, with $\eta=1$). Then, there exists a constant $K_2$ and a sequence $\epsilon_n$ converging to zero at infinity such that 
\begin{equation} \label{VFCV.eq.oracle_traj_non-as}
\perte{\ERM_{\mh}} \leq \paren{ 1 + \epsilon_n } \inf_{\mM_n} \set{ \perte{\ERM_m} } \end{equation}
with probability at least $1 - K_2 n^{-2}$.
Moreover, we have the oracle inequality
\begin{equation} \label{VFCV.eq.oracle_class_non-as}
\E \croch{\perte{\ERM_{\mh}}} \leq \paren{ 1 + \epsilon_n } \E \croch{ \inf_{\mM_n} \set{ \perte{\ERM_m} } } + \frac{A^2 K_2}{n^2} \enspace .
\end{equation} 

The constant $K_2$ may depend on $V$ and constants in \hypAb, \hypAn, \hypAp, \hypArXl\ and \hypPtoutVF, but not on $n$. The term $\epsilon_n$ is smaller than $\ln(n)^{-1/5}$ for instance; it can also be taken smaller than $n^{-\delta}$ for any $0 < \delta < \delta_0 (\beta_1, \beta_2)$, at the price of enlarging $K_2$.
\end{theorem}

We first make a few comments on our assumptions.
\begin{enumerate}
\item When assumption \hypPconst\ is satisfied with $\eta>1$, the same result holds with a leading constant $2\eta - 1 + \epsilon_n$ instead of $1 + \epsilon_n$ in \eqref{VFCV.eq.oracle_traj_non-as} and \eqref{VFCV.eq.oracle_class_non-as}.
\item In Thm.~\ref{VFCV.the.histos.penVFCV}, we assume that $V$ is fixed when $n$ grows. A careful look at the proof shows that we only need $V \leq \ln(n)$ for $n$ large enough. With a few more work, we could go up to $V$ of order $n^{\delta}$ for some $\delta>0$ depending on the assumptions of Thm.~\ref{VFCV.the.histos.penVFCV}, but we can not handle the leave-one-out case ($V=n$). This is probably a technical restriction, since a similar result for several exchangeable weights (including leave-one-out) is proven in Chap.~6 of \cite{Arl:2007:phd}.
\item \hypAb\ and \hypAn\ are rather mild (and neither $A$ nor $\sigmin$ need to be known from the statistician). In particular, they allow quite general heteroscedastic noises. They can even be relaxed, for instance thanks to results proven in Chap.~6 and Sect.~8.3 of \cite{Arl:2007:phd}, allowing the noise to vanish or to be unbounded.
\item \hypArXl\ is satisfied for ``almost regular'' histograms when $X$ has a lower bounded density w.r.t. $\Leb$, as for instance all the simulation experiments of Sect.~\ref{VFCV.sec.simus}.
\item The upper bound in \hypAp\ holds when $(\Il)_{\lamm}$ is regular and $\bayes$ $\alpha$-h\"olderian with $\alpha \in (0,1]$. The lower bound may seem more surprising, since it means that $\bayes$ is not too well approximated by the models $S_m$. 
However, it is classical to assume that $\perte{\bayes_m}>0$ for every $\mM_n$ for proving the asymptotic optimality of Mallows' $C_p$ (\latin{e.g.} by Shibata \cite{Shi:1981}, Li \cite{KCLi:1987} and Birg\'e and Massart \cite{Bir_Mas:2006}). We here make a stronger assumption because we need a non-asymptotic lower bound on the dimension of both the oracle and selected models. 
The reason why it is not too restrictive is that non-constant $\alpha$-h\"olderian functions satisfy \hypAp\ with \[ \beta_1 = k^{-1} + \alpha^{-1} - (k-1) k^{-1} \alpha^{-1} \qquad \mbox{and} \qquad \beta_2 = 2 \alpha k^{-1} \enspace , \] when $(\Il)_{\lamm}$ is regular and $X$ has a lower-bounded density w.r.t. the Lebesgue measure on $\X \subset \R^k$ (\latin{cf.} Sect.~8.10 in \cite{Arl:2007:phd} for more details). Notice also that Stone \cite{Sto:1985} and Burman \cite{Bur:2002} used the same assumption in the density estimation framework. 
\end{enumerate}

\medskip

Theorem~\ref{VFCV.the.histos.penVFCV} has at least two major consequences. 
First, $V$-fold penalties provide an {\em asymptotically optimal} model selection procedure, at least in the histogram regression framework, as soon as $C \sim V-1$. This should be compared to Thm.~\ref{VFCV.the.subopt}, where we proved that $V$-fold cross-validation is suboptimal for a rather mild homoscedastic problem.
Notice that a slight modification of the proof of Thm.~\ref{VFCV.the.histos.penVFCV} shows that several other cross-validation like methods (even with the same computational cost) have 
similar theoretical properties. We discuss this point in Sect.~\ref{VFCV.sec.disc}.

Second, Thm.~\ref{VFCV.the.histos.penVFCV} can handle several kinds of {\em heteroscedastic noises}, while Algorithm~\ref{VFCV.def.proc.his} does not need any knowledge about $\sigma$, $\norm{Y}_{\infty}$ or the smoothness of $\bayes$. Even the tuning of $C$ and $V$ can be made (at least at first order) without any information on the distribution $P$ of the data.
This shows that $V$-fold penalization is a {\em naturally adaptive algorithm}, as long as $\M_n$ allows adaptation. The point here is that 
when $\bayes$ belongs to some h\"olderian ball $\mathcal{H}(\alpha,R)$ (with $\alpha \in (0,1]$ and $R>0$), we can choose $\M_n$ as the family of regular histograms on $\X \subset \R^k$ to obtain such an adaptivity result. Then, from Thm.~\ref{VFCV.the.histos.penVFCV}, we can build an estimator adaptive to $(\alpha,R)$ in a heteroscedastic framework (see \cite{Arl:2007:phd} for more details). If moreover the noise-level $\sigma$ satisfies some regularity assumption, we can show that this estimator attains the minimax estimation rate, up to some numerical constant, when $\alpha=k=1$.

Notice also that a similar adaptation result could be obtained with $V$-fold cross-validation, which also satisfies \eqref{VFCV.eq.oracle_traj_non-as} and \eqref{VFCV.eq.oracle_class_non-as} with leading constants $K(V)>1$, under similar assumptions. 
The advance with $V$-fold penalization is that we have simultaneously the adaptivity property of $V$-fold cross-validation, its mild computational cost (when $V$ is chosen small), and asymptotic optimality (contrary to VFCV).

Finally, we would like to emphasize that building such estimators is not the final goal of penVF. As a matter of fact, there are several procedures that are adaptive to the smoothness of $\bayes$ and the heteroscedasticity of the noise (\latin{e.g.} by Efromovich and Pinsker \cite{Efr_Pin:1996} or Galtchouk and Pergamenshchikov \cite{Gal_Per:2005}), and they may have better performances than both VFCV and penVF in this particular framework. Contrary to these \latin{ad hoc} procedures, particulary built for dealing with heteroscedasticity, VFCV and penVF are general-purpose devices. What our theoretical results show is that they behave quite well in this framework, for which they were not built in particular.

\section{Simulation study} \label{VFCV.sec.simus}
As an illustration of the results of the two previous sections, we compare the performances of VFCV, penVF (for several values of $V$) and   Mallows' $C_p$ on some simulated data. 

\subsection{Experimental setup} \label{VFCV.sec.simus.setup}
\begin{figure}
\begin{minipage}[b]{.48\linewidth}   \centerline{\epsfig{file=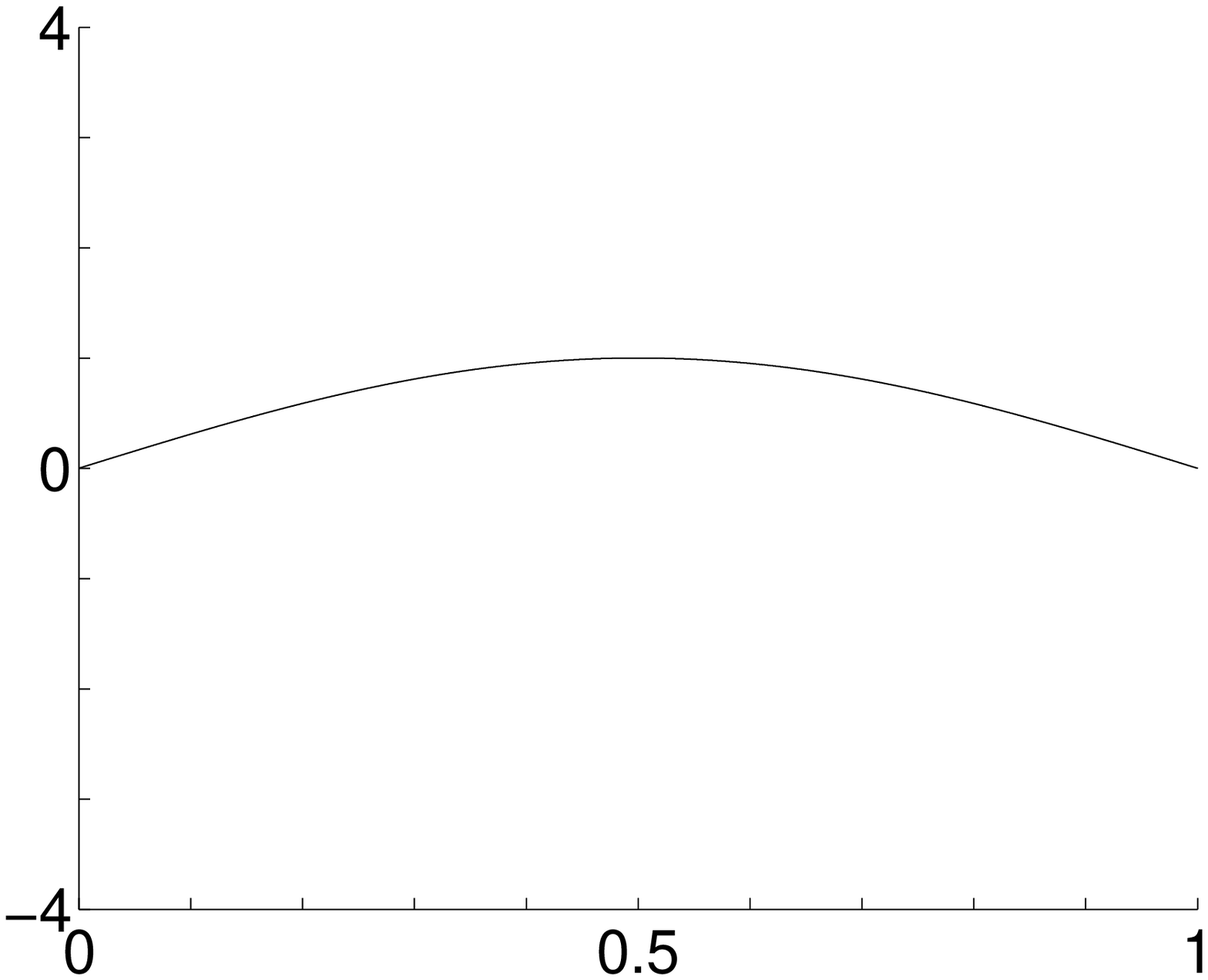,width=0.8\textwidth}}
   \caption{$\bayes(x) = \sin(\pi x)$\label{VFCV.fig.sin.fonc}}
\end{minipage} \hfill
 \begin{minipage}[b]{.48\linewidth}   \centerline{\epsfig{file=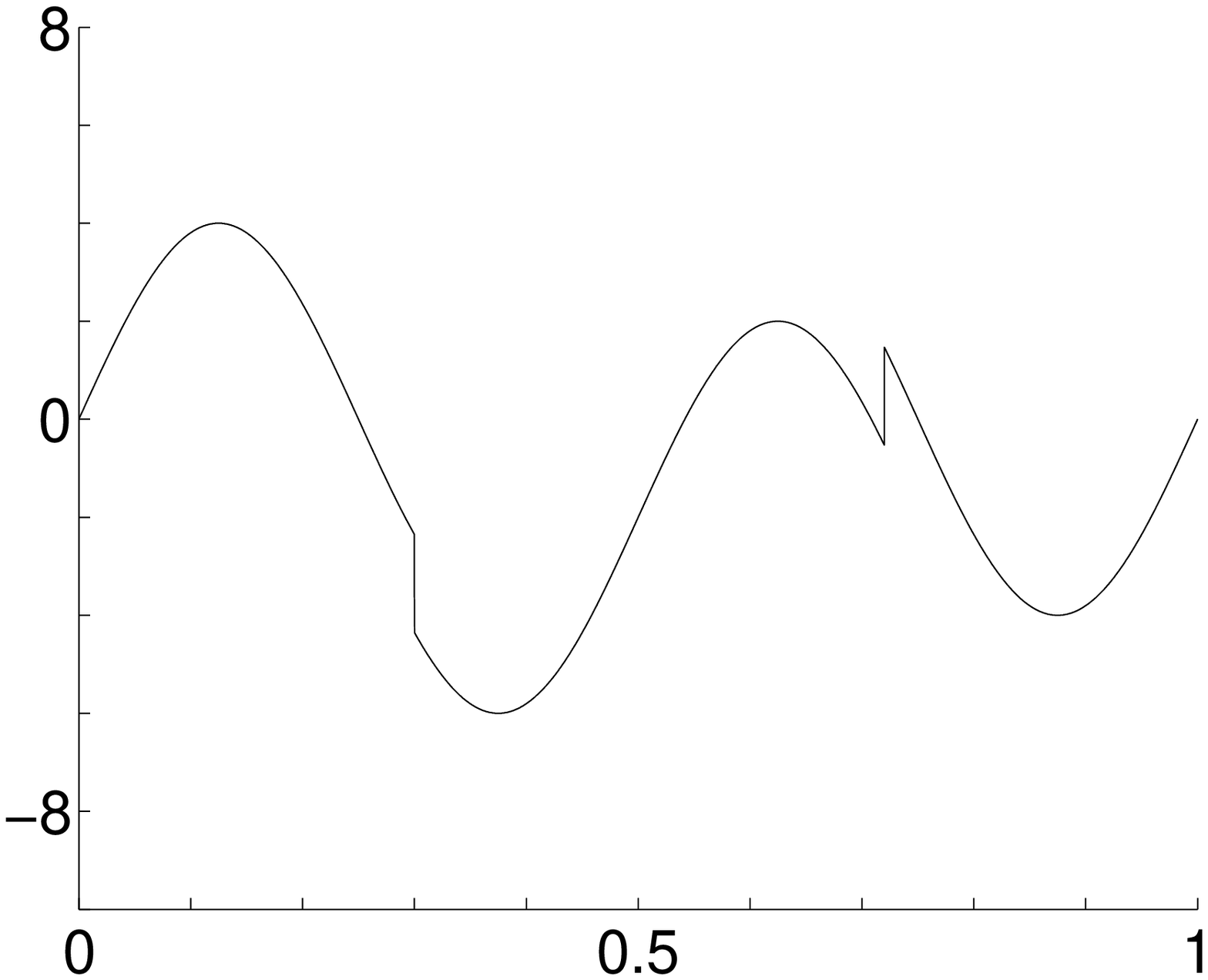,width=0.8\textwidth}}
   \caption{$\bayes(x) = \mathrm{HeaviSine}(x)$ (see \cite{Don_Joh:1995})\label{VFCV.fig.Hsin.fonc}}
\end{minipage}
\end{figure}

\begin{figure}
\begin{minipage}[b]{.48\linewidth}   \centerline{\epsfig{file=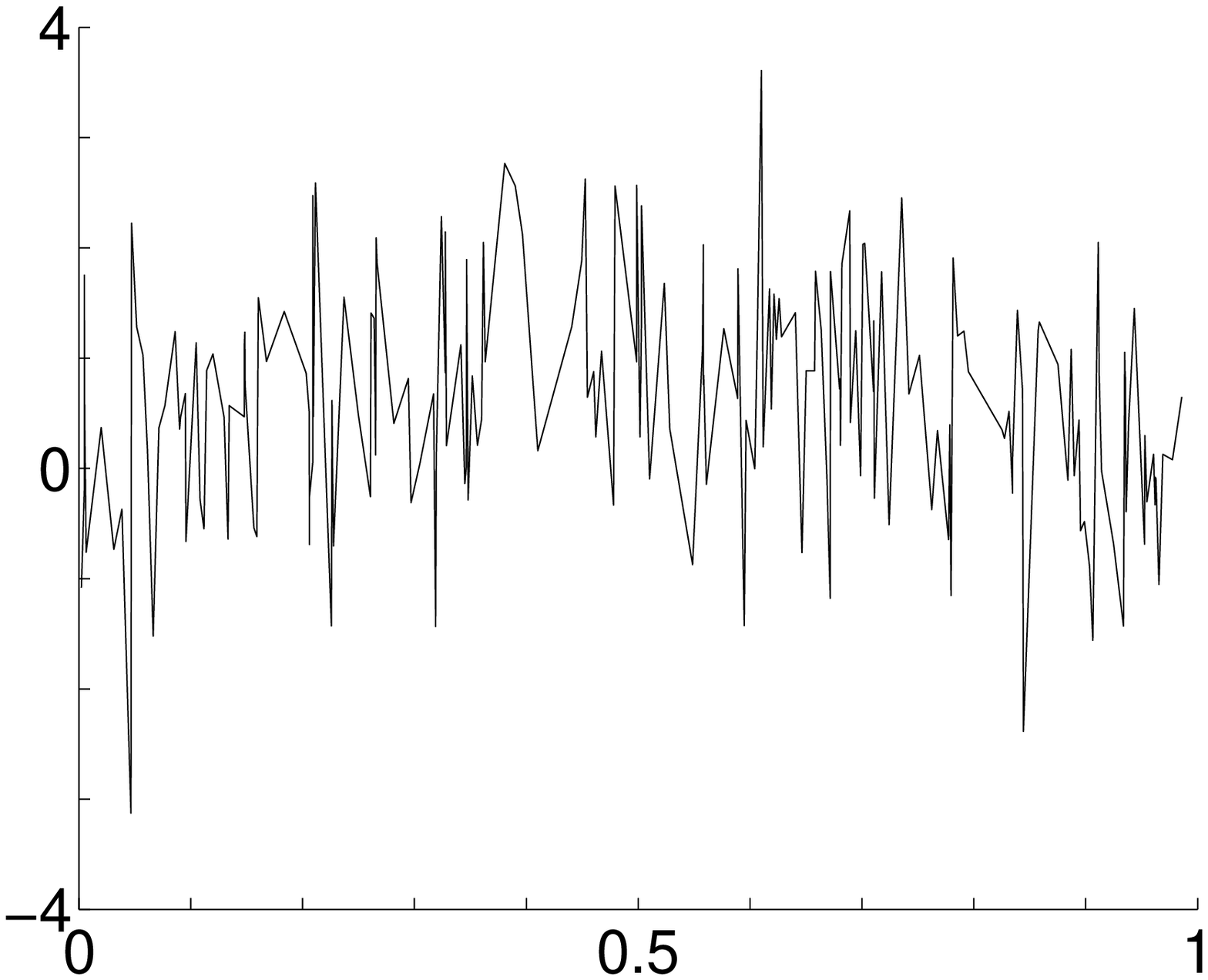,width=0.8\textwidth}}
   \caption{S1: $\bayes(x)=\sin(\pi x)$, $\sigma \equiv 1$, $n=200$}
   \label{VFCV.fig.S1.data}
\end{minipage} \hfill
 \begin{minipage}[b]{.48\linewidth}   \centerline{\epsfig{file=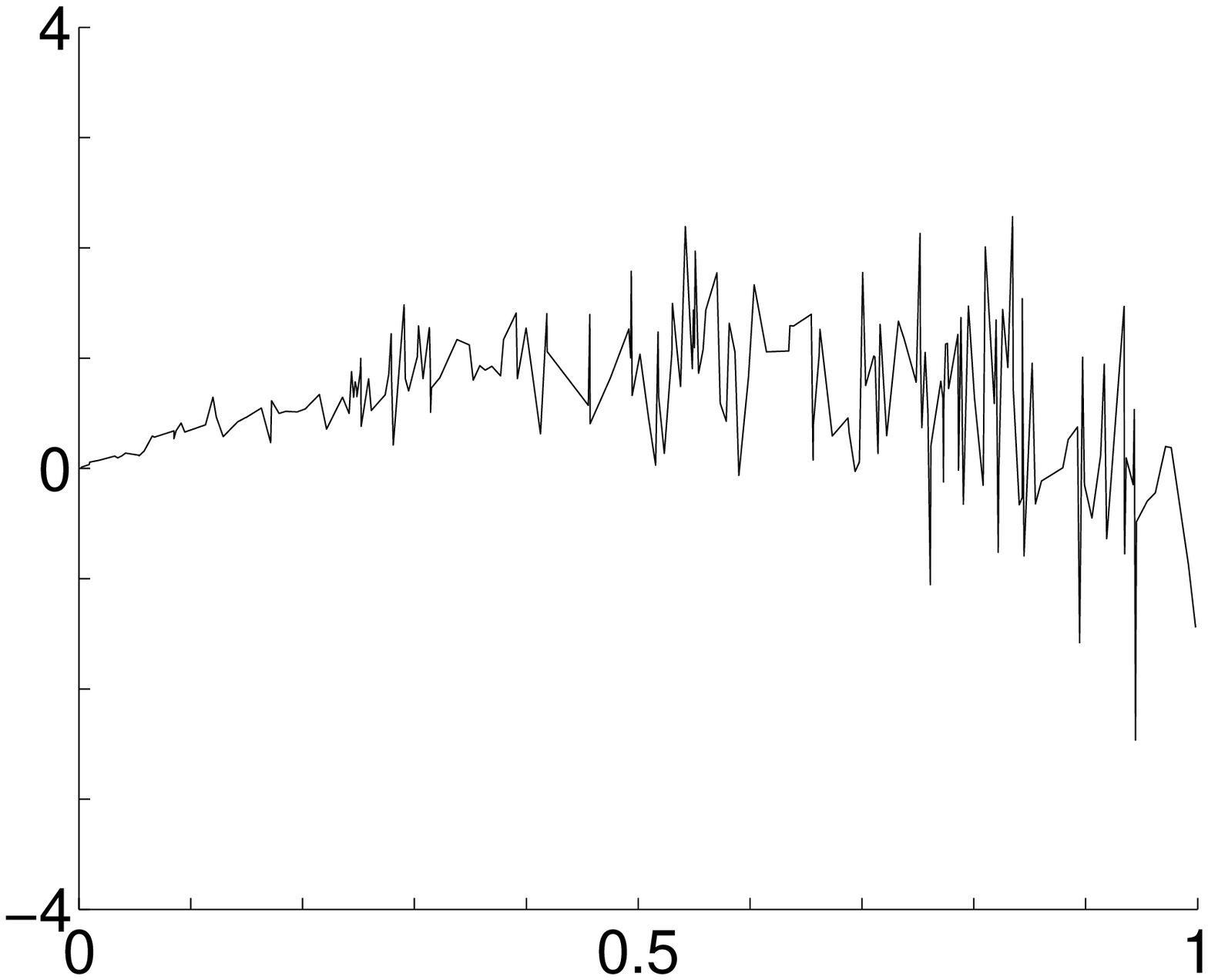,width=0.8\textwidth}}
   \caption{S2: $\bayes(x)=\sin(\pi x)$, $\sigma(x)=x$, $n=200$}
   \label{VFCV.fig.S2.data}
\end{minipage}
\end{figure}

\begin{figure}
\begin{minipage}[b]{.48\linewidth}   \centerline{\epsfig{file=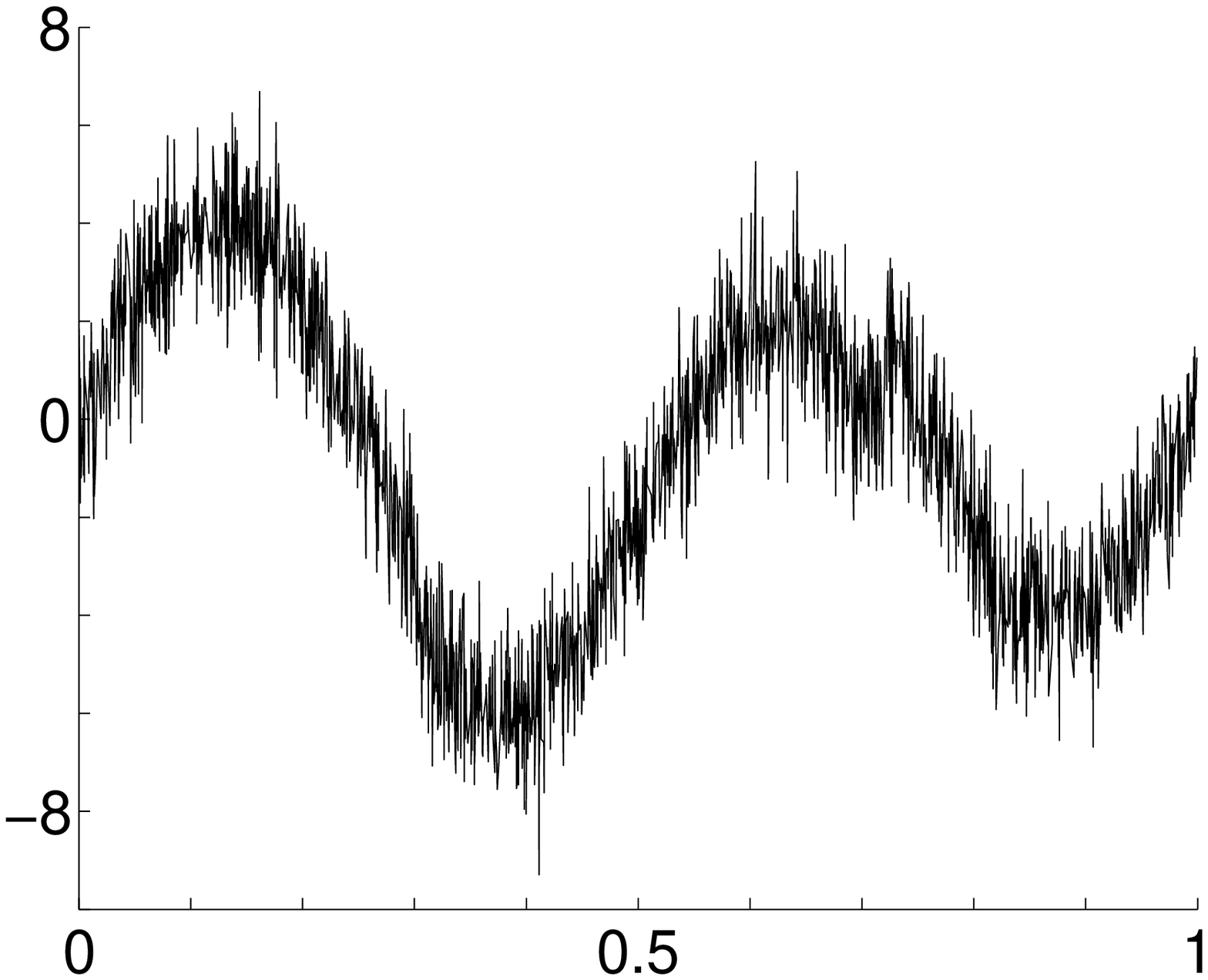,width=0.8\textwidth}}
   \caption{HSd1: HeaviSine, $\sigma\equiv 1$, $n=2048$}
   \label{VFCV.fig.HSd1.data}
\end{minipage} \hfill
 \begin{minipage}[b]{.48\linewidth}   \centerline{\epsfig{file=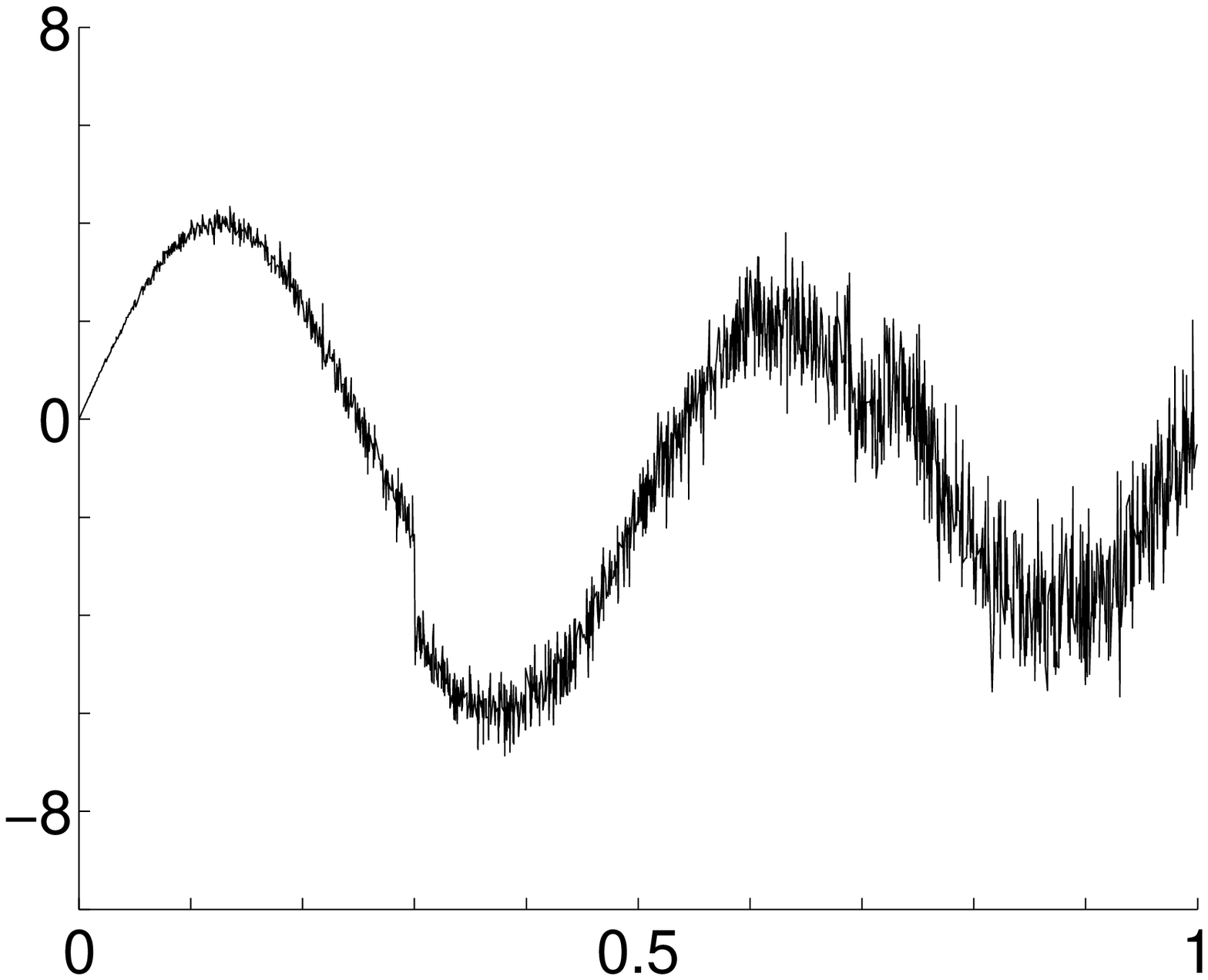,width=0.8\textwidth}}
   \caption{HSd2: HeaviSine, $\sigma(x)=x$, $n=2048$}
   \label{VFCV.fig.HSd2.data}
\end{minipage}
\end{figure}

We consider four experiments, called S1, S2, HSd1 and HSd2. Data are generated according to 
\[ Y_i = \bayes(X_i) + \sigma(X_i) \epsilon_i \]
with $X_i$ i.i.d. uniform on $\X=[0;1]$ and $\epsilon_i \sim \mathcal{N}(0,1)$ independent from $X_i$. 
The experiments differ from the regression function $\bayes$ (smooth for S, see Fig.~\ref{VFCV.fig.sin.fonc}; smooth with jumps for HS, see Fig.~\ref{VFCV.fig.Hsin.fonc}), the noise type (homoscedastic for S1 and HSd1, heteroscedastic for S2 and HSd2) and the number $n$ of data. Instances of data sets are given by Fig.~\ref{VFCV.fig.S1.data} 
to~\ref{VFCV.fig.HSd2.data}.
Their last difference lies in the families of models. 
Defining \begin{gather*} 
\forall k, k_1, k_2 \in \N \backslash \set{0}, \quad \paren{\Il}_{\lambda \in \Lambda_k} = \paren{ \left[ \frac{j}{k}; \frac{j+1}{k} \right) }_{0 \leq j \leq k-1} \mbox{ and }  \\
\paren{\Il}_{\lambda \in \Lambda_{(k_1,k_2)}} = \paren{ \left[ \frac{j}{2k_1}; \frac{j+1}{2k_1} \right) }_{0 \leq j \leq k_1-1} \cup \paren{ \left[ \frac{1}{2} + \frac{j}{2k_2}; \frac{1}{2} + \frac{j+1}{2k_2} \right) }_{0 \leq j \leq k_2-1} \enspace ,
\end{gather*}
the four model families are indexed by $m \in \M_n \subset \paren{\N \backslash \set{0}} \cup \paren{\N \backslash \set{0}}^2$:
\begin{enumerate}
\item[S1] regular histograms with $1 \leq D \leq n (\ln(n))^{-1}$ pieces, \latin{i.e.} \[ \M_n = \set{ 1, \ldots, \left\lfloor \frac{n}{\ln(n)} \right\rfloor } \enspace . \]
\item[S2] histograms regular on $\croch{0;1/2}$ (resp. on $\croch{1/2;1}$), with $D_1$ (resp. $D_2$) pieces, $1 \leq D_1,D_2 \leq n (2 \ln(n))^{-1}$. The model of constant functions is added to $\M_n$, \latin{i.e.} \[ \M_n = \set{1} \cup \set{ 1, \ldots, \left\lfloor \frac{n}{2\ln(n)} \right\rfloor }^2 \enspace . \]
\item[HSd1] dyadic regular histograms with $2^k$ pieces, $0 \leq k \leq \ln_2(n) - 1$, \latin{i.e.} \[ \M_n = \set{ 2^k \telque 0 \leq k \leq \ln_2(n)-1 } \enspace .\]
\item[HSd2] dyadic regular histograms with bin sizes $2^{-k_1}$ and $2^{-k_2}$, $0 \leq k_1,k_2 \leq \ln_2(n) - 2$ (dyadic version of S2). The model of constant functions is added to $\M_n$, \latin{i.e.} \[ \M_n = \set{1} \cup \set{ 2^k \telque 0 \leq k \leq \ln_2(n)-2 }^2 \enspace .\]
\end{enumerate}
Notice that we choose models that can approximately fit the true shape of $\sigma(x)$ in experiments S2 and HSd2. This choice makes the oracle  model even more efficient, hence the model selection problem more challenging.

We compare the following algorithms: 
\begin{enumerate}
\item[VFCV] Classical $V$-fold cross-validation, defined by \eqref{VFCV.def.mh.VFCvclass}, with $V \in \{ 2,5,10,20 \}$.
\item[LOO] Classical Leave-one-out (\latin{i.e.} VFCV with $V=n$).
\item[penVF] $V$-fold penalty, with $V \in \{ 2,5,10,20 \}$. $C= \CWinf =V-1$. The partition $(B_j)$ is chosen once, as in Algorithm~\ref{VFCV.def.proc.gal.VFold}, and $\penVF$ is defined by \eqref{VFCV.def.penVFCV.his.2}. In practice, this is almost the same as Algorithm~\ref{VFCV.def.proc.his}.
\item[penLoo] $V$-fold penalty, with $V = n$. $C= \CWinf =n-1$.
\item[Mal] Mallows' $C_p$ penalty: $\pen(m) = 2 \widehat{\sigma}^2 D_m n^{-1}$, where $\widehat{\sigma}^2= 2 n^{-1} d^2 \paren{ Y\nupl , S_{n/2} } $ is the classical variance estimator ($d$ being the Euclidean distance on $\R^n$, $S_{n/2}$ any vector space of dimension $n/2$ of $\R^n$ and $Y\nupl = (Y_1, \ldots, Y_n) \in \R^n$).
The non-asymptotic validity of this procedure for model selection in homoscedastic regression has been assessed by Baraud \cite{Bar:2000}.
\item[$\E\croch{\penid}$] Ideal deterministic penalty: $\pen(m) = \E\croch{\penid(m)}$. We use it as a witness of what is a good performance in each experiment.
\end{enumerate}
For each penalization procedure, we also consider the same penalty multiplied by $5/4$ (denoted by a $+$ symbol added after its shortened name). This intends to test for overpenalization (the choice of the factor $5/4$ being arbitrary and certainly not optimal).

In each experiment, for each simulated data set, we replace $\M_n$ by $\Mh_n$ as in step~1 of Algorithm~\ref{VFCV.def.proc.his}. Then, we compute the least-squares estimators $\ERM_m$ for each $\mMh_n$. Finally, we select $\mh \in \Mh_n$ using each algorithm and compute its true excess loss $\perte{\ERM_{\mh}}$ (and the excess loss $\perte{\ERM_m}$ for every $m \in \Mh_n$). 
We simulate $N=1000$ data sets, from which we can estimate the model selection performance of each procedure, through the two following benchmarks:
\begin{equation*} \label{VFCV.def.Cor}
C_{\mathrm{or}} = \frac{ \E\croch{ \perte{\ERM_{\mh}} }} {\E\croch{ \inf_{\mM_n}  \perte{\ERM_m} }} \qquad \mbox{and} \qquad C_{\mathrm{path-or}} = \E\croch{\frac{  \perte{\ERM_{\mh} }} {\inf_{\mM_n} \perte{\ERM_m} } } \enspace .
\end{equation*}
Basically, $C_{\mathrm{or}}$ is the constant that should appear in an oracle inequality like \eqref{VFCV.eq.oracle_class_non-as}, and $C_{\mathrm{path-or}}$ corresponds to a pathwise oracle inequality like \eqref{VFCV.eq.oracle_traj_non-as}. As $C_{\mathrm{or}}$ and $C_{\mathrm{path-or}}$ approximatively give the same rankings between algorithms, we only report $C_{\mathrm{or}}$ in Tab.~\ref{VFCV.tab.un}. 

\begin{table} 
\caption{Accuracy indexes $C_{\mathrm{or}}$ for each algorithm in four experiments, $\pm$ a rough estimate of uncertainty of the value reported (\latin{i.e.} the empirical standard deviation divided by $\sqrt{N}$). In each column, the more accurate algorithms (taking the uncertainty into account; $\E\croch{\penid}$ and $\E\croch{\penid}+$ are not taken into account there) are bolded.
\label{VFCV.tab.un}}
\begin{center}
\begin{tabular}
{p{0.16\textwidth}@{\hspace{0.025\textwidth}}p{0.16\textwidth}@{\hspace{0.025\textwidth}}p{0.16\textwidth}@{\hspace{0.025\textwidth}}p{0.17\textwidth}@{\hspace{0.025\textwidth}}p{0.17\textwidth}}
\hline\noalign{\smallskip}
Experiment  & S1      & S2        & HSd1            & HSd2 \\
\noalign{\smallskip}
\hline
\noalign{\smallskip}
$\bayes$    & $\sin(\pi\cdot)$ & $\sin(\pi\cdot)$ & HeaviSine & HeaviSine \\
$\sigma(x)$       & 1       & $x$       & 1               & $x$ \\
$n$ (sample size) 
				          & 200     & 200       & 2048            & 2048 \\
$\M_n$            & regular & 2 bin sizes & dyadic, regular & dyadic, 2 bin sizes\\
\noalign{\smallskip}
\hline
\noalign{\smallskip}
$\E\croch{\penid}$ & $1.919\pm 0.03$ & $2.296\pm 0.05$& $1.028\pm 0.004$& $1.102\pm 0.004$\\
$\E\croch{\penid}+$ & $1.792\pm 0.03$ & $2.028\pm 0.04$& $1.003\pm 0.003$& $1.089\pm 0.004$\\
Mal     & $1.928 \pm 0.04$ & $3.687 \pm 0.07$ & $1.015 \pm 0.003$& $1.373\pm 0.010$ \\
Mal+    & $\meil{1.800 \pm 0.03}$ & $3.173 \pm 0.07$ & $\meil{1.002\pm 0.003}$ & $1.411 \pm 0.008$ \\
\noalign{\smallskip}
\hline
\noalign{\smallskip}
$2-$FCV & $2.078 \pm 0.04$ & $2.542 \pm 0.05$ & $\meil{1.002\pm 0.003}$& $1.184\pm 0.004$ \\
$5-$FCV & $2.137 \pm 0.04$ & $2.582 \pm 0.06$ & $1.014\pm 0.003$& $1.115\pm 0.005$ \\
$10-$FCV& $2.097 \pm 0.05$ & $2.603 \pm 0.06$ & $1.021\pm 0.003$& $1.109\pm 0.004$ \\
$20-$FCV& $2.088 \pm 0.04$ & $2.578 \pm 0.06$ & $1.029\pm 0.004$& $1.105\pm 0.004$ \\
LOO     & $2.077 \pm 0.04$ & $2.593 \pm 0.06$ & $1.034\pm 0.004$& $1.105\pm 0.004$ \\
\noalign{\smallskip}
\hline
\noalign{\smallskip}
pen$ 2-$F& $2.578\pm 0.06$ & $3.061\pm 0.07$ & $1.038\pm 0.004$& $\meil{1.103\pm 0.005}$\\
pen$ 5-$F& $2.219\pm 0.05$ & $2.750\pm 0.06$ & $1.037\pm 0.004$& $1.104\pm 0.004$\\
pen$10-$F& $2.121\pm 0.05$ & $2.653\pm 0.06$ & $1.034\pm 0.004$& $1.104\pm 0.004$\\
pen$20-$F& $2.085\pm 0.04$ & $2.639\pm 0.06$ & $1.034\pm 0.004$& $1.105\pm 0.004$\\
penLoo     & $2.080\pm 0.05$ & $2.593\pm 0.06$ & $1.034\pm 0.004$& $1.105\pm 0.004$\\
\noalign{\smallskip}
\hline
\noalign{\smallskip}
pen$ 2-$F+& $2.175\pm 0.05$ & $2.748\pm 0.06$& $1.011\pm 0.003$& $1.106\pm 0.004$\\
pen$ 5-$F+& $1.913\pm 0.03$ & $2.378\pm 0.05$& $\meil{1.006\pm 0.003}$& $\meil{1.102\pm 0.004}$\\
pen$10-$F+& $1.872\pm 0.03$ & $\meil{2.285\pm 0.05}$& $\meil{1.005\pm 0.003}$& $\meil{1.098\pm 0.004}$\\
pen$20-$F+& $1.898\pm 0.04$ & $\meil{2.254\pm 0.05}$& $\meil{1.004\pm 0.003}$& $\meil{1.098\pm 0.004}$\\
penLoo+     & $\meil{1.844\pm 0.03}$ & $\meil{2.215\pm 0.05}$& $\meil{1.004\pm 0.003}$& $\meil{1.096\pm 0.004}$ \\
\hline
\end{tabular}
\end{center}
\end{table}

\subsection{Results and comments}
First of all, our experiments show the interest of both penVF and VFCV in several difficult framework, with relatively small sample sizes. Although it can not compete with simple procedures such as Mallows' $C_p$ from the computational viewpoint, it is much more efficient when the noise is heteroscedastic (S2 and HSd2). In these hard frameworks, the performances of penVF and VFCV are comparable to those of the ``ideal deterministic penalty'' $\E\croch{\penid}$.
On the other hand, they perform slighlty worse than Mallows' for the easier problems (S1 and HSd1), which we interpretate as the unavoidable price for robustness.

Secondly, in the four experiments, the best procedures are always the overpenalizing ones: many of them even beat the perfectly unbiased $\E\croch{\penid}$, showing the crucial need to overpenalize. This is mainly due to the small sample size compared to the high noise-level, since it is no the case when $\sigma$ is smaller, and less obvious when $n$ is larger (see respectively experiments S0.1 and S1000 in Chap.~5 of \cite{Arl:2007:phd}). We would like to insist on the importance of this phenomenon, which is seldom mentioned because it it vanishes in the asymptotic framework, and it is quite hard to find from theoretical results. 

\medskip

We can now come back to the discussion of Sect.~\ref{VFCV.sec.subopt.VF.choix} on the choice of $V$ for VFCV, which is enlightened by the results of Tab.~\ref{VFCV.tab.un}.
In the first three experiments, and more clearly in HSd1, $V=2$ has comparable or better performances than $V \in \set{5, 10, 20, n}$. This is highly non intuitive, unless we consider the need for overpenalization in those experiments where the signal-to-noise ratio is quite low. It appears that the variability issue is less important in those three cases. This is not because the variance of $\critVF$ is negligible in front of its bias, but mainly because its dependence on $V$ is only mild. Hence, whatever $V$, it has to be compensate by overpenalizing.
On the contrary, the best choices are $V=20$ and $V=n$ in experiment HSd2, where overpenalization seems to be less needed. 
The main conclusion here should be that one really has to take into account both overpenalization and variance for choosing an optimal $V$. 
The larger $V$ is not always the better one, so that a larger computation time does not always improve the accuracy. The main difficulty here is that it does not seem straightforward to choose $V$ from the data only.

\medskip

Finally, let us compare the performances of $V$-fold cross-validation and $V$-fold penalization in Tab.~\ref{VFCV.tab.un}. At first glance, it seems that penVF with $V<20$ performs worse than VFCV in the first three experiments, and not clearly better in the last one. The point is that it matches exactly with the experiments for which overpenalization is crucial. 
But looking at the performance of penVF+, we have evidence for the advantage conferred to penVF by its flexibility. In three over four experiments, penVF+ with any $V \in \set{5, 10, 20, n}$ does better than VFCV with any choice of $V$; and it is almost the case for HSd1. 
This comes from the overpenalizing ability of $V$-fold penalization, which is crucial in such non-asymptotic situations.

Moreover, choosing the optimal $V$ for penVF or penVF+ is much simpler than for VFCV: it is always the largest $V$. Remark that $V=n$ does not always perform significantly better than $V=20$ or $V=10$, which can be considered as almost optimal choices.
For the practical user, the choice of $V$ thus reduces to a trade-off between computational complexity and performance (the latter being governed by the variability of the $V$-fold penalties). Then, once $V$ is chosen, $C$ has to be taken equal to $(V-1)$ times the overpenalization factor (and estimating it from the data remains an open question). 

\medskip

We conclude this section by some additional remarks, concerning some particular points of our simulation study.
\begin{itemize}
\item We also performed Mallows' $C_p$ (and its overpenalized version Mal+) with the true mean variance $\E\croch{\sigma^2(X)}$ instead of  $\widehat{\sigma}^2$ (which would not be possible on a real data set). It gave worse performance for all experiments but S2, in which $C_{\mathrm{or}}(\mathrm{Mal})=2.657 \pm 0.06$ and $C_{\mathrm{or}}(\mathrm{Mal}+)=2.437 \pm 0.05$. This shows that overpenalization is really crucial in experiment S2, even more than the shape of the penalty itself. But once we overpenalize, penVF+ remains significantly better than Mallows' $C_p$ ($\critVF$ being too variable for small $V$ to do better than Mallows). The ability to overpenalize with penVF while keeping the variability low (\latin{i.e.} $V$ large) thus appears to be crucial in this case.
In addition, it can be proved that Mallows' $C_p$ penalty (and, more generally, any penalty of the form $\widehat{K} D_m$) leads to suboptimal model selection in some heteroscedastic framework. See \cite{Arl:2007:phd}, Chap.~4. This should be compared to Thm.~\ref{VFCV.the.histos.penVFCV}, which can be applied in that framework.
\item In experiment HSd1, 2-fold cross-validation appears to be among the best model selection procedures overall. This should be linked with the fact that $\M_n$ only consists on histograms on dyadic partitions of $[0,1]$, so that the assumptions of Thm.~\ref{VFCV.the.subopt} are not fulfilled. More precisely, our computations may show that the model which minimize $\E\croch{\critVF(m)}$ with $V=2$ is the oracle model for arbitrarily large values of $n$. 
This emphasizes the fact that VFCV is not universally suboptimal for model selection for prediction. It is only unable to make the right choice among estimators whose excess losses are within a constant factor smaller than some $K(V)>1$.
\item Eight additional experiments are reported in Chap.~5 of \cite{Arl:2007:phd}, showing similar results with various $n$, $\sigma$ and $\bayes$ (the assumptions of Thm.~\ref{VFCV.the.histos.penVFCV} not being always satisfied). Notice that overpenalization is not always necessary, in particular when the signal-to-noise ratio is larger. In such situations, $V=20$ or $V=n$ is generally optimal for VFCV.
\end{itemize}

\section{Discussion} \label{VFCV.sec.disc}
\subsection{$V$-fold cross-validation \latin{vs.} $V$-fold penalties}
\label{VFCV.sec.disc.VFCV-vs-penVF}
Time has come for us to give an accurate answer to this practical (but quite  hard) question: how to use $V$-fold?

Firstly, the classical $V$-fold cross-validation is biased and asymptotically suboptimal for prediction in some ``easy framework'' (\latin{i.e.} with a smooth regression function and an homoscedastic Gaussian noise). It thus has to be corrected, and we suggest a $V$-fold penalization algorithm that provides such a correction. This algorithm is asymptotically optimal in theory, quite efficient on some simulated data, and has the same computational cost as VFCV.

Secondly, a non-asymptotic phenomenon is likely to arise, that make the problem harder: when the sample size is small and the noise-level large, overpenalizing procedures are more efficient than unbiased ones. Then, our $V$-fold penalization method allows to choose an overpenalizing factor, whereas VFCV imposes it (through $V$) and a corrected VFCV forbids it. This flexibility is the main reason why we suggest to use penVF instead of VFCV or Burman's corrected VFCV. Otherwise, $V$ has to be chosen very carefully, taking into account variability, bias and the possible need for some bias.

\medskip

We shall now explain how to use $V$-fold penalties. It depends on two tuning parameters: the number $V$ of folds and the overpenalization factor $C/(V-1)$. The choice of $V$ depends on the trade-off between variability and computational complexity. If the latter one does not matter, the optimal choice is close to $V=n$ (at least for least-squares regression). Otherwise, the choice has to be done by the final user. We refer to asymptotic computations of Burman \cite{Bur:1989,Bur:1990} (in linear regression) and the recent work of Celisse and Robin \cite{Cel_Rob:2006} (in density estimation) for quantitative measures of variability according to $V$. Further research in that direction would be very useful for practical use of $V$-fold model selection criteria.

The question of choosing the overpenalization factor is probably harder to solve. According to our simulation study, the optimal one depends at least on the sample size, the noise level and the smoothness of the regression function. Since the first criterion is that the penalty almost never underestimates the ideal one, a wise choice of $C$ depends on the fluctuations of both the $V$-fold penalty and the ideal penalty. We thus need a better understanding of the variability of penVF. 
Another idea would be to replace the conditional expectation in \eqref{VFCV.def.penid} by a quantile, in order to build a simultaneous confidence region for the prediction errors $\paren{P\gamma\paren{\ERM_m}}_{\mM_n}$. Then, we could deduce a confidence set, to which the oracle model should belong. Defining $\mh$ as the more parcimonious model in this confidence set, we would have done the work of overpenalization by choosing the probability coverage of the confidence region. 
We refer to \cite{Arl:2007:phd} (Sect.~6.6 and 11.3.3) for further discussions about overpenalization.

\subsection{Other cross-validation methods} \label{VFCV.sec.disc.otherCV}
In this paper, we focused on VFCV and penVF, among many other cross-validation like methods: hold-out, repeated learning-testing methods \cite{Bre_etal:1984}, leave-$p$-out, \latin{etc.} However, it follows from our proofs that the asymptotic performances of these methods mainly depends on their bias, which is itself a function of the ratio between the size of the learning set and the sample size. It is thus possible to have asymptotic optimality with any complexity cost, even without using penVF.

Let us fix for instance the computational complexity to the one of $2$-fold cross-validation. We may use $2$-fold cross-validation, Burman's corrected $2$-fold CV, $2$-fold penalization or repeated learning-testing methods (with 2 splits of the data and a learning set of size equivalent to the sample size $n$). Asymptotically, the first one is suboptimal (Thm.~\ref{VFCV.the.subopt}), while the three other ones are optimal (Thm~\ref{VFCV.the.histos.penVFCV} and the proof of Thm.~\ref{VFCV.the.subopt}). We have already seen in Sect.~\ref{VFCV.sec.disc.VFCV-vs-penVF} that Burman's corrected $2$-fold can not overpenalize when needed, which can be a serious drawback in non-asymptotic situations. 
Repeated learning-testing does not have this drawback, since it is possible to overpenalize within any factor $C\geq 1$ by choosing a learning set of size $\sim n/(2C-1)$. 

However, there remains a strong argument in favour of $2$-fold penalization. When $C$ has to be taken close to 1 (which is the asymptotic situation), repeating learning-testing requires the size of the learning set to be very close to $n$. Hence, if we can only make two splits, most of the data remains in both learning sets. This makes the final criterion much variable, since it strongly depends on the few data which belong to the union of the two training sets. On the contrary, with $2$-fold penalization (as well as $2$-fold cross-validation and its corrected version), each data point belongs is used once for learning and once for training.

Finally, it seems to us that $V$-fold penalization should be preferred, because of its versatility: it is asymptotically optimal, quite flexible (for non-asymptotic situations) and makes use of all the data for both learning and training.

\subsection{Prediction in other frameworks} \label{VFCV.sec.disc.pred}
In order to make theoretical computations feasible, we restricted ourselves to the histogram regression framework in this article. Of course, this is only a first step towards a more general study of $V$-fold methods for model selection. 
Although all our proofs strongly rely on some particular features of histograms (in particular for computing expectations), we conjecture than most of our conclusions stay valid much more generally. 
The main argument supporting this claim is that part of our concentration inequalities are still valid in a general framework, including bounded regression and binary classification. Accurate statements and proofs are to be found in Chap.~7 of \cite{Arl:2007:phd}.
In addition, penVF is built upon the same general heuristics as VFCV, and was never designed particularly for the heteroscedastic histogram regression problem. 
Hence, it should have at least the same robustness and adaptivity properties as VFCV, while its flexibility should allow better performance in terms of multiplicative constants (which may be crucial, when the sample size is small).

\medskip

Let us now point out some expected changes in our analysis in the general case. 
First, the no-overpenalization constant $\CWinf$ may not stay equal to $V-1$. Although me mentioned an asymptotic theoretical argument, it may break down when one considers models with a large number of parameters (that is, dependent from $n$). If this occurs, we suggest to use a data-dependent procedure for estimating $\CWinf$, based upon the so-called ``slope heuristics'' \cite{Bir_Mas:2006,Arl_Mas:2008}. Basically, it states that $\CWinf$ is twice the constant under which $D_{\mh}$ blows up dramatically. We refer to the above papers for a detailed statement of this algorithm, as well as theoretical insights.

Second, the influence of $V$ on variability may also be quite different. For instance, in classification, it is often noticed that the leave-one-out is much more variable than VFCV with smaller values of $V$ \cite{Has_Tib_Fri:2001}. According to Molinaro, Simon and Pfeiffer \cite{Mol_Sim_Pfe:2005}, this seems to disappear when the algorithm producing $\ERM_m$ is stable. In addition, in the density estimation framework, Celisse and Robin \cite{Cel_Rob:2006} also report that the variance of $\critVF$ increases for large $V$. We believe that an extensive study of this variability issue in all those frameworks should be made, considering that it is a crucial point for choosing $V$ for VFCV. It would also be quite interesting to determine whether the variability of penVF depends on $V$ in the same way or not.

\subsection{Consistency} \label{VFCV.sec.disc.consist}
We focused in this article on prediction, but one often uses model selection for identification. In this framework, one assumes that $\bayes \in S_{\mo}$ (and maybe also to some more complex models), and the goal of a model selection procedure is to catch $\mo$ as often as possible, whatever the prediction risk of $\ERM_{\mo}$. Asymptotic optimality there become consistency, \latin{i.e.} 
\[ \P \paren{\mh = \mo} \xrightarrow[n \rightarrow \infty]{} 1 \enspace . \]
There is a huge amount of papers about model selection for identification; we refer to the introduction of papers by Yang \cite{Yan:2006,Yan:2007b} for references about the consistency of cross-validation in the regression and classification settings. 

The main point for consistency is that overpenalization is needed, even from the asymptotic viewpoint. This is the main reason why BIC is roughly the AIC criterion multiplied by a constant times $\ln(n)$. See also Aerts, Claeskens and Hart \cite{Aer_Cla_Har:1999} about this question.
Our penalization interpretation of VFCV (and more generally, any cross-validation like method) then enlightens several theoretical and empirical results about the consistency issue.  

With VFCV, the overpenalization factor is bounded from above by $3/2$ (which corresponds to $V=2$). Hence, $V$-fold cross-validation may be inconsistent in general for any $V$ (although it can sometimes be used, when one compares sufficiently different models, see Yang \cite{Yan:2007b}). Moreover, the better choice is often $V=2$ as remarked by Zhang \cite{Zha:1993}, Dietterich \cite{Die:1998} and Alpaydin \cite{Alp:1999}. 
On the contrary, $V$-fold penalties could work, by choosing $C \propto (V-1) \ln(n)$ (for instance). We conjecture that such a method would be consistent, whatever $V$.

More generally, it has been noticed several times that the consistency of cross-validation requires the size of the learning set to be chosen negligible in front of the sample size. In the linear regression framework, this has be shown by Shao \cite{Sha:1993,Sha:1997}. In the classification setting, this is called the ``cross-validation paradox'' by Yang \cite{Yan:2006}. 
With penVF, we believe that we may have proposed a way of solving this paradox, by allowing to choose the overpenalization factor independently from the size of the learning set.
\appendix

\section{Probabilistic tools} \label{VFCV.sec.proba}
In this section, we give some probability theory results that we need to prove our main result, while being of self-interest.
In the rest of the paper, for any $a,b \in \R$, we denote by $\mini{a}{b}$ the minimum of $a$ and $b$, and by $\maxi{a}{b}$ the maximum of $a$ and $b$. 
\subsection{Expectations of inverses of binomials} \label{VFCV.sec.proba.Einv}
For any non-negative random variable $Z$, define
\[ \einv{Z} = \einv{\mathcal{L}(Z)} \egaldef \E \croch{Z} \E \croch{ Z^{-1} \sachant Z>0} \enspace . \]
Non-asymptotic bounds on this quantity when $Z$ has a binomial distribution are required in the proof of Prop.~\ref{VFCV.pro.EcritVFCV-EcritID}, which is at the core of our main results. Former results concerning $\einv{Z}$ can be found in papers by Lew \cite{1976:Lew} (for general $Z$) or Znidaric \cite{2005:Zni} (for the binomial case), but they are either asymptotic or not accurate enough. The following lemma solves this issue.
\begin{lemma} \label{VFCV.le.einv.binom}
For any $n \in \N \backslash\set{0}$ and $p \in (0;1]$, $\mathcal{B}(n,p)$ denotes the binomial distribution with parameters $(n,p)$, $\kappa_3 = 5.1$ and $\kappa_4 = 3.2$. Then, if $np \geq 1$,
\begin{equation} \label{VFCV.eq.einv-binom.inf-sup}  \mini{ \kappa_4 }{ \paren{ 1+\kappa_3 (np)^{-1/4} } } \geq \einv{\mathcal{B}(n,p)} \geq 1 - e^{-np}  \enspace . \end{equation}
\end{lemma}
In particular, $\einv{\mathcal{B}(n,p)} \rightarrow 1$ when $np \rightarrow \infty$, which can be derived from \cite{2005:Zni}.

\subsection{Concentration of inverses of multinomials} \label{VFCV.sec.app.conc.inv-mult}
Let $(X_{\lambda})_{\lamm} \sim \mathcal{M} (n; (\pl)_{\lamm})$ be a multinomial random vector, $(a_{\lambda})_{\lamm}$ a family of non-negative real numbers, and define for every $T \in (0,1]$
\[ Z_{m,T} \egaldef \sum_{\lamm} a_{\lambda} \min\paren{ T , X_{\lambda}^{-1} } \enspace .\]
Such a quantity naturally appears in our setting, mainly because of the randomness of the design.
Unfortunately, classical concentration inequalities for sums of random variables can not be applied to $Z_{m,T}$ because the $X_{\lambda}$ are not independent. 
Using that they are negatively associated \cite{Joa:1983}, we can use the Cram\'er-Chernoff method \cite{Dub_Ran:1998} to obtain the following lemma. Its complete proof can be found in Sect.~8.8 of \cite{Arl:2007:phd}.

\begin{lemma} \label{RP.le.conc.inv.mult}
Assume that $\min_{\lamm} \set{n \pl} \geq B_n \geq 1$ and $T \in (0,1]$. Define $c_1=0.184$, $c_2 = 0.28$, $c_3 = 9.6$, $c_4 = 0.09$, $c_5 = 10.5$, and for every $t \geq 0$, $\varphi_1(t) = \max(t,1) e^{-\max(t,1)}$.
\begin{enumerate}
\item Lower deviations: for every $x \geq 0$, with probability at least $1 - e^{-x}$, {\small
\begin{equation}
\label{RP.eq:conc:som_inv_multin:inf}
\E\croch{Z_{m,1}} - Z_{m,1} 
\leq \frac{\varphi_1(c_1 B_n)}{c_1} \sum_{\lamm} \frac{a_{\lambda}}{n\pl} 
\quad + 3 \sqrt{2} \sqrt{\sum_{\lamm} \frac{a_{\lambda}^2}{(n\pl)^2}}  \sqrt{4 D_m \exp(- c_1 B_n)+x}
\end{equation}}
\item Upper deviations: for every $x \geq 0$, with probability at least $1 - e^{-x}$, {\small
\begin{equation}
\label{RP.eq:conc:som_inv_multin:sup}
\begin{split}
Z_{m,T} - \E\croch{Z_{m,T}} &\leq \frac{\varphi_1 \left(c_2 B_n  \right)}{c_2} \sum_{\lamm} \left( \frac{a_{\lambda}}{n \pl  } \right) 
\\ &\hspace{-3cm} +   \sqrt{ \sum_{\lamm} \left( \frac {a_{\lambda} } {n \pl} \right)^2 \left(D_m  e^{- c_4 B_n } + x \right)} 
\times  \maxi{c_3}{ \croch{  \frac{ c_5  T \sqrt{x+ e^{ - c_4 B_n }} } {n \min_{\lamm} \set{ \frac{\pl}{a_{\lambda}} } \sqrt{ \sum_{\lamm} \paren{\frac{a_{\lambda}}{n\pl}}^2} }} }
\end{split}
\end{equation} }
\end{enumerate}
\end{lemma}

\subsection{Moment inequalities for some U-statistics} \label{VFCV.sec.proba.conc-Ustats}
There are several papers about concentration or moment inequalities for U-statistics, \latin{e.g.} \cite{Gin_Lat_Zin:2000,Ada:2005}. 
It appears that our main results strongly rely on concentration properties for a particular kind of U-statistics of order 2, which are given by the following lemma. It can be derived either from the aforementioned papers, or from \cite{Bou_Bou_Lug_Mas:2005}, as we did in Sect.~8.9 of \cite{Arl:2007:phd}. 

\begin{lemma} \label{RP.le.mom.U-stat}
Let $(a_{\lambda})_{\lamm}$ and $(b_{\lambda})_{\lamm}$ be two families 
of real numbers, $(r_\lambda)_{\lamm}$ a family of integers. For all $\lamm$, let $(\xi_{\lambda,i})_{1 \leq i \leq
r_{\lambda}}$ be independent centered random variables admitting $2q$-th moments
$m_{2q,\lambda,i}$ for some $q \geq 2$. We define $S_{\lambda,1}$, $S_{\lambda,2}$ and $Z$ as follows:
\begin{equation} \label{RP.eq.U-stat2.pen} 
Z = \sum_{\lamm} \left(a_{\lambda} S_{\lambda,2} + b_{\lambda} S_{\lambda,1}^2 \right)
\quad \mbox{with} \quad S_{\lambda,1} = \sum_{i=1}^{r_{\lambda}} \xi_{\lambda,i} \mbox{ and } S_{\lambda,2} = \sum_{i=1}^{r_{\lambda}} \xi_{\lambda,i}^2  \enspace .\end{equation} 
Then, there is a numerical constant $\kappa \leq 1.271 $ such that, for every $q \geq 2$,
{ \small
\begin{align*} \norm{Z - \E[Z]}_q &\leq
4 \sqrt{\kappa} \sqrt{q} \sqrt{ \sum_{\lamm}
\paren{ \carre{ a_{\lambda} + b_{\lambda} } \sum_{i=1}^{r_{\lambda}} m_{2q,\lambda,i}^4 } } 
 +  8 \sqrt{2} \kappa q \sqrt{
\sum_{\lamm} \left( b_{\lambda}^2 \sum_{1 \leq i\neq j \leq r_{\lambda}} m_{2q,\lambda,i}^2 m_{2q,\lambda,j}^2 \right) } \enspace . \end{align*}
 }
\end{lemma}

\section{Proofs} \label{VFCV.sec.proofs} 

\subsection{Notations} \label{VFCV.sec.proofs.nota}
Before starting the proofs, we introduce some notations or conventions:
\begin{itemize}
\item The letter $L$ will be used to design ``some positive numerical constant, possibly different from some place to another''. 
In the same way, a constant which depends on $c_1, \ldots, c_k$ will be denoted $L_{c_1, \ldots, c_k}$, and if \hyp{A} denotes a set of assumptions, $L_{\hyp{A}}$ will be any constant that depends on the parameters appearing in \hyp{A}.
\item For any non-negative random variable $Z$, we define $\einvz{\mathcal{L}(Z)} \egaldef \E\croch{Z} \E \croch{Z^{-1} \1_{Z>0}} $.
\item For every model $\mM_n$, and every $j \in \set{1, \ldots, V}$,
\begin{align*}
 p_1(m) &\egaldef P \paren{\gamma(\ERM_m) - \gamma(\bayes_m)} &\qquad&   &p_2(m) &\egaldef P_n \paren{\gamma(\bayes_m) - \gamma(\ERM_m)} \\
 p_1^{(-j)}(m) &\egaldef P \paren{\gamma(\ERM_m^{(-j)}) - \gamma(\bayes_m)} &\qquad&  &p_2^{(-j)}(m) &\egaldef P_n^{(-j)} \paren{\gamma(\bayes_m) - \gamma(\ERM_m^{(-j)})} \\
\delc(m) &\egaldef (P_n - P) \paren{ \gamma(\bayes_m) - \gamma(\bayes) } &\qquad& &\delc^{(j)} (m) &\egaldef (P_n^{(j)} - P) \paren{ \gamma(\ERM_m^{(-j)}) - \gamma(\bayes) } \enspace .
\end{align*}
\item Histograms-specific notations: for any random variable $Z$, $q>0$, $\mM_n$ and $\lamm$:
\begin{align*}
\El\croch{Z} &\egaldef \E \croch{Z \sachant \paren{\1_{X_i \in \Il}}_{1 \leq i \leq n, \, \lamm} } 
\qquad \qquad
\norm{Z}_{q}^{(\Lambda_m)} \egaldef \El \croch{ \absj{Z}^q }^{1/q} 
\\ 
S_{\lambda,1} &\egaldef \sum_{X_i \in \Il} \paren{Y_i - \betl} \qquad \mbox{and} \qquad  S_{\lambda,2} \egaldef \sum_{X_i \in \Il} \carre{Y_i - \betl} \enspace . \end{align*}
\item Conventions for $p_1$ and $p_2$ when $\ERM_m$ is not well-defined (in the histogram framework):
{\small
\begin{equation} 
\label{RP.def.punmin}
\punmin (m) = \punzero(m) + \sum_{\lamm}
\pl \carre{\sigl} \1_{\phl=0} \qquad \mbox{with} \qquad \punzero (m) = \sum_{\lamm}
\frac{\pl \1_{\phl>0} }{(n \phl)^2} S_{\lambda,1}^2 \enspace .
\end{equation}
}
\[ \widetilde{p}_2(m) \egaldef p_2(m) + \frac{1}{n} \sum_{\lamm} \carre{\sigl} \1_{n \phl = 0} \]
Notice that whatever the convention we choose (and even if we keep their original definition), $p_1$ and $p_2$ have the same value when $\ERM_m$ is uniquely defined, and we will always remove from $\M_n$ the other models. The choice we make here is only important when writing expectations, so it is merely technical.
In the following, we will often write simply $p_1$ (resp. $p_2$) instead of $\punmin$ (resp. $\widetilde{p}_2$).
\end{itemize}

\subsection{Proof of Thm.~\ref{VFCV.the.subopt}}
The idea of the proof is to show that $\crit_1(m) = P\gamma\paren{\ERM_m}$ and $\crit_2(m) = \critVF(m) - \widehat{c} $ (for some random quantity $\widehat{c}$ independent from $m$) satisfy the assumptions of Lemma~\ref{VFCV.le.determ} below, on an event of large probability. 
To this aim, we will use Prop.~\ref{VFCV.pro.EcritVFCV-EcritID} as well as concentration inequalities of Sect.~\ref{VFCV.sec.proofs.conc}.

First, we have to be more precise about what we do with models $m$ such that $\ERM_m^{(-j)}$ is not well defined for at least one $j \in \set{1, \ldots, V}$. Denote $E_n(m)$ this event.
By \eqref{VFCV.eq.threshold} in Lemma~\ref{le.bernstein.binom}, $E_n(m)$ has a probability smaller than $n^{-2}$ as soon as $D_m \leq L n (\ln(n))^{-1}$, so that all the reasonable conventions will have the same effect. For the sake of simplicity, we choose in this proof is to eliminate such models from $\M_n$. Notice that this removes automatically models such that $\min_{\lamm} \set{n\phl} \leq 1$, in particular all models of dimension strictly larger than $n/2$.

Denote $\widehat{c} = V^{-1} \sum_{j=1}^V P_n^{(j)} \gamma\paren{\bayes}$. Then, for every $\mM_n$,
\begin{equation} \label{eq.critVF.decoupe}
\crit_2(m) \egaldef \critVF(m) - \widehat{c} 
= \perte{\bayes_m} + \frac{1}{V} \sum_{j=1}^V \paren{ p_1^{(-j)}(m) + \delc^{(j)}(m) } + \infty \1_{E_n(m)}  \enspace .\end{equation}
First, notice that for every $j$, conditionally to $(X_i,Y_i)_{i \notin B_j}$, $\ERM_m^{(-j)}$ is deterministic. In addition, $\norm{Y}_{\infty} \leq A \egaldef 1 + \sigma \norm{\epsilon}_{\infty} < \infty$ by assumption.
So, Lemma~\ref{VFCV.le.conc.delta.borne} can be applied  with $t=\ERM_m^{(-j)}$ and $n$ changed into $\card(B_j) \geq L n / V$. 
More precisely, for every $\mM_n$ such that $E_n(m)$ does not hold, for every $j \in \set{1, \ldots, V}$, taking $x = 4 \ln(n)$ and $\eta = \ln(n)^{-1}$, there is an event of probability $1 - L n^{-4}$ on which 
\begin{equation} \label{eq.VFCV-sub.delc} \absj{\delc^{(j)}(m)} \leq \frac{\perte{\ERM_m^{(-j)}}}{\ln(n)} + \frac{L V A^2 \ln(n)^2}{n} \enspace . \end{equation}
A union bound shows that these inequalities hold uniformly over $j$ and $m$ on an event of probability at least $1 - L n^{-2}$.
Combined with \eqref{eq.critVF.decoupe}, this gives 
\begin{equation} \label{eq.VFCV-sub.critVF.min} \crit_2 (m) \geq \paren{1 - \ln(n)^{-1}} \croch{  \perte{\bayes_m} + \frac{1}{V} \sum_{j=1}^V p_1^{(-j)}(m) } - \frac{L V A^2 \ln(n)^2}{n} + \infty \1_{E_n(m)}
\end{equation} and a similar upper bound.

A second key remark is that for every $j$, $p_1^{(-j)}$ has the distribution of $p_1$ with a sample size $n - \card(B_j)$ instead of $n$. We can then apply Prop.~\ref{VFCV.pro.conc.penid} (with $\gamma=4$) to get that on an event of probability $1-Ln^{-2}$, for every $j \in \set{1, \ldots, V}$ and  $\mM_n$ such that $E_n(m)$ does not hold,
\begin{align} 
\label{eq.VFCV-sub.p1.maj}
p_1^{(-j)}(m) &\leq \E \croch{p_1^{(-j)}(m)} + L_{A,\sigma}  \croch{  \ln(n)^2 D_m^{-1/2}  +  \sqrt{D_m} e^{-L n D_m^{-1}} } \E \croch{p_2^{(-j)}(m)} \\
\label{eq.VFCV-sub.p1.min.pt-moy}
p_1^{(-j)}(m) &\geq \E \croch{p_1^{(-j)}(m)} - L_{A,\sigma}  \croch{ \ln(n)^2 D_m^{-1/2} +  e^{-L n D_m^{-1}} } \E \croch{p_2^{(-j)}(m)}
\\
\label{eq.VFCV-sub.p1.min.gros} p_1^{(-j)}(m) &\geq \paren{ L \ln(n)^{-1} - L_{A,\sigma} \ln(n)^2 D_m^{-1} } \E \croch{p_2^{(-j)}(m)} \enspace . 
\end{align}

Finally, since $\bayes(x)=x$, $X$ is uniform and the models are regular histograms on $\X=[0,1]$, we can compute exactly for each model the bias and the variance term (when the sample size is $n$):
\begin{equation} \label{eq.VFCV-sub.bias+var}
\perte{\bayes_m} = \frac{1}{12 D_m^2} \qquad \mbox{and} \qquad \E\croch{p_2(m)} = \frac{\sigma^2 D_m}{n} + \frac{1}{12 D_m n} \enspace .
\end{equation}
\medskip

We now explain how this can be used to check the assumptions of Lemma~\ref{VFCV.le.determ}. 
Let $c_1$ and $\kappa_1$ be positive constants to be chosen later.
\paragraph{Small models}
First, assume that $D_m < \ln(n)^{\kappa_1}$. 
Combining \eqref{eq.VFCV-sub.critVF.min}, \eqref{eq.VFCV-sub.p1.min.pt-moy},  \eqref{eq.VFCV-sub.bias+var} and using that $\E\croch{p_1^{(-j)}(m)} \geq 0$, $\crit_2(m)$ is roughly of the order of the bias term. Hence, condition \eqref{VFCV.le.determ.eq.crit2.petit} holds with $c_3 = L$ and $\kappa_3 = 2 \kappa_1$ when $n \geq L_{A,\sigma,V,\kappa_1}$. Notice that this holds for every $\kappa_1>0$.

\paragraph{Intermediate models}
We now consider models of dimension $\ln(n)^{\kappa_1} \leq D_m \leq c_1 n (\ln(n))^{-1}$. As already noticed, $E_n(m)$ does not hold true for any of them, with a large probability.

From \eqref{eq.VFCV-sub.critVF.min} (and the similar upper bound), \eqref{eq.VFCV-sub.p1.min.pt-moy}, \eqref{eq.VFCV-sub.p1.maj} and \eqref{eq.VFCV-sub.bias+var}, it follows that condition \eqref{VFCV.le.determ.eq.crit2.moyen} holds with $a = 1/12$, $b = \sigma^2$, $C = V/(V-1)$, $c_2 = L_{A,V,\sigma}$ and $\kappa_2 = 1$, as soon as $n \geq L_{A,\sigma,V}$, $c_1 \leq L$ and $\kappa_1 \geq 6$.
Very similar (and somehow simpler) arguments prove that the condition \eqref{VFCV.le.determ.eq.crit1.moyen} holds with the same parameters.

\paragraph{Large models}
Finally, let $\mM_n$ be such that $D_m > c_1 n (\ln(n))^{-1}$.
Combining \eqref{eq.VFCV-sub.critVF.min}, \eqref{eq.VFCV-sub.p1.min.gros} and  \eqref{eq.VFCV-sub.bias+var}, $\crit_2(m)$ is roughly of the order of the variance term $L \E\croch{p_2^{(-j)}(m)}$ when $n \geq L_{A,\sigma,V,c_1}$.  
As a result, condition \eqref{VFCV.le.determ.eq.crit2.grand} holds with $c_4 = L c_1 \sigma^2$ and $\kappa_4 = 2$, for $n \geq L_{A,\sigma,V,c_1}$. 

\medskip

Choosing now $c_1 \leq L$ and $\kappa_1=6$, the conclusion directly follows from  Lemma~\ref{VFCV.le.determ} below. Notice that we have assumed several times that $n \geq n_0 = L_{A,\sigma,V}$. These conditions can be dropped by choosing $K_1 \geq n_0^2$. \qed

\begin{lemma}\label{VFCV.le.determ}
Let $a,b,\paren{c_i}_{1 \leq i \leq 4},\paren{\kappa_i}_{1 \leq i \leq 4},c_{\mathrm{rich}}>0$ and $C>1$ be some constants, $n \in \N$ and $\M_n$ a set of indexes. 
Assume that for every $\mM_n$, $D_m \in [1,n]$, and moreover that $\forall x \in [1, n - c_{\mathrm{rich}}]$, $\exists m \in \M_n$ such that $D_m \in [x,x+c_{\mathrm{rich}}]$. 
Let $\crit_1$ and $\crit_2$ be some functions $\M_n \flens \R$ satisfying the following conditions:
\begin{enumerate}
\item[(i)] for every $\mM_n$, 
\begin{align} \label{VFCV.le.determ.eq.crit1.moyen}
\crit_1(m) &= \paren{ \frac{a}{D_m^{2}} + \frac{b D_m}{n}} \paren{ 1 + \epsilon_{1,m}}  \\ \label{VFCV.le.determ.eq.crit2.moyen}
\crit_2(m) &= \paren{ \frac{a}{D_m^{2}} + \frac{C b D_m}{n}} \paren{ 1 + \epsilon_{2,m}}
\end{align}
with $\max_{i=1,2} \sup_{\mM_n \telque \ln(n)^{\kappa_1} \leq D_m \leq 
\frac{c_1 n}{\ln(n)} } \absj{\epsilon_{i,m}} \leq c_2 \ln(n)^{-\kappa_2} $.
\item[(ii)] for every $\mM_n$ such that $D_m < \ln(n)^{\kappa_1}$, 
\begin{equation} \label{VFCV.le.determ.eq.crit2.petit}
\crit_2(m) \geq c_3 \paren{\ln(n)}^{-\kappa_3} \enspace .
\end{equation}
\item[(iii)] for every $\mM_n$ such that $D_m \geq \frac{c_1 n}{\ln(n)}$, 
\begin{equation} \label{VFCV.le.determ.eq.crit2.grand}
\crit_2(m) \geq c_4 \paren{\ln(n)}^{-\kappa_4} \enspace .
\end{equation}
\end{enumerate}

Then, there is some constant $K(C)=2^{2/3} \times 3^{-1} \paren{ C^{-1/3} - 1}^2 >0$ and some $n_0 >0$ (depending on $a$, $b$, $\paren{c_i}_{1 \leq i \leq 4}$, $\paren{\kappa_i}_{1 \leq i \leq 4}$, $c_{\mathrm{rich}}$ and $C$) such that, if $n \geq n_0$, for every $\mh \in \arg\min_{\mM_n} \crit_2 (m)$, 
\begin{equation} \label{VFCV.le.determ.eq.res.minoration}
\crit_1(\mh) \geq \paren{1 + K(C) - \ln(n)^{-\kappa_2 / 5}} \inf_{\mM_n} \set{\crit_1(m)} \enspace .
\end{equation}
\end{lemma}

\begin{proof}[sketch of the proof of Lemma~\ref{VFCV.le.determ}]
We skip this proof which is only technical. The main arguments are the following. 
First, there is a model $m_1$ of dimension close to $\paren{2 a n}^{1/3} b^{-1/3}$, so that $\crit_1(m_1)$ is close to $3 \times 2^{-2/3} a^{1/3} b^{2/3} n^{-2/3}$.
Second, any model $\mh$ which minimizes $\crit_2(m)$ must have a dimension close to $\paren{2 a n}^{1/3} \paren{b C}^{-1/3}$.
This implies that $\crit(\mh)$ is larger than $(1 + K(C) - \ln(n)^{-\kappa_2 / 5}) \crit_1(m_1)$, and the result follows.
\end{proof}

\subsection{Proof of Thm.~\ref{VFCV.the.histos.penVFCV}} \label{VFCV.sec.proofs.Thm1}
In this section, $L_{\hypVF}$ denotes a constant that depends only on the set of assumptions of Thm.~\ref{VFCV.the.histos.penVFCV}, including $V$. 
For every $\mM_n$, define $\penid^{\prime}(m) = p_1(m) + p_2(m) - \delc(m) = \penid(m) + (P - P_n) \gamma(\bayes)$. Then, by definition of $\penid$ and $\mh$, we have for every $\mMh_n$,
\begin{equation} 
\label{VFCV.eq.oracle.pr.1}
\perte{\ERM_{\mh}} - \paren{\penid^{\prime}(\mh) - \pen(\mh)} \leq \perte{\ERM_m} + \paren{\pen(m) - \penid^{\prime}(m)} \enspace . \end{equation}
The idea of the proof is to show that $\pen - \penid^{\prime}$ is negligible in front of $\perte{\ERM_m}$ for ``reasonable'' models (\latin{i.e.}, those which are likely to be either selected by penVF, or an oracle model) with a large probability. We will prove it by using Prop.~\ref{VFCV.pro.EcritVFCV-EcritID} and~\ref{VFCV.pro.EcritpenVFCV}, as well as the concentration inequalities of Sect.~\ref{VFCV.sec.proofs.conc}.

For every $\mM_n$, define $A_n(m) = \min_{\lamm} \set{ n \phl }$ and $B_n(m) = \min_{\lamm} \set{n\pl}$. 
We now define the event $\Omega_{n}$ on which the concentration inequalities of Prop.~\ref{VFCV.pro.conc.penid} and~\ref{VFCV.pro.conc.penVFCV} and Lemma~\ref{VFCV.le.conc.delta.borne} and~\ref{le.bernstein.binom}, hold with $\gamma  = \aM+2$ (or similarly $x = (\aM+2) \ln(n)$), for every $\mM_n$. Using assumption \hypPpoly, the union bound gives $\Prob\paren{\Omega_{n}} \geq 1 - L_{\cM} n^{-2}$.

\medskip

First, let $c>0$ be a constant to be chosen later, and consider $\widetilde{\M}_n$, the set of models $\mM_n$ such that $\ln(n)^6 \leq D_m \leq c n (\ln(n))^{-1}$. According to \hypArXl, this implies $B_n(m) \geq \crXl c^{-1} \ln(n)$, so that \eqref{VFCV.eq.threshold} ensures that $A_n(m) \geq \ln(n)$ if $c \leq L_{\crXl,\aM}$. In particular, 
$m \in \Mh_n$ on $\Omega_n$.
Now, using both bounds on $D_m$, by construction of $\Omega_n$, 
\[ \max \set{ \absj{\punmin(m) - \E\croch{\punmin(m)}} , \absj{p_2(m) - \E\croch{p_2(m)}} , \absj{\delc(m)} , \absj{\pen(m) - \El\croch{\pen(m)}} } \] is smaller than $L_{\hypVF} \ln(n)^{-1} \paren{ \perte{\bayes_m} + \E\croch{p_2(m)}}$ on this event, at least if $c \leq L_{\crXl}$ (to ensure that $B_n(m)$ is large enough). We now fix $c=L_{\crXl,\aM}$ that satisfies those two conditions.
Using Prop.~\ref{VFCV.pro.EcritpenVFCV}, Lemma~\ref{VFCV.le.calc.p2} and the lower bound on $B_n(m)$, we have for every $m \in \widetilde{\M}_n$
\[ \frac{- L_{\hypVF}}{\ln(n)^{1/4}} \perte{\ERM_m} \leq (\pen-\penid^{\prime})(m) \leq \croch{ 2 (\eta - 1) + \frac{L_{\hypVF}}{\ln(n)^{1/4}} }\perte{\ERM_m} \enspace .\]
as soon as $n \geq L_{\hypVF}$ (this restriction is necessary because the bounds are in terms of excess loss of $\ERM_m$ instead of $\perte{\bayes_m} + \E\croch{p_2}$).
Combined with \eqref{VFCV.eq.oracle.pr.1}, this gives: if $n \geq L_{\hypVF}$ and $c \leq L_{\crXl,\aM}$,
\begin{equation} \label{VFCV.eq.oracle.pr.2} 
\perte{\ERM_{\mh}} \1_{\mh \in \widetilde{\M}_n} \leq \croch{ 2 \eta - 1 + \frac{L_{\hypVF}}{\ln(n)^{1/4}} } \\
\times \inf_{m \in \widetilde{\M}_n} \set{ \perte{\ERM_m} } \enspace .
\end{equation}

\medskip

Second, we prove that any minimizer $\mh$ of $\crit$ belongs to $\widetilde{\M}_n$ on the event $\Omega_n$. Define, for every $\mM_n$, $\crit^{\prime}(m) = \crit(m) - P_n\gamma\paren{\bayes}$, which has the same minimizers over $\Mh_n$ as $\crit$.
According to \hypPrich, there exists $m_0 \in \M_n$ such that $\sqrt{n} \leq D_{m_0} \leq c_{\mathrm{rich}} \sqrt{n}$. If $n \geq L_{\hypVF}$, $m_0 \in \widetilde{\M}_n$, from which we deduce (using \hypAp)
\begin{equation} \label{eq.penVF.critm0} \crit^{\prime}(m_0) \leq \perte{\bayes_{m_0}} + \absj{\delc(m_0)} + \pen(m_0) \leq L_{\hypVF} \paren{n^{-\beta_2/2} + n^{-1/2}} \enspace . \end{equation}
On the other hand, if $D_m < \ln(n)^6$, we have 
\begin{equation} \label{eq.penVF.critm-pt} \crit^{\prime}(m) \geq \perte{\bayes_m} - \absj{\delc(m)} - p_2(m) \geq \cbiasmin \paren{\ln(n)}^{- 6 \beta_1} - L_A \sqrt{\frac{\ln(n)}{n}} - L_{\hypVF} \frac{\ln(n)^7}{n} \end{equation}
on $\Omega_n$. 
In addition, if $D_m > c n (\ln(n))^{-1}$ and $\mMh_n$, by Prop.~\ref{VFCV.pro.EcritpenVFCV}, $\El\croch{\pen(m) - p_2(m)} \geq \El\croch{p_2(m)}$. As a consequence, by construction of $\Omega_n$, we have $\pen(m) - p_2(m) \geq (1 - L_{\hypVF} n^{-1/4}) \E\croch{p_2(m)}$ on it, so that 
\begin{equation} \label{eq.penVF.critm-grd} \crit^{\prime}(m) \geq \pen(m) - p_2(m) - \absj{\delc(m)} \geq L_{\hypVF} \ln(n)^{-1}  \end{equation}
when $n \geq L_{\hypVF}$.
Comparing \eqref{eq.penVF.critm0}, \eqref{eq.penVF.critm-pt} and \eqref{eq.penVF.critm-grd}, it follows that $\mh \in \widetilde{\M}_n$ on $\Omega_n$, provided that $n \geq L_{\hypVF}$.

\medskip

Finally, we show that the infimum can be extended to $\M_n$ in the right-hand side of \eqref{VFCV.eq.oracle.pr.2}, with the convention $\perte{\ERM_m} = + \infty$ if $A_n(m) = 0$. Using similar arguments as above (as well as the definition of $\Omega_n$, in particular \eqref{VFCV.eq.conc.p1.min.2} for large models), we have $\perte{\ERM_{m_0}} \leq L_{\hypVF} \paren{n^{-\beta_2/2} + n^{-1/2}}$ on $\Omega_n$. 
On the other hand, for every $\mM_n$, if $D_m < \ln(n)^6$,  $\perte{\ERM_m} \geq \perte{\bayes_m} \geq L_{\hypVF} \ln(n)^{-6 \beta_1}$ while if $D_m > c n (\ln(n))^{-1}$, $\perte{\ERM_m} \geq L_{\hypVF} \ln(n)^{-2}$ on $\Omega_n$ as soon as $n \geq L_{\hypVF}$. Hence, if $n \geq L_{\hypVF}$, no model $m \notin \widetilde{\M}_n$ can contribute to the infimum in the right-hand side of \eqref{VFCV.eq.oracle.pr.2}. 

\medskip

To conclude the proof of \eqref{VFCV.eq.oracle_traj_non-as}, we notice that $L_{\hypVF} \ln(n)^{-1/4} \leq \epsilon_n = \ln(n)^{-1/5}$ if $n \geq L_{\hypVF}$. All the conditions of the kind $n \geq n_0$ can finally be removed by enlarging $K_1$ so that $K_1 n_0^{-2} \geq 1$.
The final remark concerning $\epsilon_n$ holds true because we can replace the threshold dimensions $\ln(n)^6$ and $c n (\ln(n))^{-1}$ for ``small'' and ``large'' models by some powers of $n$, as soon as the exponents are not taken too far from 0 (resp. 1).

We now get the more classical oracle inequality  \eqref{VFCV.eq.oracle_traj_non-as} by noticing that $\perte{\ERM_m} \leq A^2$ a.s., so that 
\begin{align*}
\E \croch{ \perte{\ERM_{\mh}} } &\leq  \E \croch{ \perte{\ERM_{\mh}} \1_{\Omega_n} } + \norm{ \perte{\ERM_{\mh}}}_{\infty} \Prob\paren{ \Omega_n^c } \\
&\leq \croch{ 2\eta - 1 + \ln(n)^{-1/5} }  \E \croch{ \inf_{\mM_n} \set{ \perte{\ERM_m} } } +  \frac{ A^2 K_1} {n^2} \enspace . \qed \end{align*}

\subsection{Expectations} \label{VFCV.sec.proofs.E}
\subsubsection{Proof of Prop.~\ref{VFCV.pro.EcritVFCV-EcritID}}
\paragraph{Ideal criterion}
We have to compute $\E\croch{P \gamma\paren{\ERM_m} - P \gamma\paren{\bayes_m}} = \E\croch{p_1(m)}$. Assume that $\ERM_m$ is well-defined, \latin{i.e.} $\min_{\lamm} \phl >0$. Using that $\bayes_m$ minimizes $P \gamma(t)$ over $t \in S_m$, we have {\small
\begin{equation} \label{VFCV.eq.p1.hist}
p_1(m) = \sum_{\lamm} \pl \paren{\betl - \bethl}^2 = \sum_{\lamm} \frac{1}{n^2 \phl} \frac{\pl}{\phl} S_{\lambda,1}^2 \quad \mbox{so that} \quad \El \croch{p_1(m)} 
= \frac{1}{n} \sum_{\lamm} \frac{\pl}{\phl} \carre{\sigl} \enspace .
\end{equation} }
The result \eqref{VFCV.eq.EcritID} follows, with $\delta_{n, \pl} = \einvz{\mathcal{B}(n,\pl)} -1$ if $p_1 = \punzero$, or $\delta_{n, \pl} = \einvz{\mathcal{B}(n,\pl)} -1 + n\pl (1-\pl)^n$ if $p_1 = \punmin$.
In each case, the proof of Lemma~\ref{VFCV.le.einv.binom} gives non-asymptotic bounds on $\delta_{n,\pl}$.

\paragraph{$V$-fold criterion}
By definition \eqref{VFCV.def.mh.VFCvclass}, 
on the event on which $\ERM_m^{(-j)}$ is well-defined for every $j$,
\[ \critVF(m) = \frac{1}{V} \sum_{j=1}^V \croch{ p_1^{(-j)}(m) + \paren{ P_n^{(j)} - P} \gamma\paren{ \ERM_m^{(-j)}} + P\gamma\paren{\bayes_m}} \enspace . \]
The second term is centered conditionally to $(X_i,Y_i)_{i \notin B_j}$, so that we only have to compute $\E\croch{p_1^{(-j)}}$ for every $j$. Since $(X_i,Y_i)_{i \notin B_j}$ is an i.i.d. sample of size $n - \card(B_j)$, we can apply the above computation of $\E\croch{p_1}$.
Using a convention similar to $\punzero$ (which can be used on real data, since it does not depend on $P$), the result \eqref{VFCV.eq.EcritVFCV} holds with 
\[ \delta_{n, \pl}^{(VF)} =  \frac{1}{V} \sum_{j=1}^V \croch{  \frac{n-n/V}{n-\card(B_j) } \paren{ \einvz{\mathcal{B}(n-\card(B_j),\pl)} - 1 } + \frac{\card(B_j)}{n - \card(B_j)} - \frac{1}{V-1} }\enspace . \]
From Lemma~\ref{VFCV.le.einv.binom}, we deduce that if $n^{-1} \max_j \card(B_j) \leq c_B < 1$, then 
\[ \frac{-1}{1 - c_B} e^{- n \pl \paren{1 - c_B} } - L \epsilon_n^{reg} \leq \delta_{n, \pl}^{(VF)} \leq \frac{L}{\paren{ 1-c_B }^{5/4}} (n \pl)^{-1/4} +  L \epsilon_n^{reg} \enspace . \qed \]

Similarly to the computation of $p_1(m)$, when $\min_{\lamm} \phl >0$, we have {\small
\begin{equation} \label{VFCV.eq.p2.hist} 
p_2(m) = \sum_{\lamm} \phl \paren{\betl - \bethl}^2 = \sum_{\lamm} \frac{S_{\lambda,1}^2  \1_{n\phl > 0}}{n^2 \phl} \quad \mbox{so that} \quad \El \croch{p_2(m)} = \frac{1}{n} \sum_{\lamm} \carre{\sigl} \enspace . 
\end{equation} }
Notice that  $\El\croch{p_2(m)}=\E\croch{\widetilde{p}_2(m)}$ on this event. Using Lemma~\ref{VFCV.le.einv.binom}, this proves the following.
\begin{lemma}\label{VFCV.le.calc.p2}
If $\min_{\lamm} \set{n \pl} \geq B \geq 1$,
\begin{equation} \label{VFCV.eq.comp.Ep1.Ep2}
\paren{ 1 - e^{-B} } \E \croch{\widetilde{p}_2(m)} \leq \E \croch{\punzero(m)} \leq \E \croch{\punmin(m)} \leq \paren{1 + \sup_{np \geq B} \delta_{n,p} } \E \croch{\widetilde{p}_2(m)}
\end{equation}
where $\delta_{n,p}$ comes from Prop.~\ref{VFCV.pro.EcritVFCV-EcritID}. A similar result holds with $p_2$ instead of $\widetilde{p}_2$ inside the expectation.
\end{lemma}

\subsubsection{Proof of Prop.~\ref{VFCV.pro.EcritpenVFCV}}
%
First of all, notice that all this proof is made {\em conditionally to $\paren{ \1_{X_i \in \Il} }_{ 1 \leq i \leq n, \, \lamm}$}.
The outline of the proof is to prove that $\El\croch{\penVF}$ can be derived from the case where $W$ satisfies an exchangeability condition, for which we can use Lemma~\ref{VFCV.le.cal.pen.Wech} below.
This is why we consider more generally the penalty $\pen_W\paren{m,\paren{X_i,Y_i}_{1 \leq i \leq n}}$, defined by \eqref{VFCV.def.penVFCV.his.2} for a general weight vector $W \in \R^n$, strengthening its dependence on the distribution of $W$ and the data. When $W$ is the subsampling weight vector of interest, $\pen_W$ coincides with the definition of $\penVF$ in Algorithm~\ref{VFCV.def.proc.his}.

Let $\sigma$ be a random permutation of $\set{1, \ldots , n}$, independent from $W$ and the data, and uniform over the permutations that leave $\paren{ \1_{X_i \in \Il} }_{ 1 \leq i \leq n, \, \lamm}$ invariant. 
Defining $\widetilde{W} = \paren{ W_{\sigma(i)} }_{1 \leq i \leq n}$,
\begin{align*} \El \croch{ \pen_{\widetilde{W}} \paren{m,\paren{X_i,Y_i}_{1 \leq i \leq n}} } 
&= \El \croch{ \pen_W \paren{m,\paren{X_{\sigma^{-1}(i)},Y_{\sigma^{-1}(i)}}_{1 \leq i \leq n}} } \\
&= \El \croch{ \pen_W \paren{m,\paren{X_i,Y_i}_{1 \leq i \leq n}} }  \end{align*} since the penalty does not depend on the order of $(W_i,X_i,Y_i)_{X_i \in \Il}$ (for the first equality), and $(X_i,Y_i)_{X_i \in \Il}$ is exchangeable (for the second equality). 
Moreover, for every $\lamm$, $(\widetilde{W}_i)_{X_i \in \Il}$ is exchangeable and independent from $(X_i,Y_i)_{X_i \in \Il}$. 
We can thus use Lemma~\ref{VFCV.le.cal.pen.Wech} to compute $\pen_{\widetilde{W}} (m)$. Then,
\begin{align*}
\El \croch{ \pen(m) } = \frac{C}{n} \sum_{\lamm} \paren{ R_{1,\widetilde{W}}(n,\phl) + R_{2,\widetilde{W}}(n,\phl) } \carre{\sigl} \enspace .
\end{align*}
It now remains to compute $R_{1,\widetilde{W}}$ and $R_{2,\widetilde{W}}$. If $V$ divides $n \phl$, then $\Wl=1$ a.s. and $R_{1,\widetilde{W}} = R_{2,\widetilde{W}} = (V-1)^{-1}$. For the general case, see the proof of Prop.~5.2 in \cite{Arl:2007:phd} (Sect.~5.7.2). \qed

%
\begin{lemma}[Lemma~5.7 of \cite{Arl:2007:phd}]\label{VFCV.le.cal.pen.Wech}
Let $S_m$ be the model of histograms adapted to some partition $\paren{\Il}_{\lamm}$, $W \in [0;\infty)^n$ be a random vector such that for every $\lamm$, $(W_i)_{X_i \in \Il}$ is exchangeable and independent from $(X_i,Y_i)_{X_i \in \Il}$. Define the Resampling Penalty for histograms as \eqref{VFCV.def.penVFCV.his.2}, and assume $\min_{\lamm}\set{n\phl} \geq 2$.
%
Then, 
\begin{gather}
\label{VFCV.eq.pen.Wech}
 \pen(m) = \frac{C}{n} \sum_{\lamm} \paren{R_{1,W}(n,\phl) + R_{2,W}(n,\phl)}  \frac{n\phl S_{\lambda,2} - S_{\lambda,1}^2 }{n\phl - 1} \enspace,  \qquad \mbox{where} \\
\label{VFCV.def.R1}
R_{1,W}(n,\phl) = \El\croch{ \frac{(W_{i_{\lambda}} - \Wl)^2}{\Wl^2} \sachant \Wl >0 } 
\quad 
R_{2,W}(n,\phl) = \El\croch{ \frac{(W_{i_{\lambda}} - \Wl)^2}{\Wl} }  \enspace. 
\end{gather}
and $i_{\lambda}$ is any index such that $X_{i_{\lambda}} \in \Il$.
\end{lemma}

\subsection{Concentration results} \label{VFCV.sec.proofs.conc}
In order to prove Thm.~\ref{VFCV.the.subopt} and~\ref{VFCV.the.histos.penVFCV}, we need to combine Prop.~\ref{VFCV.pro.EcritVFCV-EcritID} and~\ref{VFCV.pro.EcritpenVFCV} with concentration inequalities, which are the purpose of the present section. 
Let $S_m$ be the model of histograms associated with some partition $\paren{\Il}_{\lamm}$, and assume that both \hypAb\ and \hypAn\ are satisfied (see the statement of Thm.~\ref{VFCV.the.histos.penVFCV}).

Our first result has to deal with $p_1$ and $p_2$, which are the main components of the ideal penalty. Whereas concentration for $p_2$ can be obtained in a general framework (see \cite{Arl:2007:phd}, Chap.~7), lower bounds on $p_1$ are completely new, up to our best knowledge.
\begin{proposition} \label{VFCV.pro.conc.penid}
Let $\gamma>0$ and assume that $\min_{\lamm} \set{n \pl} \geq B_n $.
Then, if $B_n \geq 1$, on an event of probability at least $1 - L n^{-\gamma}$, 
\begin{gather} \label{VFCV.eq.conc.p1.min}
\punmin(m) \geq \E \croch{\punmin(m)} - L_{A,\sigmin,\gamma}  \croch{   \ln(n)^2 D_m^{-1/2} +  e^{-L B_n} }  \E\croch{p_2(m)}
\\
\label{VFCV.eq.conc.p1.maj}
\punmin(m) \leq \E \croch{\punmin(m)} + L_{A,\sigmin,\gamma}  \croch{  \ln(n)^2 D_m^{-1/2}  +  \sqrt{D_m} e^{-L B_n} }  \E\croch{p_2(m)}
\\
\label{VFCV.eq.conc.p2}
\absj{p_2(m) - \E [p_2(m)]} \leq L_{A,\sigmin,\gamma} D_m^{-1/2} \ln(n) \E\croch{p_2(m)} \enspace .
\end{gather}
%
%
In addition, if $B_n>0$, there is an event of probability at least $1-Ln^{-\gamma}$ on which
\begin{equation}
\label{VFCV.eq.conc.p1.min.2}
\punmin(m) \geq \paren{ \frac{1}{2 + (\gamma+1) B_n^{-1} \ln(n) } -  L_{A,\sigmin,\gamma} \ln(n)^2 D_m^{-1/2} } \E\croch{\widetilde{p}_2(m)} \enspace .
\end{equation}
\end{proposition}
\begin{proof}[proof of Prop.~\ref{VFCV.pro.conc.penid}]
According to the explicit expressions \eqref{VFCV.eq.p1.hist} and \eqref{VFCV.eq.p2.hist}, $\punmin(m)$ and $p_2(m)$ are both U-statistics of order 2 conditionally to $(\1_{X_i \in \Il})_{(i,\lambda)}$. Then, we use Lemma~\ref{RP.le.mom.U-stat}, with $\xi_{i,\lambda} = Y_i - \betl$, $a_{\lambda} = 0$,  $b_{\lambda} = \pl (n \phl)^{-2}$ for $\punmin$ and $b_{\lambda} = (n^2 \phl)^{-1}$ for $p_2$.
This proves, for all $q \geq 2$, 
\begin{align} \label{RP.eq.mom.p1}
\norm{\punmin(m) - \El [\punmin(m)]}^{(\Lambda_m)}_q &\leq \max_{\lamm} \set{ \frac{\pl}{\phl} \1_{\phl >0}} L_{A,\sigmin} D_m^{-1/2} q \E\croch{p_2(m)} \\ \label{RP.eq.mom.p2}
\norm{p_2(m) - \E [p_2(m)]}^{(\Lambda_m)}_q &\leq L_{A,\sigmin} D_m^{-1/2} q \E\croch{p_2(m)} \enspace .
\end{align}

We deduce conditional concentration inequalities from those moment inequalities (for instance by Lemma~8.9 of \cite{Arl:2007:phd}), with a deterministic probability bound $1 - L e^{-x} = 1 - n^{-\gamma}$. Hence, we deduce unconditional concentration inequalities, and the result follows for $p_2$. 
To control the remainder term for $\punmin$, we use \ref{VFCV.eq.maj.Qmd} in Lemma~\ref{le.bernstein.binom}.

We now have to control the distance between $\El\croch{\punmin}$ and $\E\croch{\punmin}$. First, if $B_n \geq 1$, we can use Lemma~\ref{RP.le.conc.inv.mult}: taking $X_{\lambda} = n\phl$ and $a_{\lambda}= \pl \carre{\sigl}$,  according to \eqref{VFCV.eq.p1.hist}, we have $\punmin(m) = Z_{m,1}$ and the concentration inequality for $\punmin$ follows.
On the other hand, if we only know that $B_n >0$, instead of using Lemma~\ref{RP.le.conc.inv.mult}, we remark that 
\[ \El \croch{\punmin(m)} \geq \min_{\lamm} \set{ \frac{\pl}{\phl} } \El \croch{p_2(m)} \enspace , \]
and the result follows thanks to \eqref{VFCV.eq.min.Qmd}  in Lemma~\ref{le.bernstein.binom}.
\end{proof}

We mention here a much classical result, which is a consequence of Bernstein's inequality, since it deals with sums of independent variables. We refer to \cite{Arl_Mas:2008} for a detailed proof.
\begin{lemma}[Prop.~3, \cite{Arl_Mas:2008}]
 \label{VFCV.le.conc.delta.borne}
Let $t$ be any deterministic predictor. For every $x \geq 0$, there is an event of probability at least $1 - 2 e^{-x}$ on which
\begin{gather} 
 \label{VFCV.eq.conc.delta.borne.2} \forall \eta >0, \quad \absj{(P-P_n)\paren{\gamma\paren{t} - \gamma\paren{s}}} \leq \eta \perte{t} + \left( \frac{4}{\eta} + \frac{8}{3} \right) \frac{A^2 x}{n} 
\enspace .\end{gather}
\end{lemma}

Finally, we consider the $V$-fold penalties defined by Algorithm~\ref{VFCV.def.proc.his}.
\begin{proposition} \label{VFCV.pro.conc.penVFCV}
Let $\pen(m)$ be defined by \eqref{VFCV.def.penVFCV.his} with the weights $W$ defined in Algorithm~\ref{VFCV.def.proc.his} and $\gamma>0$. 
There is an event of probability at least $1 - n^{-\gamma}$ on which, if $\min_{\lamm} \phl >0$,
\begin{equation} \label{VFCV.eq.conc.penVFCV}
\begin{split}
\absj{\pen(m) - \El \croch{\pen(m)} } 
\leq C \paren{\maxi{\frac{1}{\min_{\lamm} \set{n\phl}}}{\frac{1}{V}}} L_{A,\sigmin,\gamma} D_m^{-1/2} \ln(n) \E \croch{p_2(m)} \enspace .
\end{split}
\end{equation}
\end{proposition}
\begin{proof}[proof of Prop.~\ref{VFCV.pro.conc.penVFCV}]
By definition \eqref{VFCV.def.penVFCV.his}, $\pen(m) = \Es \croch{Z}$ with
\begin{align} 
Z &= \sum_{\lamm} \paren{\phl + \phlW} \paren{\bethl - \bethlW}^2 
\label{VFCV.eq.conc.penVFCV.1}
&= \sum_{\lamm} \frac{1 + \Wl}{n^2 \phl \Wl^2} \croch{\sum_{X_i \in \Il} \paren{\Wl - W_i} \paren{Y_i - \betl} }^2 \enspace .
\end{align}

For every $q \geq 1$, using Jensen inequality and the independence between $W$ and the data (conditionally to $\paren{\1_{X_i \in \Il}}_{i,\lambda}$),
\begin{align} \notag \norm{ \pen(m) - \El \croch{\pen(m)} }_q^{(\Lambda_m)} &\leq \norm{ Z - \El \croch{Z \sachant W} }_q^{(\Lambda_m)} \\
&\leq  
\sup_{W_0 \in \supp(W)} \set{ \norm{ Z - \El\croch{Z \sachant W=W_0} }_q^{(W_0,\Lambda_m)} } \label{VFCV.eq.conc.penVFCV.2}
\end{align}
where $\supp(W)$ is the support of the resampling weight vector $W$ distribution (conditionally to $\paren{\1_{X_i \in \Il}}_{i,\lambda}$) and $\norm{\cdot}_q^{(W_0,\Lambda_m)}$ denotes the $q$-th moment conditionally to $\paren{\1_{X_i \in \Il}}_{(i,\lambda)}$ and $W=W_0$.
In other words, the deviations of $\pen$ are smaller than those of the worse case with a deterministic weight vector $W_0 \in \supp(W)$.

From now on, we work conditionally to $\paren{\1_{X_i \in \Il}}_{(i,\lambda)}$ and assume that $W \in \R^n$ is deterministic, among those authorized by Algorithm~\ref{VFCV.def.proc.his}. Denote by $X_{(1,\lambda)}, \ldots, X_{(n \phl,\lambda)}$ the data such that $X_i \in \Il$. 
According to \eqref{VFCV.eq.conc.penVFCV.1}, Lemma~\ref{RP.le.mom.U-stat} with $r_{\lambda} = n \phl$, $a_{\lambda} = 0$, $b_{\lambda} = (1 + \Wl)(n^2 \phl \Wl^2)^{-1}$ and $\xi_{i,\lambda} = \paren{W_{(i,\lambda)} - \Wl} \paren{Y_{(i,\lambda)} - \betl}$ shows that 
\begin{align*}
\norm{ Z - \El\croch{Z} }_q^{(W,\Lambda_m)} 
&\leq 
\frac{L A^2 q}{n} \sqrt{ \sum_{\lamm} \paren{ \frac{1+\Wl}{n \phl \Wl^2} \sum_{i=1}^{n \phl} \paren{W_{(i,\lambda)} - \Wl}^2}^2 } \enspace .
\end{align*}

We now fix some $\lamm$ and write $n \phl = a V + b \geq 1$ with $a,b \in \N$ and $0 \leq b \leq V-1$. Since $W$ is in the support of the $V$-fold weights distribution of Algorithm~\ref{VFCV.def.proc.gal.VFold}, there is an $\epsilon \in \set{0,1}$ such that 
\begin{align*}
\set{ W_i \telque X_i \in \Il } &= \set{ 0 \mbox{ repeated } a + \epsilon \mbox{ times}, \frac{V}{V-1} \mbox{ repeated } r_{\lambda} - a - \epsilon \mbox{ times} } \enspace . \end{align*}
Hence, 
\[ \Wl 
= 1 + \frac{b - V \epsilon} {(V-1) (aV+b)}  
\quad \mbox{and} \quad 
\sum_{i=1}^{r_{\lambda}} \paren{W_{(i,\lambda)} - \Wl}^2 
\leq L \times \croch{ \maxi{1}{\frac{n \phl}{V}} } \enspace , \]
so that for every $q \geq 2$,
\begin{align*}
\norm{ \pen(m) - \El \croch{\pen(m)} }_q^{(\Lambda_m)} 
&\leq L A^2 q \croch{ \maxi { \frac{1}{\min_{\lamm} \set{n\phl} } } {\frac{1}{V}}} \El \croch{p_2(m)} \enspace .
\end{align*}
The classical link between moment and concentration inequalities (\latin{e.g.} Lemma~8.9 in \cite{Arl:2007:phd}) gives \eqref{VFCV.eq.conc.penVFCV} conditionally to $\paren{\1_{X_i \in \Il}}_{i,\lambda}$. We can remove this conditioning since the probability bound $1 - n^{-\gamma}$ is deterministic.
\end{proof}

\subsection{Expectation of inverses of binomials (proof of Lemma~\ref{VFCV.le.einv.binom})}
Let $Z \sim \mathcal{B}(n,p)$. By Jensen inequality,
\[ \einv{Z} \geq \Prob(Z>0) = 1-(1-p)^n \geq 1-e^{-np} \enspace . \] 
For the upper bound, define
\begin{equation} \label{VFCV.def.einfz}
\einvz{\mathcal{L}(Z)} \egaldef \E \croch{Z} \E \croch{ Z^{-1} \1_{Z>0}} = \einv{Z} \Prob(Z>0) \enspace ,
\end{equation} 
so that we can focus on $\einvz{\mathcal{B}(n,p)}$. 

The bound by $\kappa_4$ follows from Lemma~$4.1$ of \cite{Gyo_etal:2002}, according to which \begin{equation} \label{RP.eq.einvz.2} \forall n \in \N, \, \forall p \in [0,1], \quad \einvz{\mathcal{B}(n,p)} \leq \frac{2 n p}{(n+1)p} \leq 2 \enspace .\end{equation}
We can now assume that $np \geq A \geq 29.17$ since otherwise, $1 + \kappa_3 (np)^{-1/4} \geq \kappa_4$. 
Using that $\Prob(1 > Z > 0)=0$, we have for every $\alpha>0$,
\begin{align*}
\einvz{\mathcal{B}(n,p)} = np \E \croch{Z^{-1} \1_{\alpha \E[Z] > Z>0}} + np  \E \croch{Z^{-1} \1_{Z \geq \alpha \E[Z]}}  
\leq  np \Prob \paren{\alpha np > Z>0} + \alpha^{-1} \enspace .
\end{align*}
We now bound the probability on the right-hand side thanks to Bernstein's inequality (\latin{e.g.} Prop.~$2.9$ of \cite{Mas:2003:St-Flour}): 
\[ \forall \theta>0, \quad  \Prob \paren{Z \leq \paren{1 - \sqrt{2 \theta} - \frac{\theta}{3}} np} \leq e^{-\theta np} \enspace ,\]
and $\theta = A^{-1/2}$. Straightforward computations shows that 
\begin{align*} \sup_{np\geq A} \{ \einv{\mathcal{B}(n,p)} \} &\leq \left[ \frac{1}{1 - \sqrt{2} A^{-1/4} -
\frac{1}{3} A^{-1/2}} + A e^{- \sqrt{A}} \right] \frac{1}{1 - e^{-A}} \enspace ,
\end{align*} 
from which the result follows.

%

\subsection{A technical lemma}
Because of the randomness of the design, we have to ensure that the empirical frequencies $n \phl$ are not too far from the expected ones $n \pl$.
\begin{lemma} \label{le.bernstein.binom}
Let $(\pl)_{\lamm}$ be non-negative real numbers of sum 1, $(n\phl)_{\lamm}$ a multinomial vector of parameters
$(n;(\pl)_{\lamm})$, $\gamma>0$. Assume that $\card(\Lambda_m) \leq n$ and $\min_{\lamm} \set{n \pl} \geq B_n >0$. There is an event of probability at least $1 - L n^{-\gamma}$ on which the following three inequalities hold.
\begin{align}
\label{VFCV.eq.maj.Qmd}
\max_{\lamm} \set{ \frac{\pl}{\phl} \1_{\phl >0}}  \leq L \times (\gamma+1) \ln(n) \\
\label{VFCV.eq.min.Qmd}
\min_{\lamm} \set{ \frac{\pl}{\phl}}  \geq \frac{1}{2 + (\gamma+1) B_n^{-1} \ln(n) }  \\
\label{VFCV.eq.threshold} 
\min_{\lamm} \set{n \phl} \geq \frac{\min_{\lamm} \set{n \pl}}{2} - 2 (\gamma+1) \ln(n) 
\end{align}
\end{lemma}
\begin{proof}[proof of Lemma~\ref{le.bernstein.binom}]
Those three results come from Bernstein's inequality (\latin{e.g.} Prop.~$2.9$ of \cite{Mas:2003:St-Flour}) applied to $n \phl$: for every $\lamm$, there is a set of probability $1 - 2n^{- (\gamma+1)}$ on which
\[ n\pl - \sqrt{2 n \pl (\gamma+1) \ln(n)} - \frac{(\gamma+1) \ln(n)}{3} \leq n \phl \leq n \pl + \sqrt{2 n \pl (\gamma+1) \ln(n)} + \frac{(\gamma+1) \ln(n)}{3} \enspace . \]
For \eqref{VFCV.eq.maj.Qmd}, if $n \pl \geq 8 (\gamma+1) \ln(n)$, the lower bound gives the result. Otherwise, remark only that $(\pl/\phl) \1_{\phl >0} \leq n \pl \leq 8 (\gamma+1 \ln(n)$.
For \eqref{VFCV.eq.min.Qmd}, use the upper bound and remark that $n \pl (\gamma+1) \ln(n) B_n^{-1} \geq (\gamma +1) \ln(n) $.
For \eqref{VFCV.eq.threshold}, use the lower bound and remark that $\sqrt{2n \pl (\gamma +1) \ln(n) } \leq (n \pl)/2 + (\gamma +1) \ln(n) $. 
Finally, the union bound gives the result since $\card(\Lambda_m) \leq n$.
\end{proof}

\section*{Acknowledgments}
The author would like to thank gratefully Pascal Massart for several fruitful discussions.

\bibliographystyle{alpha}
\bibliography{penVF,publi_arlot_en}

\end{document}


\begin{frontmatter}

\title{Technical appendix to ``$V$-fold cross-validation improved: $V$-fold penalization''}
\runtitle{Technical appendix}

\begin{aug}
\author{\fnms{Sylvain} \snm{Arlot}\ead[label=e1]{sylvain.arlot@math.u-psud.fr}}
\runauthor{Arlot, S.}

\affiliation{Universit\'e Paris-Sud\vspace{1cm}\\}

\address{Sylvain Arlot\\
Univ Paris-Sud, UMR 8628, \\Laboratoire de Math\'ematiques,\\ Orsay, F-91405 ;
CNRS, Orsay, F-91405 ;\\
INRIA-Futurs, Projet Select \\
\printead{e1}\\
\phantom{E-mail: sylvain.arlot@math.u-psud.fr\ }}
\end{aug}

\runauthor{Arlot, S.}

\begin{abstract}
This is a technical appendix to ``$V$-fold cross-validation improved: $V$-fold penalization''. We present some additional simulation experiments, a few remarks about expectations of inverses, and the proofs which have been skipped or shortened in the main paper.
\end{abstract}

\begin{keyword}[class=AMS]
\kwd[Primary ]{62G09}
\kwd[; secondary ]{62G08}
\kwd{62M20}
\end{keyword}

\begin{keyword}
\kwd{non-parametric statistics}
\kwd{statistical learning}
\kwd{resampling}
\kwd{non-asymptotic}
\kwd{$V$-fold cross-validation}
\kwd{model selection}
\kwd{penalization}
\kwd{non-parametric regression}
\kwd{adaptivity}
\kwd{heteroscedastic data}
\end{keyword}

\end{frontmatter}

Throughout this appendix, we use the notations of the main paper \cite{Arl:2008a}. In order to distinguish references within the appendix from references to the main paper, we denote the former ones by $(1)$ or $1$, and the latter ones by \maineqref{1} or \main{1}.

Following the ordering of \cite{Arl:2008a}, we first present the additional simulation studies mentioned in Sect.~\main{4}. Then, we add a few comments to Appendix~\main{A.1}. Finally, we give some technical proofs.

\section{Simulation study}
We consider in this section eight experiments (called S1000, S$\sqrt{0.1}$, S0.1, Svar2, Sqrt, His6, DopReg and Dop2bin) in which we have compared the same procedures as in Sect.~\main{4}, with the same benchmarks, but with only $N=250$ samples for each experiment.

\begin{figure}
\centerline{\epsfig{file=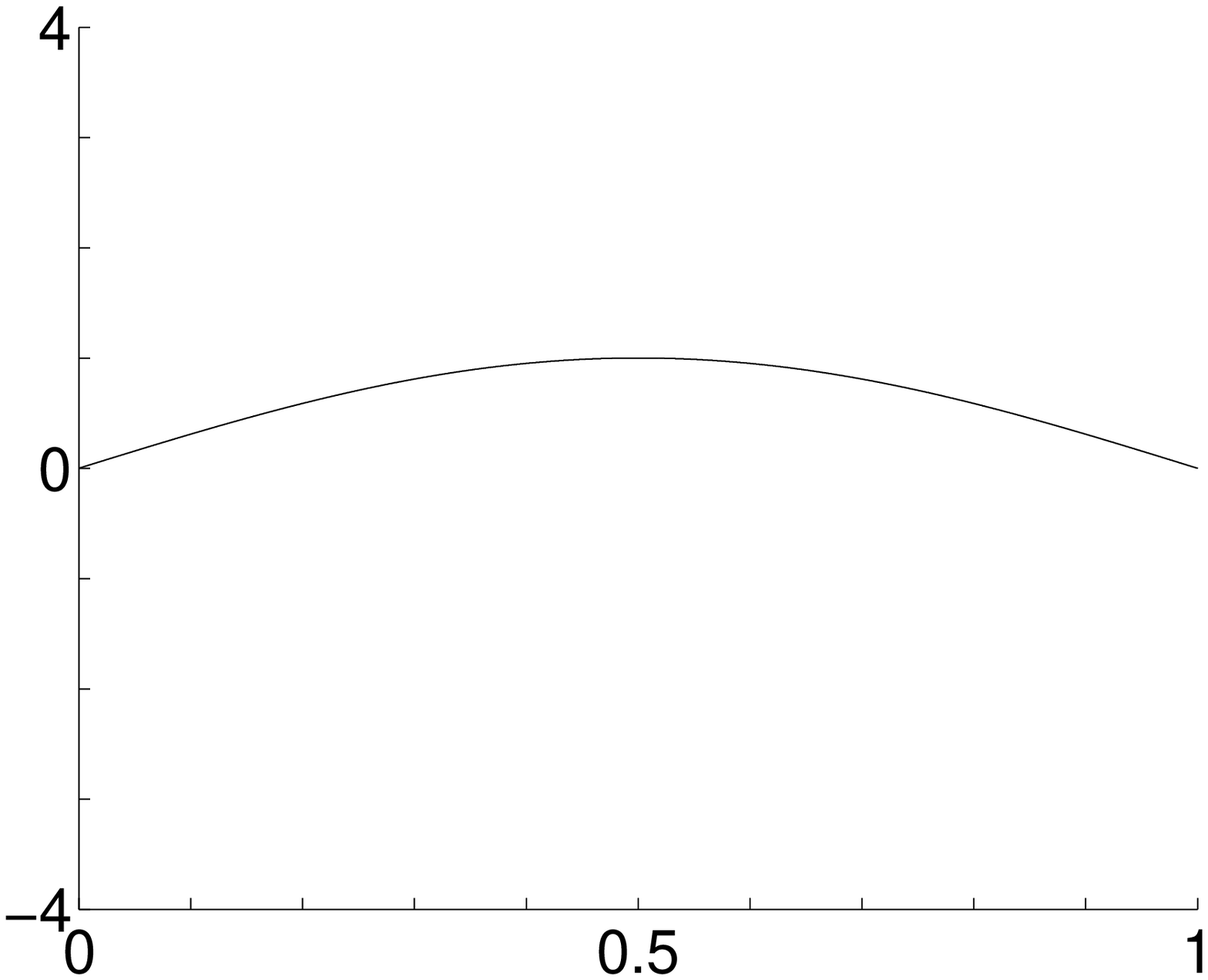,width=0.384\textwidth}}
   \caption{$\bayes(x) = \sin(\pi x)$\label{fig.S.fonc}}
\end{figure}

\begin{figure}
 \begin{minipage}[b]{.48\linewidth}   \centerline{\epsfig{file=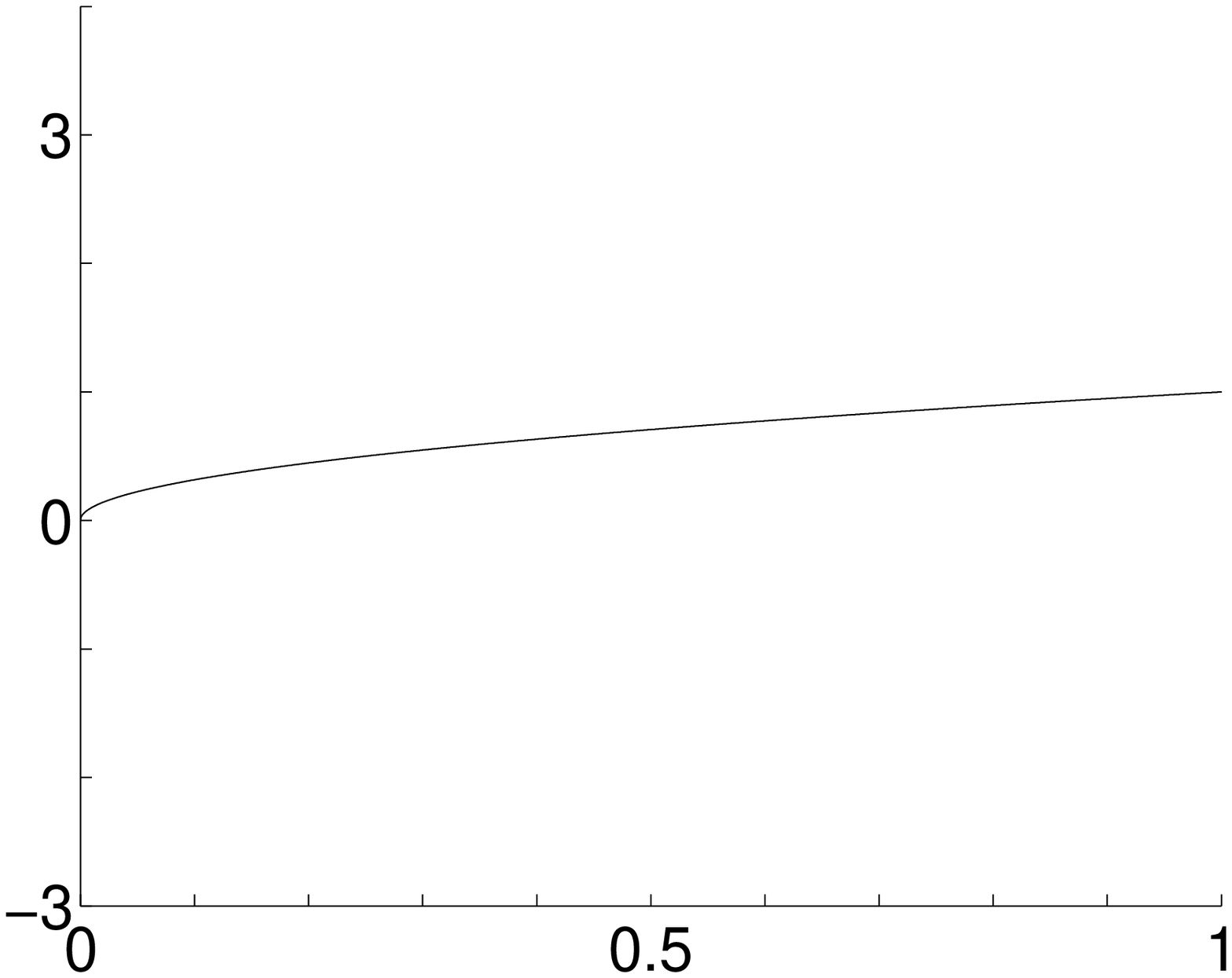,width=0.85\textwidth}}
   \caption{$\bayes(x) = \sqrt{x}$ \label{fig.Sqrt.fonc}}
\end{minipage} \hfill
\begin{minipage}[b]{.48\linewidth}   \centerline{\epsfig{file=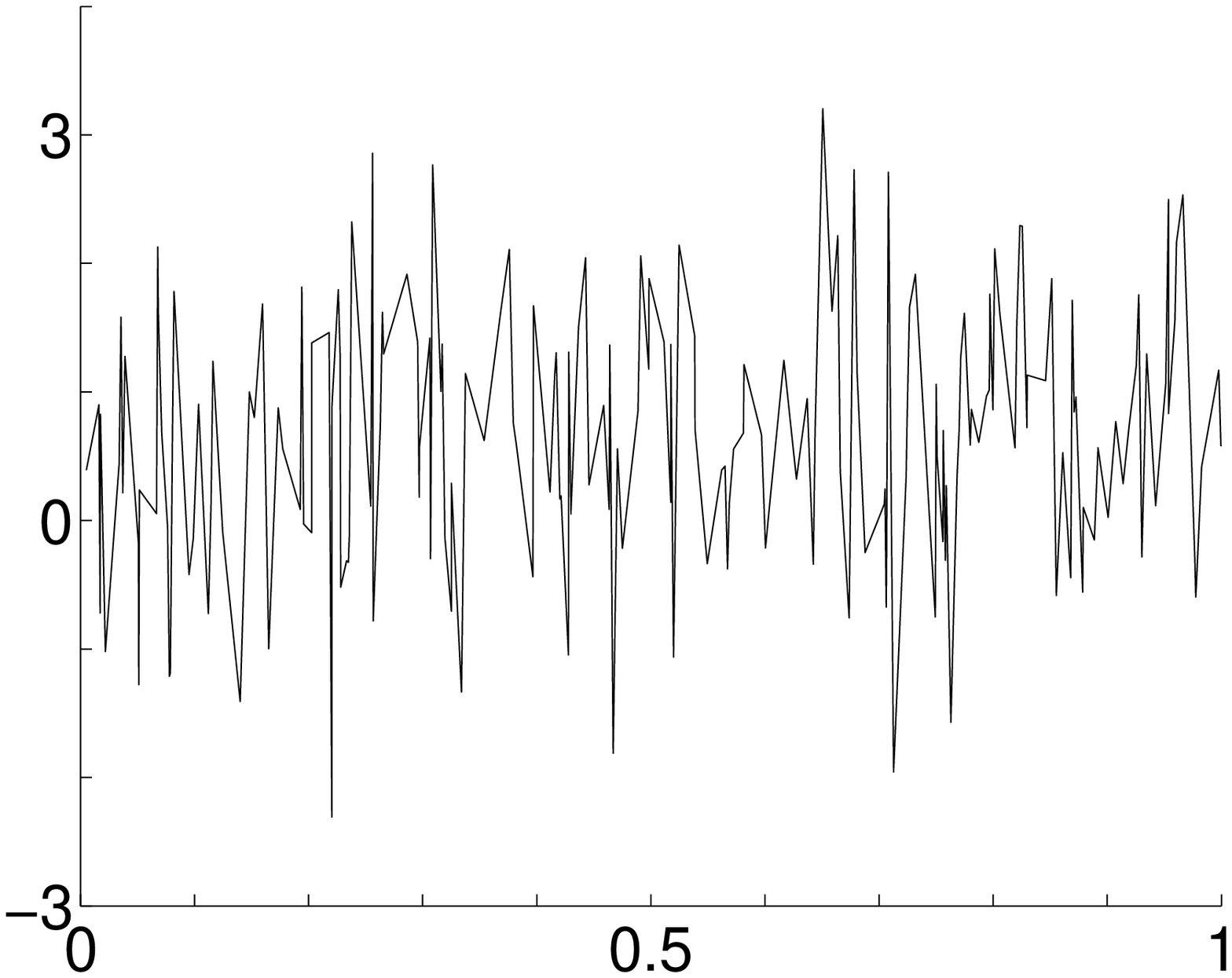,width=0.85\textwidth}}
   \caption{Data sample for Sqrt \label{fig.Sqrt.data}}
\end{minipage} 
\end{figure}

\begin{figure}
\begin{minipage}[b]{.48\linewidth}   \centerline{\epsfig{file=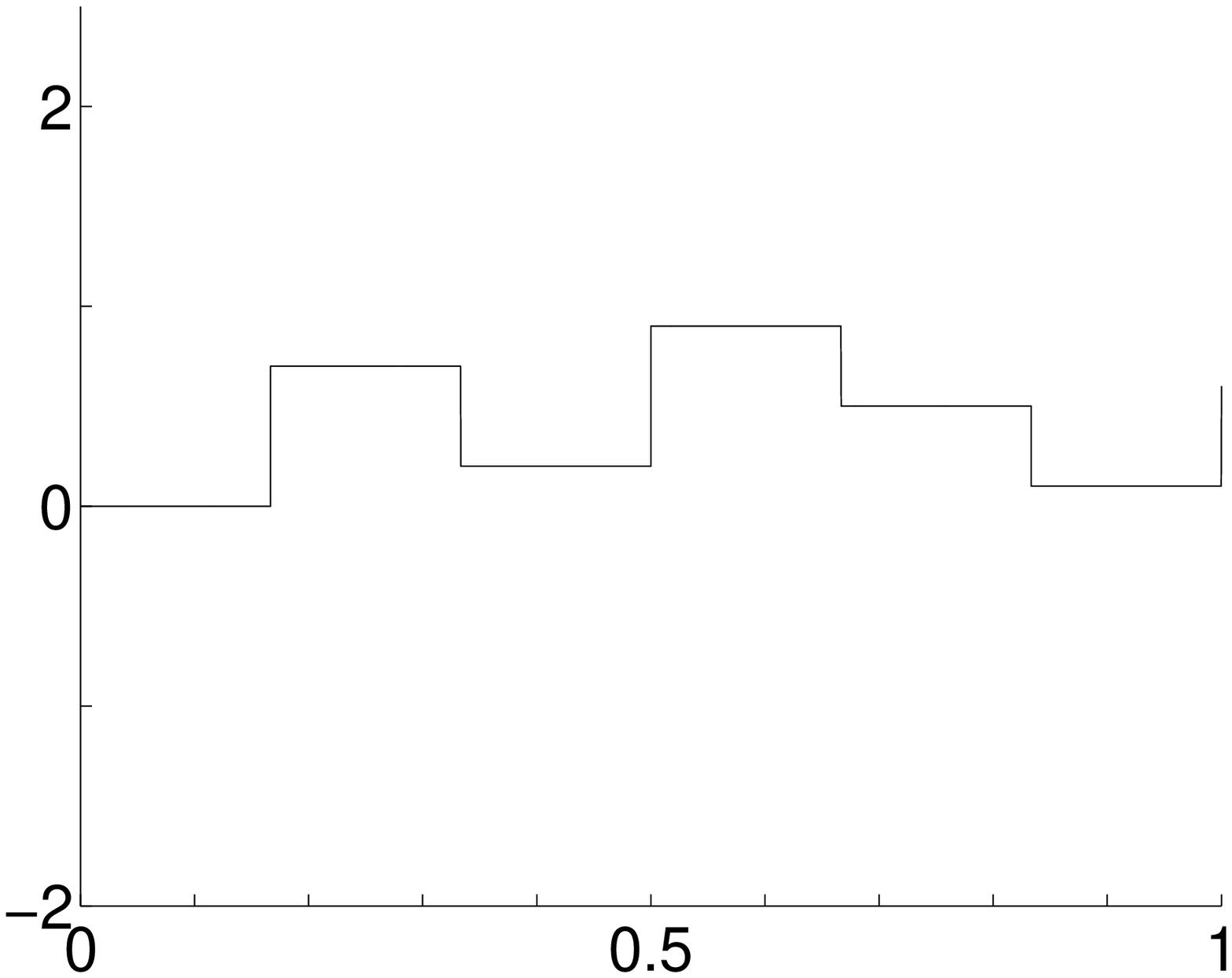,width=0.85\textwidth}}
   \caption{$\bayes (x) = \mathrm{His}_6 (x)$ \label{fig.His6.fonc}}
\end{minipage} \hfill
 \begin{minipage}[b]{.48\linewidth}   \centerline{\epsfig{file=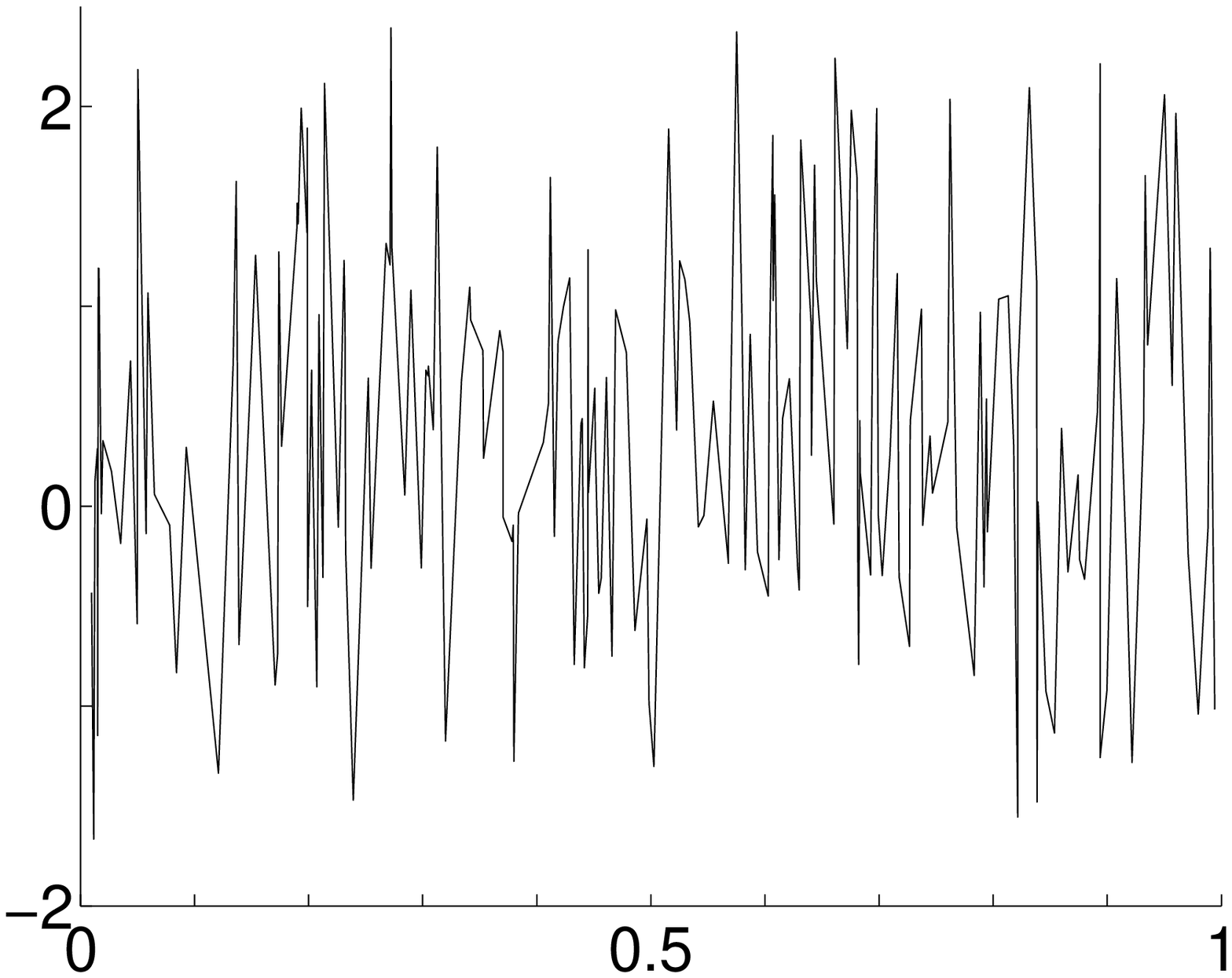,width=0.85\textwidth}}
   \caption{Data sample for His6 \label{fig.His6.data}}
\end{minipage}
\end{figure}

\begin{figure}
 \begin{minipage}[b]{.48\linewidth}   \centerline{\epsfig{file=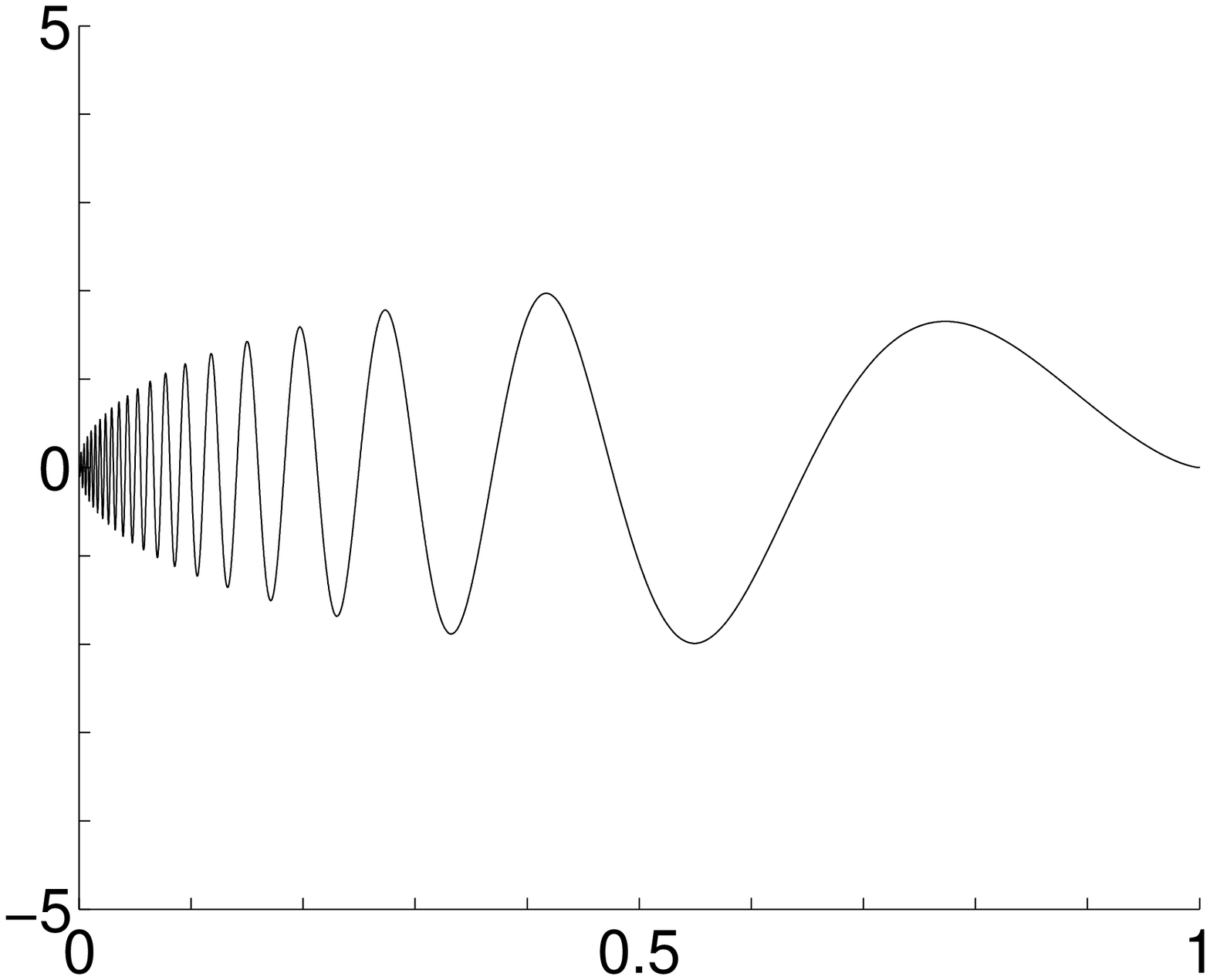,width=0.85\textwidth}}
   \caption{$\bayes(x) = \mathrm{Doppler}(x)$ (see \cite{Don_Joh:1995}) \label{fig.Dop.fonc}}
\end{minipage} \hfill
\begin{minipage}[b]{.48\linewidth}   
\centerline{\epsfig{file=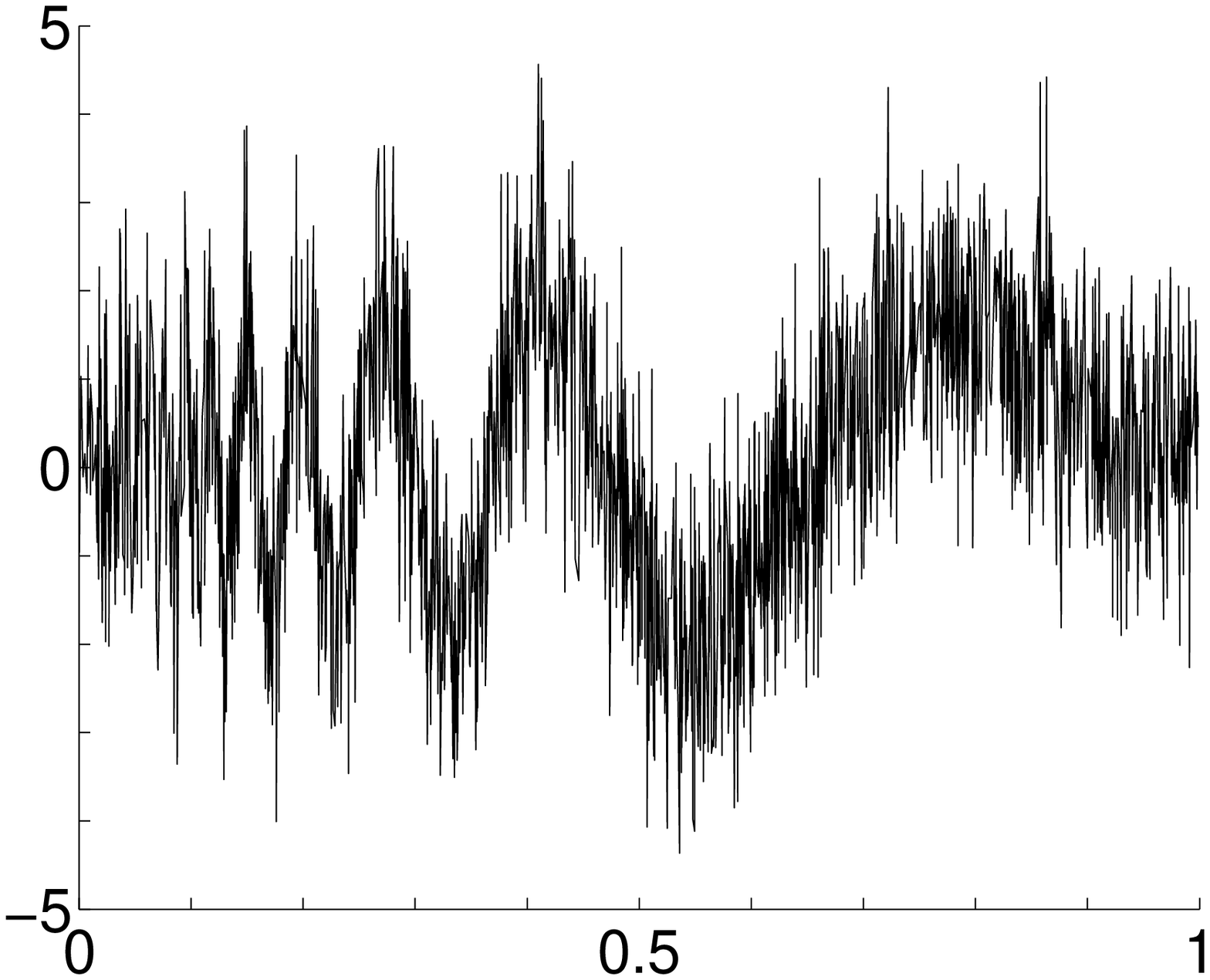,width=0.85\textwidth}}
   \caption{Data sample for DopReg and Dop2bin \label{fig.Dop.data}}
\end{minipage}
\end{figure}

Data are generated according to 
\[ Y_i = \bayes(X_i) + \sigma(X_i) \epsilon_i \]
with $X_i$ i.i.d. uniform on $\X=[0;1]$ and $\epsilon_i \sim \mathcal{N}(0,1)$ independent from $X_i$.
The experiments differ from 
\begin{itemize}
\item the regression function $\bayes$:
\begin{itemize}
\item S1000, S$\sqrt{0.1}$, S0.1 and Svar2 have the same smooth function as S1 and S2, see Fig.~\ref{fig.S.fonc}.
\item Sqrt has $\bayes(x) = \sqrt{x}$, which is smooth except around 0, see Fig.~\ref{fig.Sqrt.fonc}.
\item His6 has a regular histogram with 5 jumps (hence it belongs to the regular histogram model of dimension 6), see Fig.~\ref{fig.His6.fonc}.
\item DopReg and Dop2bin have the Doppler function, as defined by Donoho and Johnstone \cite{Don_Joh:1995}, see Fig.~\ref{fig.Dop.fonc}.
\end{itemize}
\item the noise level $\sigma$:
\begin{itemize}
\item $\sigma(x) = 1$ for S1000, Sqrt, His6, DopReg and Dop2bin.
\item $\sigma(x) = \sqrt{0.1}$ for S$\sqrt{0.1}$.
\item $\sigma(x) = 0.1$ for S0.1.
\item $\sigma(x) = \1_{x \geq 1/2}$ for Svar2.
\end{itemize}
\item the sample size $n$:
\begin{itemize}
\item $n=200$ for S$\sqrt{0.1}$, S0.1, Svar2, Sqrt and His6.
\item $n=1000$ for S1000.
\item $n=2048$ for DopReg and Dop2bin.
\end{itemize}
\item the family of models: with the notations introduced in Sect.~\main{4},
\begin{itemize}
\item for S1000, S$\sqrt{0.1}$, S0.1, Sqrt and His6, we use the ``regular'' collection, as for S1:
\[ \M_n = \set{ 1, \ldots, \left\lfloor \frac{n}{\ln(n)} \right\rfloor } \enspace . \]
\item for Svar2, we use the ``regular with two bin sizes'' collection, as for S2:
\[ \M_n = \set{1} \cup \set{ 1, \ldots, \left\lfloor \frac{n}{2\ln(n)} \right\rfloor }^2 \enspace . \]
\item for DopReg, we use the ``regular dyadic'' collection, as for HSd1:
\[ \M_n = \set{ 2^k \telque 0 \leq k \leq \ln_2(n)-1 } \enspace .\]
\item for Dop2bin, we use the ``regular dyadic with two bin sizes'' collection, as for HSd2:
\[ \M_n = \set{1} \cup \set{ 2^k \telque 0 \leq k \leq \ln_2(n)-2 }^2 \enspace .\]
\end{itemize}
\end{itemize}
Notice that contrary to HSd2, Dop2bin is an homoscedastic problem. The interest of considering two bin sizes for it is that the smoothness of the Doppler function is quite different for small $x$ and for $x \geq 1/2$.

\begin{figure}
\begin{minipage}[b]{.48\linewidth}   \centerline{\epsfig{file=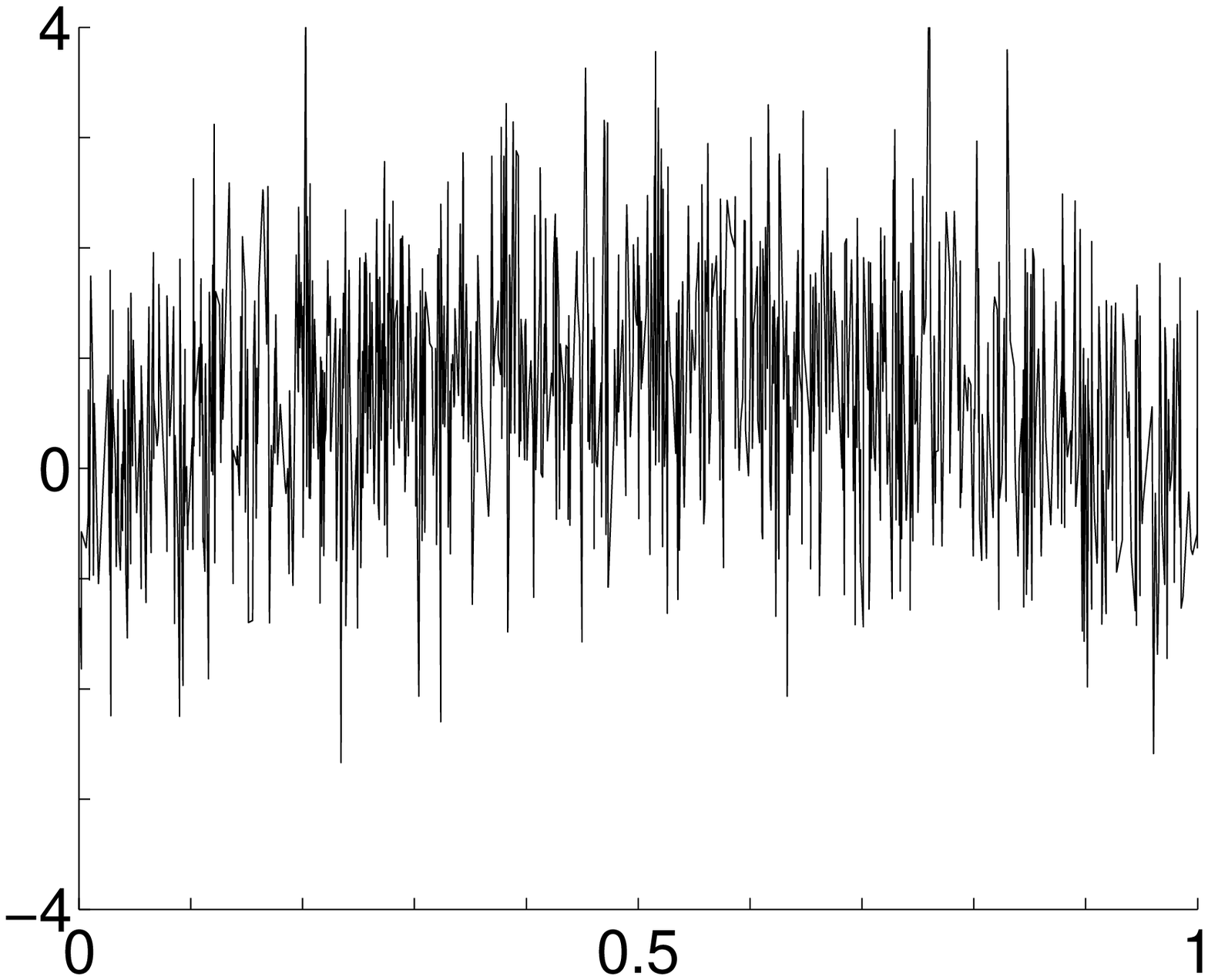,width=0.85\textwidth}}
   \caption{Data sample for S1000 \label{fig.S1000.data}}
\end{minipage} \hfill
 \begin{minipage}[b]{.48\linewidth}   \centerline{\epsfig{file=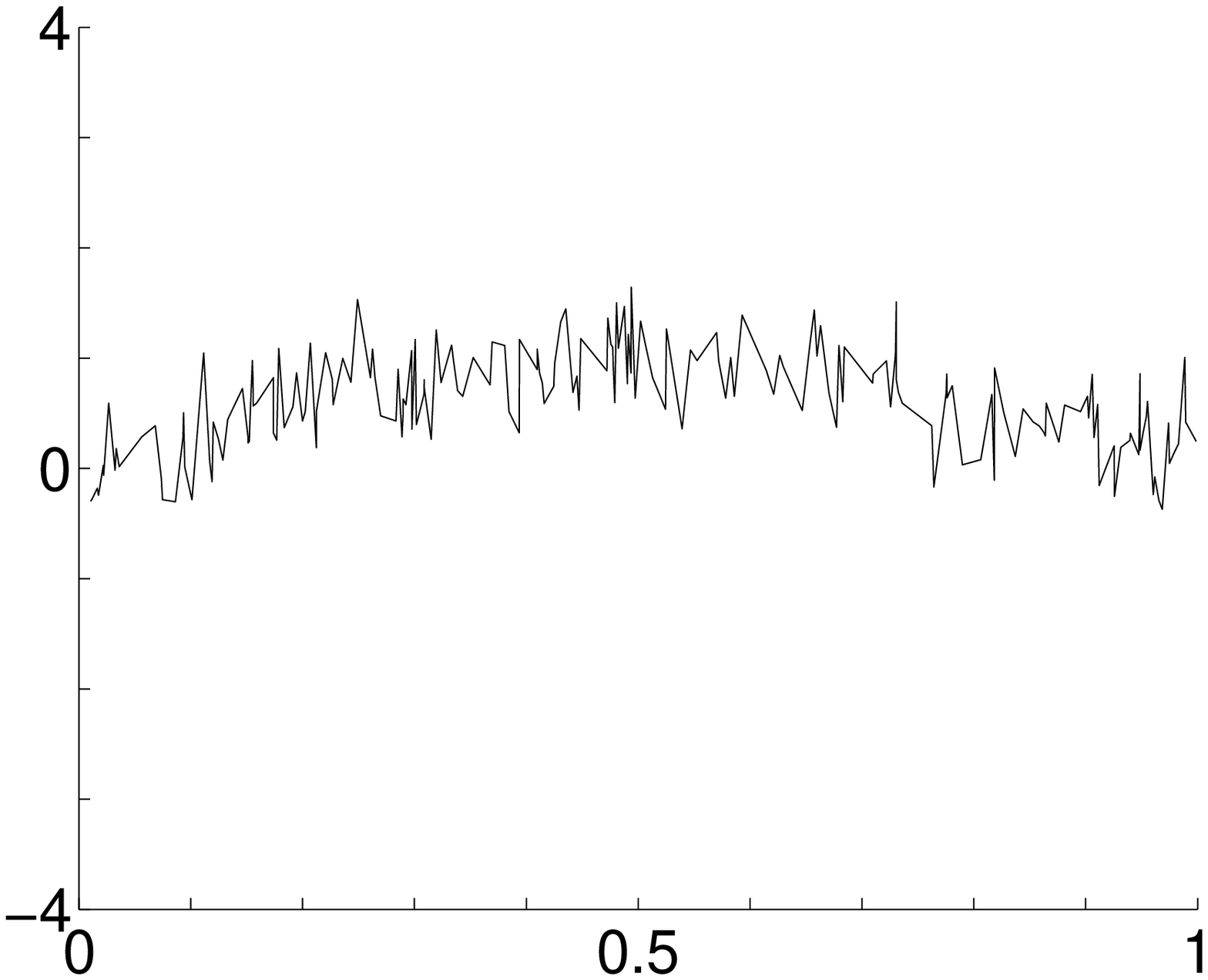,width=0.85\textwidth}}
   \caption{Data sample for S$\sqrt{0.1}$ \label{fig.Ssqrt0.1.data}}
\end{minipage}
\end{figure}

\begin{figure}
\begin{minipage}[b]{.48\linewidth}   \centerline{\epsfig{file=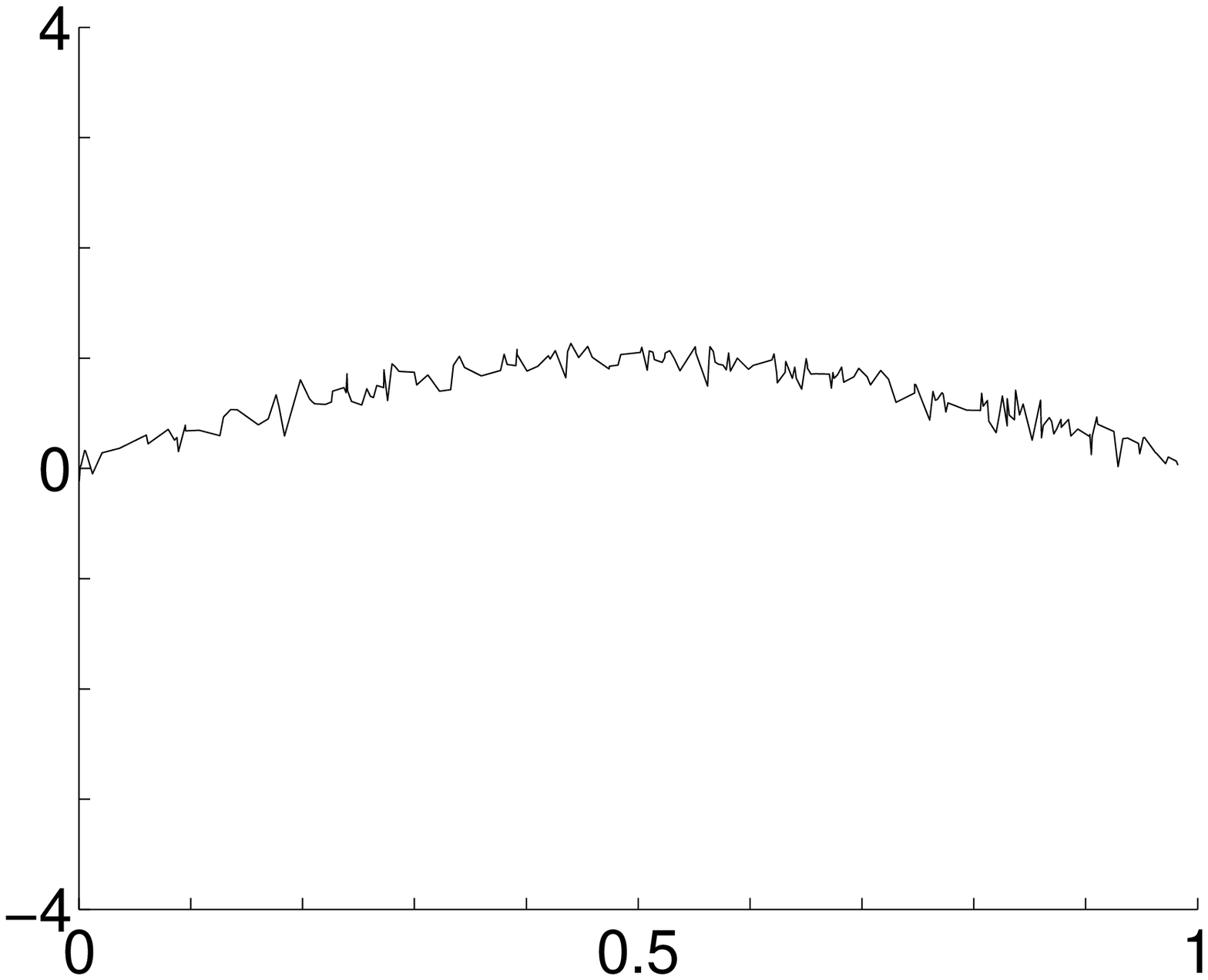,width=0.85\textwidth}}
   \caption{Data sample for S0.1 \label{fig.S0.1.data}}
\end{minipage} \hfill
 \begin{minipage}[b]{.48\linewidth}   \centerline{\epsfig{file=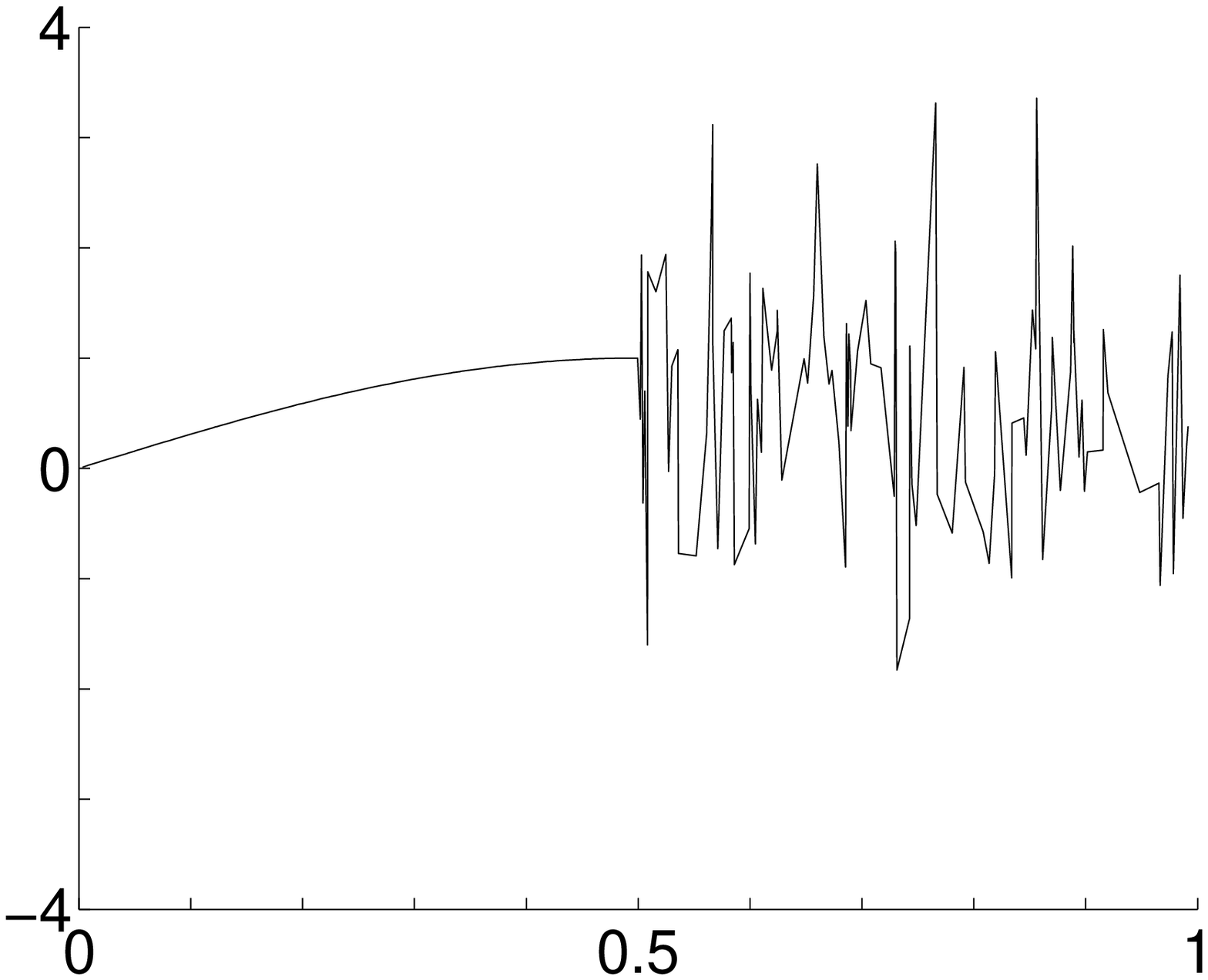,width=0.85\textwidth}}
   \caption{Data sample for Svar2 \label{fig.Svar2.data}}
\end{minipage}
\end{figure}

Instances of data sets for each experiment are given in Fig.~\ref{fig.S1000.data}--\ref{fig.Svar2.data}, \ref{fig.Sqrt.data}, \ref{fig.His6.data} and~\ref{fig.Dop.data}.


\medskip

Compared to S1, S2, HSd1 and HSd2, these eight experiments consider larger signal-to-noise ratio data (S1000, S$\sqrt{0.1}$, S0.1), another kind of heteroscedasticity (Svar2) and other regression functions, with different kinds of unsmoothness (Sqrt, His6, DopReg and Dop2bin).

We consider for each of these experiments the same algorithms as in Sect.~\main{4}, adding to them Mal$^{\star}$, which is Mallows' $C_p$ penalty with the true value of the variance: $\pen(m) = 2 \E\croch{\sigma^2(X)} D_m n^{-1}$. Although it can not be used on real data sets, it is an interesting point of comparison, which does not have possible weaknesses coming from the variance estimator $\widehat{\sigma}^2$.
Our estimates of $C_{\mathrm{or}}$ (and uncertainties for these estimates) for the procedures we consider are reported in Tab.~\ref{tab.un.Cor} to~\ref{tab.trois.Cor} (we report here again the results for S1, S2, HSd1 and HSd2 to make comparisons easier). On the last line of these Tables, we also report 
\[ \frac {\E\croch{ \inf_{\mM_n} \perte{\ERM_m} }}
{\inf_{\mM_n} \set{ \E\croch{\perte{\ERM_m}}}} 
= \frac { C_{\mathrm{or}}^{\prime}} 
{ C_{\mathrm{or}} } 
\qquad \mbox{where} \qquad 
C_{\mathrm{or}}^{\prime} \egaldef 
\frac{ \E \croch{ \perte{\ERM_{\mh}} } }
{\inf_{\mM_n} \set{ \E\croch{\perte{\ERM_m}}}} \]
is the leading constant which appear in most of the classical oracle inequalities. Notice that $C_{\mathrm{or}}^{\prime}$ is always smaller than $C_{\mathrm{or}}$.

It appears that the choice of $V$ is still difficult for VFCV: $V=2$ is optimal in S1000 and Sqrt and $V=20$ in the six other ones. On the contrary, $V = n$ is (almost) always better for penVF and penVF+, and overpenalization often improves the quality of the algorithm (but not always: see DopReg and S$0.1$). These eight experiments mainly show that the assumptions of Thm.~\main{2} are not necessary for penVF to be efficient.


\begin{table} 
\caption{Accuracy indexes $C_{\mathrm{or}}$ for experiments S1, S2, HSd1 and HSd2 ($N=1000$). Uncertainties reported are empirical standard deviations divided by $\sqrt{N}$. 
\label{tab.un.Cor}}
\begin{center}
\begin{tabular}
{p{0.16\textwidth}@{\hspace{0.025\textwidth}}p{0.16\textwidth}@{\hspace{0.025\textwidth}}p{0.16\textwidth}@{\hspace{0.025\textwidth}}p{0.17\textwidth}@{\hspace{0.025\textwidth}}p{0.17\textwidth}}
\hline\noalign{\smallskip}
Experiment & S1 & S2 & HSd1 & HSd2 \\ \noalign{\smallskip} \hline \noalign{\smallskip}  
 $\bayes$ & $\sin(\pi \cdot)$ & $\sin(\pi \cdot)$ & HeaviSine & HeaviSine \\  
 $\sigma(x)$ & 1 & $x$ & 1 & $x$ \\  
 $n$ (sample size) & 200 & 200 & 2048 & 2048 \\   
 $\M_n$ & regular & 2 bin sizes & dyadic, regular & dyadic, 2 bin sizes \\ \noalign{\smallskip} \hline \noalign{\smallskip}  
 Mal & $ 1.928 \pm 0.04 $& $ 3.687 \pm 0.07 $& $ 1.015 \pm 0.003 $& $ 1.373 \pm 0.010 $\\   
Mal+ & $ 1.800 \pm 0.03 $& $ 3.173 \pm 0.07 $& $ 1.002 \pm 0.003 $& $ 1.411 \pm 0.008 $\\   
Mal$^{\star}$ & $ 2.028 \pm 0.04 $& $ 2.657 \pm 0.06 $& $ 1.044 \pm 0.004 $& $ 1.513 \pm 0.005 $\\   
Mal$^{\star}$+ & $ 1.827 \pm 0.03 $& $ 2.437 \pm 0.05 $& $ 1.004 \pm 0.003 $& $ 1.548 \pm 0.003 $\\   
$\E\croch{\penid}$ & $ 1.919 \pm 0.03 $& $ 2.296 \pm 0.05 $& $ 1.028 \pm 0.004 $& $ 1.102 \pm 0.004 $\\   
$\E\croch{\penid}$+ & $ 1.792 \pm 0.03 $& $ 2.028 \pm 0.04 $& $ 1.003 \pm 0.003 $& $ 1.089 \pm 0.004 $\\ \noalign{\smallskip} \hline \noalign{\smallskip}  
 2-FCV & $ 2.078 \pm 0.04 $& $ 2.542 \pm 0.05 $& $ 1.002 \pm 0.003 $& $ 1.184 \pm 0.004 $\\   
5-FCV & $ 2.137 \pm 0.04 $& $ 2.582 \pm 0.06 $& $ 1.014 \pm 0.003 $& $ 1.115 \pm 0.005 $\\   
10-FCV & $ 2.097 \pm 0.04 $& $ 2.603 \pm 0.06 $& $ 1.021 \pm 0.003 $& $ 1.109 \pm 0.004 $\\   
20-FCV & $ 2.088 \pm 0.04 $& $ 2.578 \pm 0.06 $& $ 1.029 \pm 0.004 $& $ 1.105 \pm 0.004 $\\   
LOO & $ 2.077 \pm 0.04 $& $ 2.593 \pm 0.06 $& $ 1.034 \pm 0.004 $& $ 1.105 \pm 0.004 $\\ \noalign{\smallskip} \hline \noalign{\smallskip}  
 pen2-F & $ 2.578 \pm 0.06 $& $ 3.061 \pm 0.07 $& $ 1.038 \pm 0.004 $& $ 1.103 \pm 0.004 $\\   
pen5-F & $ 2.219 \pm 0.05 $& $ 2.750 \pm 0.06 $& $ 1.037 \pm 0.004 $& $ 1.104 \pm 0.004 $\\   
pen10-F & $ 2.121 \pm 0.04 $& $ 2.653 \pm 0.06 $& $ 1.034 \pm 0.004 $& $ 1.104 \pm 0.004 $\\   
pen20-F & $ 2.085 \pm 0.04 $& $ 2.639 \pm 0.06 $& $ 1.034 \pm 0.004 $& $ 1.105 \pm 0.004 $\\   
penLoo & $ 2.080 \pm 0.04 $& $ 2.593 \pm 0.06 $& $ 1.034 \pm 0.004 $& $ 1.105 \pm 0.004 $\\ \noalign{\smallskip} \hline \noalign{\smallskip}  
 pen2-F+ & $ 2.175 \pm 0.05 $& $ 2.748 \pm 0.06 $& $ 1.011 \pm 0.003 $& $ 1.106 \pm 0.004 $\\   
pen5-F+ & $ 1.913 \pm 0.03 $& $ 2.378 \pm 0.05 $& $ 1.006 \pm 0.003 $& $ 1.102 \pm 0.004 $\\   
pen10-F+ & $ 1.872 \pm 0.03 $& $ 2.285 \pm 0.05 $& $ 1.005 \pm 0.003 $& $ 1.098 \pm 0.004 $\\   
pen20-F+ & $ 1.898 \pm 0.03 $& $ 2.254 \pm 0.05 $& $ 1.004 \pm 0.003 $& $ 1.098 \pm 0.004 $\\   
penLoo+ & $ 1.844 \pm 0.03 $& $ 2.215 \pm 0.05 $& $ 1.004 \pm 0.003 $& $ 1.096 \pm 0.004 $\\ \noalign{\smallskip} \hline \noalign{\smallskip}  
 $C^{\prime}_{\mathrm{or}} / C_{\mathrm{or}}$  & 0.768 & 0.753 & 0.999 & 0.854 \\ \noalign{\smallskip} \hline  
 
\end{tabular}
\end{center}
\end{table}

\begin{table} 
\caption{Accuracy indexes $C_{\mathrm{or}}$ for experiments S1000, S$\sqrt{0.1}$, S$0.1$ and Svar2 ($N=250$). Uncertainties reported are empirical standard deviations divided by $\sqrt{N}$. 
\label{tab.deux.Cor}}
\begin{center}
\begin{tabular}
{p{0.16\textwidth}@{\hspace{0.025\textwidth}}p{0.16\textwidth}@{\hspace{0.025\textwidth}}p{0.16\textwidth}@{\hspace{0.025\textwidth}}p{0.17\textwidth}@{\hspace{0.025\textwidth}}p{0.17\textwidth}}
\hline\noalign{\smallskip}
Experiment & S1000 & S$\sqrt{0.1}$ & S$0.1$ & Svar2 \\ \noalign{\smallskip} \hline \noalign{\smallskip}  
 $\bayes$ & $\sin(\pi \cdot)$ & $\sin(\pi \cdot)$ & $\sin(\pi \cdot)$ & $\sin(\pi \cdot)$ \\  
 $\sigma(x)$ & 1 & $\sqrt{0.1}$ & $0.1$ & $\1_{x \geq 1/2}$ \\  
 $n$ (sample size) & 1000 & 200 & 200 & 200 \\   
 $\M_n$ & regular & regular & regular & 2 bin sizes \\ \noalign{\smallskip} \hline \noalign{\smallskip}  
 Mal & $ 1.667 \pm 0.04 $& $ 1.611 \pm 0.03 $& $ 1.400 \pm 0.02 $& $ 5.643 \pm 0.22 $\\   
Mal+ & $ 1.619 \pm 0.03 $& $ 1.593 \pm 0.03 $& $ 1.426 \pm 0.02 $& $ 4.647 \pm 0.22 $\\   
Mal$^{\star}$ & $ 1.745 \pm 0.05 $& $ 1.925 \pm 0.03 $& $ 3.204 \pm 0.05 $& $ 4.481 \pm 0.21 $\\   
Mal$^{\star}$+ & $ 1.617 \pm 0.03 $& $ 2.073 \pm 0.04 $& $ 3.641 \pm 0.07 $& $ 3.544 \pm 0.17 $\\   
$\E\croch{\penid}$ & $ 1.745 \pm 0.05 $& $ 1.571 \pm 0.03 $& $ 1.373 \pm 0.02 $& $ 2.409 \pm 0.13 $\\   
$\E\croch{\penid}$+ & $ 1.617 \pm 0.03 $& $ 1.554 \pm 0.03 $& $ 1.392 \pm 0.02 $& $ 2.005 \pm 0.10 $\\ \noalign{\smallskip} \hline \noalign{\smallskip}  
 2-FCV & $ 1.668 \pm 0.04 $& $ 1.663 \pm 0.04 $& $ 1.394 \pm 0.02 $& $ 2.960 \pm 0.15 $\\   
5-FCV & $ 1.756 \pm 0.07 $& $ 1.693 \pm 0.04 $& $ 1.393 \pm 0.02 $& $ 2.950 \pm 0.16 $\\   
10-FCV & $ 1.746 \pm 0.04 $& $ 1.684 \pm 0.04 $& $ 1.385 \pm 0.02 $& $ 2.681 \pm 0.14 $\\   
20-FCV & $ 1.774 \pm 0.05 $& $ 1.645 \pm 0.03 $& $ 1.382 \pm 0.02 $& $ 2.742 \pm 0.16 $\\   
LOO & $ 1.768 \pm 0.05 $& $ 1.639 \pm 0.04 $& $ 1.379 \pm 0.02 $& $ 2.641 \pm 0.15 $\\ \noalign{\smallskip} \hline \noalign{\smallskip}  
 pen2-F & $ 2.066 \pm 0.08 $& $ 1.809 \pm 0.05 $& $ 1.390 \pm 0.02 $& $ 3.209 \pm 0.18 $\\   
pen5-F & $ 1.816 \pm 0.05 $& $ 1.638 \pm 0.04 $& $ 1.400 \pm 0.02 $& $ 2.749 \pm 0.15 $\\   
pen10-F & $ 1.783 \pm 0.05 $& $ 1.706 \pm 0.04 $& $ 1.374 \pm 0.02 $& $ 2.598 \pm 0.15 $\\   
pen20-F & $ 1.801 \pm 0.05 $& $ 1.657 \pm 0.03 $& $ 1.385 \pm 0.02 $& $ 2.684 \pm 0.15 $\\   
penLoo & $ 1.776 \pm 0.05 $& $ 1.641 \pm 0.04 $& $ 1.379 \pm 0.02 $& $ 2.656 \pm 0.15 $\\ \noalign{\smallskip} \hline \noalign{\smallskip}  
 pen2-F+ & $ 1.809 \pm 0.05 $& $ 1.714 \pm 0.04 $& $ 1.416 \pm 0.02 $& $ 2.808 \pm 0.16 $\\   
pen5-F+ & $ 1.683 \pm 0.04 $& $ 1.616 \pm 0.03 $& $ 1.399 \pm 0.02 $& $ 2.460 \pm 0.14 $\\   
pen10-F+ & $ 1.627 \pm 0.04 $& $ 1.613 \pm 0.03 $& $ 1.385 \pm 0.02 $& $ 2.398 \pm 0.14 $\\   
pen20-F+ & $ 1.644 \pm 0.04 $& $ 1.583 \pm 0.03 $& $ 1.390 \pm 0.02 $& $ 2.316 \pm 0.13 $\\   
penLoo+ & $ 1.626 \pm 0.03 $& $ 1.587 \pm 0.03 $& $ 1.401 \pm 0.02 $& $ 2.349 \pm 0.13 $\\ \noalign{\smallskip} \hline \noalign{\smallskip}  
 $C^{\prime}_{\mathrm{or}} / C_{\mathrm{or}}$  & 0.8 & 0.801 & 0.816 & 0.779 \\ \noalign{\smallskip} \hline  
 
\end{tabular}
\end{center}
\end{table}

\begin{table} 
\caption{Accuracy indexes $C_{\mathrm{or}}$ for experiments Sqrt, His6, DopReg and Dop2bin ($N=250$). Uncertainties reported are empirical standard deviations divided by $\sqrt{N}$. 
\label{tab.trois.Cor}}
\begin{center}
\begin{tabular}
{p{0.16\textwidth}@{\hspace{0.025\textwidth}}p{0.16\textwidth}@{\hspace{0.025\textwidth}}p{0.16\textwidth}@{\hspace{0.025\textwidth}}p{0.17\textwidth}@{\hspace{0.025\textwidth}}p{0.17\textwidth}}
\hline\noalign{\smallskip}
Experiment & Sqrt & His6 & DopReg & Dop2bin \\ \noalign{\smallskip} \hline \noalign{\smallskip}  
 $\bayes$ & $\sqrt{\cdot}$ & His$_6$ & Doppler & Doppler \\  
 $\sigma(x)$ & 1 & 1 & 1 & 1 \\  
 $n$ (sample size) & 200 & 200 & 2048 & 2048 \\   
 $\M_n$ & regular & regular & dyadic, regular & dyadic, 2 bin sizes \\ \noalign{\smallskip} \hline \noalign{\smallskip}  
 Mal & $ 2.295 \pm 0.11 $& $ 1.969 \pm 0.11 $& $ 1.039 \pm 0.01 $& $ 1.052 \pm 0.01 $\\   
Mal+ & $ 1.989 \pm 0.08 $& $ 1.799 \pm 0.09 $& $ 1.090 \pm 0.00 $& $ 1.047 \pm 0.01 $\\   
Mal$^{\star}$ & $ 2.483 \pm 0.12 $& $ 2.021 \pm 0.11 $& $ 1.013 \pm 0.01 $& $ 1.061 \pm 0.01 $\\   
Mal$^{\star}$+ & $ 2.075 \pm 0.09 $& $ 1.836 \pm 0.10 $& $ 1.070 \pm 0.00 $& $ 1.041 \pm 0.01 $\\   
$\E\croch{\penid}$ & $ 2.365 \pm 0.11 $& $ 1.805 \pm 0.10 $& $ 1.025 \pm 0.01 $& $ 1.056 \pm 0.01 $\\   
$\E\croch{\penid}$+ & $ 2.012 \pm 0.09 $& $ 1.632 \pm 0.08 $& $ 1.083 \pm 0.00 $& $ 1.040 \pm 0.01 $\\ \noalign{\smallskip} \hline \noalign{\smallskip}  
 2-FCV & $ 2.489 \pm 0.12 $& $ 2.788 \pm 0.13 $& $ 1.097 \pm 0.00 $& $ 1.165 \pm 0.01 $\\   
5-FCV & $ 2.777 \pm 0.16 $& $ 2.316 \pm 0.12 $& $ 1.064 \pm 0.01 $& $ 1.049 \pm 0.01 $\\   
10-FCV & $ 2.571 \pm 0.13 $& $ 2.074 \pm 0.11 $& $ 1.043 \pm 0.01 $& $ 1.051 \pm 0.01 $\\   
20-FCV & $ 2.561 \pm 0.12 $& $ 2.071 \pm 0.11 $& $ 1.034 \pm 0.01 $& $ 1.053 \pm 0.01 $\\   
LOO & $ 2.695 \pm 0.14 $& $ 2.059 \pm 0.11 $& $ 1.026 \pm 0.01 $& $ 1.058 \pm 0.01 $\\ \noalign{\smallskip} \hline \noalign{\smallskip}  
 pen2-F & $ 4.088 \pm 0.23 $& $ 3.210 \pm 0.14 $& $ 1.048 \pm 0.01 $& $ 1.062 \pm 0.01 $\\   
pen5-F & $ 3.024 \pm 0.18 $& $ 2.485 \pm 0.13 $& $ 1.033 \pm 0.01 $& $ 1.055 \pm 0.01 $\\   
pen10-F & $ 3.009 \pm 0.18 $& $ 2.192 \pm 0.12 $& $ 1.029 \pm 0.01 $& $ 1.056 \pm 0.01 $\\   
pen20-F & $ 2.723 \pm 0.14 $& $ 2.150 \pm 0.12 $& $ 1.031 \pm 0.01 $& $ 1.056 \pm 0.01 $\\   
penLoo & $ 2.695 \pm 0.14 $& $ 2.063 \pm 0.12 $& $ 1.026 \pm 0.01 $& $ 1.058 \pm 0.01 $\\ \noalign{\smallskip} \hline \noalign{\smallskip}  
 pen2-F+ & $ 3.015 \pm 0.17 $& $ 2.728 \pm 0.12 $& $ 1.084 \pm 0.00 $& $ 1.084 \pm 0.01 $\\   
pen5-F+ & $ 2.409 \pm 0.13 $& $ 2.080 \pm 0.09 $& $ 1.080 \pm 0.00 $& $ 1.063 \pm 0.01 $\\   
pen10-F+ & $ 2.305 \pm 0.11 $& $ 1.869 \pm 0.09 $& $ 1.082 \pm 0.00 $& $ 1.050 \pm 0.01 $\\   
pen20-F+ & $ 2.180 \pm 0.10 $& $ 1.832 \pm 0.09 $& $ 1.079 \pm 0.00 $& $ 1.052 \pm 0.01 $\\   
penLoo+ & $ 2.152 \pm 0.10 $& $ 1.858 \pm 0.10 $& $ 1.082 \pm 0.00 $& $ 1.048 \pm 0.01 $\\ \noalign{\smallskip} \hline \noalign{\smallskip}  
 $C^{\prime}_{\mathrm{or}} / C_{\mathrm{or}}$  & 0.795 & 0.996 & 0.998 & 0.977 \\ \noalign{\smallskip} \hline  
 
\end{tabular}
\end{center}
\end{table}

\medskip

For the sake of completeness, we also reported the results for the twelve experiments in terms of the other benchmark 
\[ C_{\mathrm{path-or}} \egaldef \E\croch{ \frac{\perte{\ERM_{\mh}} } {\inf_{\mM_n} \perte{\ERM_m} }} \]
in Tab.~\ref{tab.un.Cpath-or} to Tab.~\ref{tab.trois.Cpath-or}. They are indeed quite similar to the previous ones.


\begin{table} 
\caption{Accuracy indexes $C_{\mathrm{path-or}}$ for experiments S1, S2, HSd1 and HSd2 ($N=1000$). Uncertainties reported are empirical standard deviations divided by $\sqrt{N}$. 
\label{tab.un.Cpath-or}}
\begin{center}
\begin{tabular}
{p{0.16\textwidth}@{\hspace{0.025\textwidth}}p{0.16\textwidth}@{\hspace{0.025\textwidth}}p{0.16\textwidth}@{\hspace{0.025\textwidth}}p{0.17\textwidth}@{\hspace{0.025\textwidth}}p{0.17\textwidth}}
\hline\noalign{\smallskip}
Experiment & S1 & S2 & HSd1 & HSd2 \\ \noalign{\smallskip} \hline \noalign{\smallskip}  
 $\bayes$ & $\sin(\pi \cdot)$ & $\sin(\pi \cdot)$ & HeaviSine & HeaviSine \\  
 $\sigma(x)$ & 1 & $x$ & 1 & $x$ \\  
 $n$ (sample size) & 200 & 200 & 2048 & 2048 \\   
 $\M_n$ & regular & 2 bin sizes & dyadic, regular & dyadic, 2 bin sizes \\ \noalign{\smallskip} \hline \noalign{\smallskip}  
 Mal & $ 2.064 \pm 0.04 $& $ 4.129 \pm 0.10 $& $ 1.015 \pm 0.002 $& $ 1.316 \pm 0.010 $\\   
Mal+ & $ 1.921 \pm 0.03 $& $ 3.500 \pm 0.09 $& $ 1.002 \pm 0.001 $& $ 1.354 \pm 0.008 $\\   
Mal$^{\star}$ & $ 2.168 \pm 0.04 $& $ 2.907 \pm 0.07 $& $ 1.045 \pm 0.003 $& $ 1.453 \pm 0.006 $\\   
Mal$^{\star}$+ & $ 1.941 \pm 0.03 $& $ 2.645 \pm 0.06 $& $ 1.004 \pm 0.001 $& $ 1.487 \pm 0.005 $\\   
$\E\croch{\penid}$ & $ 2.053 \pm 0.04 $& $ 2.458 \pm 0.06 $& $ 1.029 \pm 0.003 $& $ 1.050 \pm 0.002 $\\   
$\E\croch{\penid}$+ & $ 1.903 \pm 0.03 $& $ 2.142 \pm 0.04 $& $ 1.003 \pm 0.001 $& $ 1.038 \pm 0.002 $\\ \noalign{\smallskip} \hline \noalign{\smallskip}  
 2-FCV & $ 2.230 \pm 0.05 $& $ 2.755 \pm 0.06 $& $ 1.002 \pm 0.001 $& $ 1.134 \pm 0.004 $\\   
5-FCV & $ 2.290 \pm 0.05 $& $ 2.827 \pm 0.08 $& $ 1.014 \pm 0.002 $& $ 1.064 \pm 0.003 $\\   
10-FCV & $ 2.237 \pm 0.05 $& $ 2.832 \pm 0.08 $& $ 1.021 \pm 0.002 $& $ 1.057 \pm 0.002 $\\   
20-FCV & $ 2.225 \pm 0.05 $& $ 2.794 \pm 0.07 $& $ 1.029 \pm 0.003 $& $ 1.054 \pm 0.002 $\\   
LOO & $ 2.212 \pm 0.05 $& $ 2.832 \pm 0.08 $& $ 1.034 \pm 0.003 $& $ 1.053 \pm 0.002 $\\ \noalign{\smallskip} \hline \noalign{\smallskip}  
 pen2-F & $ 2.770 \pm 0.07 $& $ 3.340 \pm 0.08 $& $ 1.039 \pm 0.003 $& $ 1.052 \pm 0.003 $\\   
pen5-F & $ 2.383 \pm 0.06 $& $ 2.982 \pm 0.08 $& $ 1.038 \pm 0.003 $& $ 1.053 \pm 0.002 $\\   
pen10-F & $ 2.256 \pm 0.05 $& $ 2.867 \pm 0.07 $& $ 1.035 \pm 0.003 $& $ 1.053 \pm 0.002 $\\   
pen20-F & $ 2.219 \pm 0.05 $& $ 2.869 \pm 0.08 $& $ 1.035 \pm 0.003 $& $ 1.053 \pm 0.002 $\\   
penLoo & $ 2.215 \pm 0.05 $& $ 2.832 \pm 0.08 $& $ 1.034 \pm 0.003 $& $ 1.053 \pm 0.002 $\\ \noalign{\smallskip} \hline \noalign{\smallskip}  
 pen2-F+ & $ 2.328 \pm 0.05 $& $ 2.979 \pm 0.07 $& $ 1.011 \pm 0.002 $& $ 1.056 \pm 0.003 $\\   
pen5-F+ & $ 2.050 \pm 0.04 $& $ 2.540 \pm 0.06 $& $ 1.006 \pm 0.001 $& $ 1.052 \pm 0.002 $\\   
pen10-F+ & $ 1.997 \pm 0.03 $& $ 2.436 \pm 0.05 $& $ 1.005 \pm 0.001 $& $ 1.048 \pm 0.002 $\\   
pen20-F+ & $ 2.018 \pm 0.04 $& $ 2.416 \pm 0.06 $& $ 1.004 \pm 0.001 $& $ 1.047 \pm 0.002 $\\   
penLoo+ & $ 1.959 \pm 0.03 $& $ 2.397 \pm 0.06 $& $ 1.004 \pm 0.001 $& $ 1.045 \pm 0.002 $\\ \noalign{\smallskip} \hline \noalign{\smallskip}  
 
\end{tabular}
\end{center}
\end{table}

\begin{table} 
\caption{Accuracy indexes $C_{\mathrm{path-or}}$ for experiments S1000, S$\sqrt{0.1}$, S$0.1$ and Svar2 ($N=250$). Uncertainties reported are empirical standard deviations divided by $\sqrt{N}$. 
\label{tab.deux.Cpath-or}}
\begin{center}
\begin{tabular}
{p{0.16\textwidth}@{\hspace{0.025\textwidth}}p{0.16\textwidth}@{\hspace{0.025\textwidth}}p{0.16\textwidth}@{\hspace{0.025\textwidth}}p{0.17\textwidth}@{\hspace{0.025\textwidth}}p{0.17\textwidth}}
\hline\noalign{\smallskip}
Experiment & S1000 & S$\sqrt{0.1}$ & S$0.1$ & Svar2 \\ \noalign{\smallskip} \hline \noalign{\smallskip}  
 $\bayes$ & $\sin(\pi \cdot)$ & $\sin(\pi \cdot)$ & $\sin(\pi \cdot)$ & $\sin(\pi \cdot)$ \\  
 $\sigma(x)$ & 1 & $\sqrt{0.1}$ & $0.1$ & $\1_{x \geq 1/2}$ \\  
 $n$ (sample size) & 1000 & 200 & 200 & 200 \\   
 $\M_n$ & regular & regular & regular & 2 bin sizes \\ \noalign{\smallskip} \hline \noalign{\smallskip}  
 Mal & $ 1.704 \pm 0.04 $& $ 1.654 \pm 0.03 $& $ 1.407 \pm 0.02 $& $ 7.212 \pm 0.40 $\\   
Mal+ & $ 1.670 \pm 0.03 $& $ 1.636 \pm 0.03 $& $ 1.436 \pm 0.02 $& $ 5.740 \pm 0.34 $\\   
Mal$^{\star}$ & $ 1.793 \pm 0.04 $& $ 2.018 \pm 0.04 $& $ 3.273 \pm 0.06 $& $ 5.597 \pm 0.33 $\\   
Mal$^{\star}$+ & $ 1.664 \pm 0.03 $& $ 2.175 \pm 0.05 $& $ 3.719 \pm 0.08 $& $ 4.284 \pm 0.25 $\\   
$\E\croch{\penid}$ & $ 1.793 \pm 0.04 $& $ 1.611 \pm 0.03 $& $ 1.378 \pm 0.01 $& $ 2.785 \pm 0.19 $\\   
$\E\croch{\penid}$+ & $ 1.194 \pm 0.02 $& $ 1.177 \pm 0.02 $& $ 1.128 \pm 0.01 $& $ 1.337 \pm 0.07 $\\ \noalign{\smallskip} \hline \noalign{\smallskip}  
 2-FCV & $ 1.721 \pm 0.04 $& $ 1.723 \pm 0.04 $& $ 1.400 \pm 0.02 $& $ 3.507 \pm 0.19 $\\   
5-FCV & $ 1.801 \pm 0.06 $& $ 1.740 \pm 0.04 $& $ 1.399 \pm 0.02 $& $ 3.486 \pm 0.24 $\\   
10-FCV & $ 1.802 \pm 0.05 $& $ 1.735 \pm 0.04 $& $ 1.388 \pm 0.02 $& $ 3.149 \pm 0.20 $\\   
20-FCV & $ 1.832 \pm 0.05 $& $ 1.687 \pm 0.03 $& $ 1.388 \pm 0.02 $& $ 3.257 \pm 0.23 $\\   
LOO & $ 1.815 \pm 0.05 $& $ 1.685 \pm 0.04 $& $ 1.385 \pm 0.01 $& $ 3.127 \pm 0.24 $\\ \noalign{\smallskip} \hline \noalign{\smallskip}  
 pen2-F & $ 2.108 \pm 0.07 $& $ 1.864 \pm 0.05 $& $ 1.394 \pm 0.02 $& $ 3.839 \pm 0.27 $\\   
pen5-F & $ 1.852 \pm 0.05 $& $ 1.675 \pm 0.04 $& $ 1.404 \pm 0.02 $& $ 3.237 \pm 0.23 $\\   
pen10-F & $ 1.812 \pm 0.05 $& $ 1.767 \pm 0.04 $& $ 1.381 \pm 0.01 $& $ 3.093 \pm 0.23 $\\   
pen20-F & $ 1.839 \pm 0.05 $& $ 1.706 \pm 0.03 $& $ 1.391 \pm 0.01 $& $ 3.123 \pm 0.23 $\\   
penLoo & $ 1.825 \pm 0.05 $& $ 1.687 \pm 0.04 $& $ 1.385 \pm 0.01 $& $ 3.152 \pm 0.24 $\\ \noalign{\smallskip} \hline \noalign{\smallskip}  
 pen2-F+ & $ 1.852 \pm 0.05 $& $ 1.765 \pm 0.05 $& $ 1.420 \pm 0.02 $& $ 3.336 \pm 0.23 $\\   
pen5-F+ & $ 1.732 \pm 0.04 $& $ 1.664 \pm 0.03 $& $ 1.408 \pm 0.02 $& $ 2.890 \pm 0.22 $\\   
pen10-F+ & $ 1.663 \pm 0.04 $& $ 1.657 \pm 0.03 $& $ 1.394 \pm 0.02 $& $ 2.810 \pm 0.21 $\\   
pen20-F+ & $ 1.680 \pm 0.04 $& $ 1.623 \pm 0.03 $& $ 1.397 \pm 0.01 $& $ 2.657 \pm 0.19 $\\   
penLoo+ & $ 1.673 \pm 0.03 $& $ 1.624 \pm 0.03 $& $ 1.409 \pm 0.02 $& $ 2.659 \pm 0.18 $\\ \noalign{\smallskip} \hline \noalign{\smallskip}  
 
\end{tabular}
\end{center}
\end{table}

\begin{table} 
\caption{Accuracy indexes $C_{\mathrm{path-or}}$ for experiments Sqrt, His6, DopReg and Dop2bin ($N=250$). Uncertainties reported are empirical standard deviations divided by $\sqrt{N}$. 
\label{tab.trois.Cpath-or}}
\begin{center}
\begin{tabular}
{p{0.16\textwidth}@{\hspace{0.025\textwidth}}p{0.16\textwidth}@{\hspace{0.025\textwidth}}p{0.16\textwidth}@{\hspace{0.025\textwidth}}p{0.17\textwidth}@{\hspace{0.025\textwidth}}p{0.17\textwidth}}
\hline\noalign{\smallskip}
Experiment & Sqrt & His6 & DopReg & Dop2bin \\ \noalign{\smallskip} \hline \noalign{\smallskip}  
 $\bayes$ & $\sqrt{\cdot}$ & His$_6$ & Doppler & Doppler \\  
 $\sigma(x)$ & 1 & 1 & 1 & 1 \\  
 $n$ (sample size) & 200 & 200 & 2048 & 2048 \\   
 $\M_n$ & regular & regular & dyadic, regular & dyadic, 2 bin sizes \\ \noalign{\smallskip} \hline \noalign{\smallskip}  
 Mal & $ 2.557 \pm 0.12 $& $ 2.356 \pm 0.18 $& $ 1.040 \pm 0.00 $& $ 1.049 \pm 0.00 $\\   
Mal+ & $ 2.232 \pm 0.10 $& $ 2.041 \pm 0.12 $& $ 1.094 \pm 0.00 $& $ 1.045 \pm 0.01 $\\   
Mal$^{\star}$ & $ 2.838 \pm 0.15 $& $ 2.533 \pm 0.21 $& $ 1.013 \pm 0.00 $& $ 1.057 \pm 0.00 $\\   
Mal$^{\star}$+ & $ 2.349 \pm 0.11 $& $ 2.168 \pm 0.16 $& $ 1.073 \pm 0.00 $& $ 1.038 \pm 0.00 $\\   
$\E\croch{\penid}$ & $ 2.678 \pm 0.14 $& $ 2.182 \pm 0.17 $& $ 1.026 \pm 0.00 $& $ 1.053 \pm 0.00 $\\   
$\E\croch{\penid}$+ & $ 1.348 \pm 0.07 $& $ 1.230 \pm 0.06 $& $ 1.050 \pm 0.00 $& $ 1.038 \pm 0.00 $\\ \noalign{\smallskip} \hline \noalign{\smallskip}  
 2-FCV & $ 2.974 \pm 0.17 $& $ 3.713 \pm 0.25 $& $ 1.100 \pm 0.00 $& $ 1.164 \pm 0.01 $\\   
5-FCV & $ 3.209 \pm 0.21 $& $ 2.977 \pm 0.24 $& $ 1.066 \pm 0.00 $& $ 1.046 \pm 0.00 $\\   
10-FCV & $ 2.912 \pm 0.16 $& $ 2.639 \pm 0.21 $& $ 1.045 \pm 0.00 $& $ 1.047 \pm 0.00 $\\   
20-FCV & $ 2.889 \pm 0.15 $& $ 2.584 \pm 0.20 $& $ 1.035 \pm 0.00 $& $ 1.050 \pm 0.00 $\\   
LOO & $ 3.061 \pm 0.17 $& $ 2.568 \pm 0.21 $& $ 1.027 \pm 0.00 $& $ 1.055 \pm 0.00 $\\ \noalign{\smallskip} \hline \noalign{\smallskip}  
 pen2-F & $ 5.062 \pm 0.37 $& $ 4.462 \pm 0.30 $& $ 1.050 \pm 0.00 $& $ 1.059 \pm 0.01 $\\   
pen5-F & $ 3.595 \pm 0.25 $& $ 3.458 \pm 0.28 $& $ 1.034 \pm 0.00 $& $ 1.052 \pm 0.00 $\\   
pen10-F & $ 3.445 \pm 0.22 $& $ 2.744 \pm 0.21 $& $ 1.031 \pm 0.00 $& $ 1.053 \pm 0.00 $\\   
pen20-F & $ 3.120 \pm 0.17 $& $ 2.670 \pm 0.21 $& $ 1.032 \pm 0.00 $& $ 1.053 \pm 0.00 $\\   
penLoo & $ 3.063 \pm 0.17 $& $ 2.571 \pm 0.21 $& $ 1.027 \pm 0.00 $& $ 1.055 \pm 0.00 $\\ \noalign{\smallskip} \hline \noalign{\smallskip}  
 pen2-F+ & $ 3.723 \pm 0.29 $& $ 3.777 \pm 0.26 $& $ 1.087 \pm 0.00 $& $ 1.082 \pm 0.01 $\\   
pen5-F+ & $ 2.790 \pm 0.18 $& $ 2.698 \pm 0.19 $& $ 1.083 \pm 0.00 $& $ 1.061 \pm 0.01 $\\   
pen10-F+ & $ 2.653 \pm 0.14 $& $ 2.364 \pm 0.20 $& $ 1.085 \pm 0.00 $& $ 1.047 \pm 0.01 $\\   
pen20-F+ & $ 2.497 \pm 0.13 $& $ 2.318 \pm 0.20 $& $ 1.082 \pm 0.00 $& $ 1.049 \pm 0.01 $\\   
penLoo+ & $ 2.437 \pm 0.12 $& $ 2.218 \pm 0.18 $& $ 1.085 \pm 0.00 $& $ 1.045 \pm 0.00 $\\ \noalign{\smallskip} \hline \noalign{\smallskip}  
 
\end{tabular}
\end{center}
\end{table}

\section{Addendum to Appendix~\main{A.1}}
Whereas Lemma~\main{3} is stated for the particular case of Binomial variables, it is worth noticing that ingredients of its proof can be successfully used in order to derive non-asymptotic bounds on $\einv{\mathcal{L}(Z)}$ or $\einvz{\mathcal{L}(Z)}$ for several other distributions than the Binomial one. This has for instance be used in Sect.~6.7 of \cite{Arl:2007:phd} for the Hypergeometric and Poisson case.

\medskip

First, the lower bound in \maineqref{15} comes from Jensen's inequality:
\[ \einv{Z} \geq \Prob \paren{Z > 0} \enspace . \]

Second, taking $\theta = 0.16$ in the proof of Lemma~\main{3} gives the absolute upper bound 
\[ \einvz{Z} \leq \kappa_4 = 7.8 \] instead of the smaller value given by Lemma~4.1 of \cite{Gyo_etal:2002}.
Hence, the proof of Lemma~\main{3} only uses that $\Prob(0<Z<c_Z)=0$ for some $c_Z>0$ and that $Z$ satisfies a concentration inequality similar to Bernstein's inequality. This covers a wide class of  random variables.

Finally, notice that taking $\theta = 3 \ln(A) / A$ at the end of the proof of Lemma~\main{3}, instead of $\theta = A^{-1/2}$, leads to an upper bound \[ 1 + \kappa_5 \sqrt{\frac{\ln(A)}{A}} \geq  \sup_{np \geq A} \set{\einv{\mathcal{B}(n,p)}} \] for some numerical constant $\kappa_5$, showing that the rate $A^{-1/4}$ is far from optimal.

\section{Additional proofs}

\subsection{Proof of Lemma~\main{6}}
In this proof, we denote by $L$ any constant that may depend on $a$, $b$, $\paren{c_i}_{1 \leq i \leq 4}$, $\paren{\kappa_i}_{1 \leq i \leq 4}$, $c_{\mathrm{rich}}$ and $C$,  
possibly different from one place to another.

First of all, there is a model $m_1 \in \M_n$ such that 
\[ \ln(n)^{\kappa_1} \leq  \paren{2 a n b^{-1}}^{1/3} \leq D_{m_1} \leq \paren{2 a n b^{-1}}^{1/3} + c_{\mathrm{rich}}  \leq c_1 n (\ln(n))^{-1} \] (at least for $n \geq L$). As a consequence, \maineqref{27} implies that 
\begin{equation} \label{VFCV.le.determ.eq.crit1m1}
\crit_1(m_1) \leq a^{1/3} b^{2/3} n^{-2/3} \paren{ 3 \times 2^{-2/3}  + c_{\mathrm{rich}} \paren{\frac{b}{an}}^{1/3} } \paren{1 + c_2 \ln(n)^{-\kappa_2} } \enspace .
\end{equation}
With a similar argument, for $n \geq L$, there exists a model $m_2 \in \M_n$ such that 
\begin{equation} \label{VFCV.le.determ.eq.crit2m2}
\crit_2(m_2) \leq a^{1/3} \paren{b C}^{2/3} n^{-2/3} \paren{ 3 \times 2^{-2/3} + c_{\mathrm{rich}} \paren{\frac{bC}{a n}}^{1/3} } \paren{1 + c_2 \ln(n)^{-\kappa_2} }  \enspace .
\end{equation}

\medskip

We will now derive from \eqref{VFCV.le.determ.eq.crit2m2} some tight bounds on $D_{\mh}$. First, the upper bound in \eqref{VFCV.le.determ.eq.crit2m2} is smaller than the lower bounds in both \maineqref{29} and \maineqref{30} for $n \geq L$. This proves that 
\[ \ln(n)^{\kappa_1} \leq D_{\mh} \leq \frac{c_1 n}{\ln(n)} \enspace  . \]
Then, according to \maineqref{49}, we have for every $\mM_n$ of dimension $D_m = \paren{\frac{2a n}{bC}}^{1/3} \paren{1+\delta} $ (which is between $\ln(n)^{\kappa_1}$ and $\frac{c_1 n}{\ln(n)}$ for $n \geq L$, as long as $1 \leq \delta >-1$):
\begin{align*} 
\crit_2(m) &\geq a^{1/3} \paren{b C}^{2/3} n^{-2/3} \paren{ 2^{-2/3} (1+\delta)^{-2} + 2^{1/3} (1+\delta)  } \paren{1 - c_2 \ln(n)^{-\kappa_2} } \\
&\geq \crit_2(m_2) \times \frac{1 - c_2 \ln(n)^{-\kappa_2}}{1 + c_2 \ln(n)^{-\kappa_2}} \times \frac{f(\delta)} {3 \times 2^{-2/3} + c_{\mathrm{rich}} \paren{\frac{bC}{an}}^{1/3}}
\end{align*}
with $f$ defined by $f(\delta) = 2^{-2/3} (1+\delta)^{-2} + 2^{1/3} (1+\delta)$. Using Lemma~\ref{VFCV.le.tech.1} below, we then have
\begin{align*}
\frac{\crit_2(m)}{\crit_2(m_2)} \geq \frac{1 - c_2 \ln(n)^{-\kappa_2}}{1 + c_2 \ln(n)^{-\kappa_2}} \times \frac{3 \times 2^{-2/3} + 3 \times 2^{-14/3} \paren{\mini{\delta^2}{1}}} {3 \times 2^{-2/3} + c_{\mathrm{rich}} \paren{\frac{bC}{an}}^{1/3}} \enspace .
\end{align*}
This lower bound is strictly larger than 1 as soon as $\delta^2 \geq \ln(n)^{-\kappa_2 / 2}$ and $n \geq L$, so that 
\begin{equation} \label{VFCV.le.determ.eq.Dmh}
\paren{\frac{2an}{bC}}^{1/3} \paren{1-\ln(n)^{-\kappa_2 / 4}} \leq D_{\mh} \leq \paren{\frac{2an}{bC}}^{1/3} \paren{1+\ln(n)^{-\kappa_2 / 4}} \enspace .
\end{equation}

\medskip

We can now use \maineqref{27} in order to bound $\crit_1(\mh)$. For $n \geq L$, 
using again Lemma~\ref{VFCV.le.tech.1},
\begin{align*} 
\crit_1(\mh) &\geq a^{1/3} b^{2/3} n^{-2/3} \paren{ \paren{\frac{C}{2}}^{2/3} + \paren{\frac{C}{2}}^{-1/3} } \paren{1 - L \ln(n)^{-\kappa_2 / 4} } \\
&= a^{1/3} b^{2/3} n^{-2/3} f \paren{ C^{-1/3} - 1} \paren{1 - L \ln(n)^{-\kappa_2 / 4} } \\
&\geq a^{1/3} b^{2/3} n^{-2/3} \paren{ 3 \times 2^{-2/3} + \paren{ C^{-1/3} - 1}^2 } \paren{1 - L \ln(n)^{-\kappa_2 / 4} } \\
&\geq \crit_1 (m_1) \paren{ 1 + 2^{2/3} \times 3^{-1} \paren{ C^{-1/3} - 1}^2 - \ln(n)^{-\kappa_2 / 5} } \enspace ,
\end{align*}
which proves \maineqref{31}. \qed

\begin{remark}
A similar argument proves that for $n \geq L$,
\begin{align*}
\crit_1(\mh) &\leq \crit_1 (m_1) \paren{ 1 + 2^{2/3} \times 3^{-1} \paren{ C^{-1/3} - 1}^2 + L \ln(n)^{-\kappa_2 / 4} } \enspace .
\end{align*}
Moreover, if $\crit_1$ satisfies (ii) and (iii), we prove in a similar way that if $n \geq n_0$, for every $\mh \in \arg\min_{\mM_n} \crit_2 (m)$,
\begin{equation} \label{VFCV.le.determ.eq.res.majoration}
\crit_1(\mh) \leq \paren{1 + K(C) + \ln(n)^{-\kappa_2 / 5}} \inf_{\mM_n} \set{\crit_1(m)} \enspace .
\end{equation}
This justifies our first comment behind Thm.~\main{1}.
\end{remark}

\begin{lemma} \label{VFCV.le.tech.1}
Let $f : (-1, +\infty) \flens \R$ be defined by $f(x) = 2^{-2/3} (1+x)^{-2} + 2^{1/3} (1+x)$. Then, for every $x>-1$,
\[ f(x) \geq 3 \times 2^{-2/3} + 3 \times 2^{-14/3} \paren{\mini{x^2}{1}} \enspace . \]  
\end{lemma}
\begin{proof}[proof of Lemma~\ref{VFCV.le.tech.1}]
We apply the Taylor-Lagrange theorem to $f$ (which is infinitely differentiable) at order two, between 0 and $x$. The result follows since $f(0)=3 \times 2^{-2/3}$, $f^{\prime}(0)=0$ and $f^{\prime\prime}(t) = 6 \times 2^{-2/3} \times  (1+t)^{-4} \geq 3 \times 2^{1/3 - 4} $ if $t \leq 1$.
If $t>1$, the result follows from the fact that $f^{\prime} \geq 0$ on $[0,+\infty)$.
\end{proof}

\subsection{End of the proof of Prop.~\main{2}}
We here compute $R_{1,\widetilde{W}} (n,\phl)$ and $R_{2,\widetilde{W}}(n,\phl)$ when $V$ does not divide $n \phl$, that we have skipped in Appendix~\main{B.4.2}.

\medskip

Since $(\widetilde{W}_i)_{X_i \in \Il}$ is exchangeable and $\widetilde{W}_i$ takes only two values, 
\[ \Wl = \Es \croch{ W_i \sachant \Wl} = \frac{V}{V-1} \Prob\paren{W_i = \frac{V}{V-1} \sachant \Wl} \enspace .\]
Thus, 
\[ \loi \paren{ W_i \sachant \Wl } = \frac{V}{V-1} \mathcal{B}(\kappa^{-1} \Wl) 
\]
so that
\[ R_{2,W}(n,\phl) = \frac{1}{V-1} \qquad \mbox{and} \qquad R_{1,W}(n,\phl) = \frac{V}{V-1} \E \paren{ \widetilde{W}_{\lambda}^{-1} } - 1 \enspace .\]

There exists $a,b \in \N$ such that $0 \leq b \leq V-1$ and $n \phl = aV + b$. Then,
\[ \P \paren{ \widetilde{W}_{\lambda} = \frac{V (a(V-1) + b)}{(V-1) (aV+b)}} = \frac{V-b}{V} \quad \mbox{and} \quad \P \paren{ \widetilde{W}_{\lambda} = \frac{V (a(V-1) + b - 1)}{(V-1) (aV+b)}} = \frac{b}{V} \]
so that
\begin{align*} \E \croch{ \widetilde{W}_{\lambda}^{-1} } &= \frac{V-b}{V} \frac{(V-1) (aV+b)} {V (a(V-1) + b)} + \frac{b}{V} \frac{(V-1) (aV+b)}{V (a(V-1) + b - 1)} \\
&= 1 - \frac{b}{V (a(V-1) + b)} + \frac{(V-1)(aV+b)b}{V^2 (a(V-1) + b - 1) (a(V-1) + b)} \enspace .
\end{align*}
We deduce
\[ R_{1,\widetilde{W}} (n,\phl) = \frac{1}{V-1} - \frac{b }{(V-1) (a(V-1) + b)} + \frac{ (aV+b)b}{V (a(V-1) + b - 1) (a(V-1) + b)}  \enspace .\]
The result follows with 
\[ \delta_{n,\phl}^{(\mathrm{penV})} = \frac{b}{n\phl - a} \paren{ \frac{V-1} {V} \times \frac {n\phl} {n\phl - a - 1} - 1} \in 
\croch{ 0 ; \frac{2}{n\phl - 2} } \enspace . \qed \]

\subsection{Proof of Lemma~\main{8}}
Although this lemma can be found in \cite{Arl:2007:phd} (where it is called Lemma~5.7), we recall here its proof for the sake of completeness.

\medskip

First, split the penalty (without the constant $C$) into these two terms:
\begin{align}
\label{RP.def.ph1}
\ph_1(m) &= \sum_{\lamm} \Es \croch{ \phl \paren{  \bethlW - \bethl }^2 \sachant \Wl >0} \\
\label{RP.def.ph2}
\ph_2(m) &= \sum_{\lamm} \Es \croch{ \phlW \paren{  \bethlW - \bethl }^2 } \enspace .
\end{align}
This split into two terms is the equivalent of the split of $\penid$ into $p_1$ and $p_2$ (plus a centered term). 

We first compute this quantity, which appears in both $\ph_{1}$ and $\ph_{2}$: let $\lamm$ and $\Wl>0$,
\begin{align} \notag
& \Es \croch{\phl \carre{\bethlW - \bethl} \sachant \Wl} 
= \Es \croch{ \phl \carre{ \frac{1}{n\phl}
\sum_{X_i \in \Il} (Y_i - \betl) \left(1 - \frac{W_i}{W_{\lambda}} \right) } \sachant \Wl }  \\ \label{RP.eq.calc.ph.int1}
& \quad =  \frac{1}{n^2 \phl} \Biggl[
\sum_{X_i \in \Il} \carre{Y_i - \betl} \Es \croch{ \carre{ 1 -
\frac{W_i}{W_{\lambda}} } \sachant \Wl} 
\\ \notag & \quad + \frac{1}{n^2 \phl} \sum_{i \neq j, X_i
\in \Il, X_j \in \Il} (Y_i - \betl) (Y_j - \betl) \Es \croch{ 
\left( 1 - \frac{W_i}{W_{\lambda}} \right)  \left( 1 -
\frac{W_j}{W_{\lambda}} \right) \sachant \Wl} \Biggr] \enspace .
\end{align}
Since the weights are exchangeable, $(W_i)_{X_i \in \Il}$ is also exchangeable conditionally to $\Wl$ and $(X_i)_{1 \leq i \leq n}$. Thus, the ``variance'' term
\[ R_V(n,n\phl,\Wl,\loi(W)) \egaldef \Es \croch{ \paren{W_i - \Wl}^2 \sachant \Wl } \] does not depend from $i$ (provided that $X_i \in \Il$), and the ``covariance'' term
\[ R_C(n,n\phl,\Wl,\loi(W)) \egaldef \Es \croch{ \paren{W_i - \Wl}\paren{W_j - \Wl} \sachant \Wl } \] does not depend from $(i,j)$ (provided that $i \neq j$ and $X_i, X_j \in \Il$). Moreover, 
\begin{align*} 0 &= \Es \croch{ \paren{ \sum_{X_i \in \Il} \paren{W_i - \Wl} }^2 \sachant \Wl} \\
&= n \phl R_V(n,n\phl,\Wl,\loi(W)) + n \phl \paren{ n \phl - 1} R_C(n,n\phl,\Wl,\loi(W))
\end{align*}
so that, if $n \phl \geq 2$, 
\begin{equation} \label{RP.eq.calc.ph.int2} R_C(n,n\phl,\Wl,W) = \frac{-1}{n\phl -1} R_V(n,n\phl,\Wl,\loi(W)) \qquad \mbox{and} \qquad R_V(n,1,\Wl,\loi(W)) = 0 \enspace . \end{equation}

Combining \eqref{RP.eq.calc.ph.int1} and \eqref{RP.eq.calc.ph.int2}, we obtain 
\begin{align} \label{RP.eq.calc.ph.int3}
\Es \croch{\phl \carre{\bethlW - \bethl} \sachant \Wl} 
&= \frac{R_V(n,n\phl,\Wl,\loi(W)) }{\Wl n^2 \phl} \1_{n \phl \geq 2}  \\ \notag
&\quad \times \left[ \frac{n\phl}{n\phl - 1} S_{\lambda,2} - \frac{1}{n\phl -
1} S_{\lambda,1}^2 \right]
\end{align}

Combining \eqref{RP.eq.calc.ph.int3} and \eqref{RP.def.ph1} (resp.  \eqref{RP.eq.calc.ph.int3} and \eqref{RP.def.ph2}), we have the following expressions for $\ph_1$ and $\ph_2$:
\begin{align}
\label{RP.eq.calc.his.ph1} 
\ph_1(m) &= \sum_{\lamm} \frac{R_{1,W}(n,\phl) \1_{n \phl \geq
2}}{n^2 \phl}  \croch{ \frac{n\phl}{n\phl - 1} S_{\lambda,2} - \frac{1}{n\phl -
1} S_{\lambda,1}^2 }  \\
\label{RP.eq.calc.his.ph2}
\ph_2(m) &= \sum_{\lamm} \frac{R_{2,W}(n,\phl) \1_{n \phl \geq
2}}{n^2 \phl}  \croch{ \frac{n\phl}{n\phl - 1} S_{\lambda,2} - \frac{1}{n\phl -
1} S_{\lambda,1}^2 } \enspace .
\end{align}
Remark that the terms of the sum for which $n\phl=1$ are all equal to zero, which can be ensured with the convention $0 \times \infty = 0$ since $R_{1,W}(n,n^{-1})=R_{2,W}(n,n^{-1})=0$. The result follows. \qed

\subsection{Concentration of $\punmin$: detailed proof}
Within the proof of Prop.~\main{9}, we used Lemma~\main{4} in order to control the deviations of $\El\croch{\punmin(m)}$ around its expectation. Implicitly, we used the following lemma (which is indeed a straightforward consequence of Lemma~\main{4}).

\begin{lemma} \label{VFCV.le.conc.Ep1-Emp1}
We assume that $\min_{\lamm} \set{ n \pl } \geq B_n \geq 1$.
\begin{enumerate}
\item Lower deviations: let $c_1=0.184$. For all $x \geq 0$, with probability at least $1 - e^{-x}$, 
\begin{gather}
\label{RP.eq:conc_Ep1-Emp1:inf}
\El\croch{\punmin(m)} \geq \E\croch{\punmin(m)} - \theta^-(x,B_n,D_m,A,\sigmin) \times \E \croch{p_2(m)} \\
\mbox{with} \qquad \theta^- \egaldef L \croch{ \varphi_1(c_1 B_n) + \frac{A^2}{\sigmin^2} \sqrt{e^{- c_1 B_n} + \frac{x}{D_m} } } \notag \end{gather}
\item Upper deviations: let $c_2 = 0.28$ and $c_4 = 0.09$. For every $x \geq 0$, with probability at least $1 - e^{-x}$,
\begin{gather}
\label{RP.eq:conc_Ep1-Emp1:sup}
\El\croch{\punmin(m)} \leq  \E\croch{\punmin(m)} + \theta^+(x,B_n,D_m,A,\sigmin) \E \croch{p_2(m)} \\
\mbox{with} \qquad \theta^+ \egaldef L \croch{ \varphi_1 \left(c_2 B_n  \right) + \frac{A^2}{\sigmin^2} \sqrt{ x D_m^{-1} + e^{-c_4 B_n} } \paren{ \maxi{1}{\sqrt{x + D_m e^{-c_4 B_n }} }} }  \notag \enspace .
\end{gather}
\end{enumerate}
\end{lemma}
\begin{proof}
From \maineqref{19} and \maineqref{37}, we have an explicit expression for $\punmin$. We then apply Lemma~\main{4}, with $X_{\lambda} = n\phl$ and $a_{\lambda} = \pl \carre{\sigl} \geq 0$. 
For $\theta^+$, we used the general upper bound
\[ \max_{\lamm} \paren{\sigl}^4 \paren{ \sum_{\lamm} \sigl^4 }^{-1} \leq 1 \enspace .\]
\end{proof}

\begin{remark} 
If $B_n \geq \paren{\maxi{c_1^{-1}}{c_4^{-1}}} \ln(n)$, for every $\gamma>0$, 
\[ \maxi{\theta^-}{\theta^+} \paren{ \gamma \ln(n) , B_n, D_m, A, \sigmin } \leq L_{\gamma} A^2 \sigmin^{-2} D_m^{-1/2} \ln(n) \] since $D_m \leq n$.
\end{remark}

\bibliographystyle{alpha}
\bibliography{penVF,publi_arlot_en}